\documentclass[11pt,reqno]{amsart}

\usepackage{amssymb}
\usepackage{hyperref} 
\usepackage{mathrsfs}
\usepackage{amsmath}
\usepackage{amsthm}
\usepackage{enumerate}
\usepackage{dsfont}
\usepackage{color}
\usepackage[a4paper]{geometry}
\usepackage[utf8]{inputenc}
\usepackage{stmaryrd}
\usepackage{titlesec}
\usepackage{mathabx}
\usepackage{appendix}

\usepackage[all]{xy}
\let\oldtocsection=\tocsection
\let\oldtocsubsection=\tocsubsection
\let\oldtocsubsubsection=\tocsubsubsection

\renewcommand{\tocsection}[2]{\hspace{0em}\oldtocsection{#1}{#2}}
\renewcommand{\tocsubsection}[2]{\hspace{1em}\oldtocsubsection{#1}{#2}}
\renewcommand{\tocsubsubsection}[2]{\hspace{2em}\oldtocsubsubsection{#1}{#2}}

\numberwithin{theorem}{section}
\numberwithin{equation}{section}

\renewcommand{\epsilon}{\varepsilon}

\let\div=\relax
\DeclareMathOperator{\div}{div}

\newcommand\R{\mathbb{R}}
\newcommand\E{\mathbb{E}}
\newcommand\calS{\mathcal{S}}
\newcommand\G{\mathcal{G}}

\newcommand\LL{\mathcal L}

\newcommand\SSS{{\mathcal S}_o'}

\newcommand{\mcO}{\mathcal{O}}

\newcommand{\RR}{\mathbb{R}}  
\newcommand{\EE}{\mathbb{E}}
\newcommand{\mcE}{\mathcal{E}}

\newcommand{\mcB}{\mathcal{B}} 
\newcommand{\mcS}{\mathcal{S}}   
\newcommand{\mcA}{\mathcal{A}}
\newcommand{\mcP}{\mathcal{P}}
\newcommand{\mcQ}{\mathcal{Q}}

\newcommand{\mcD}{\mathcal{D}}
\newcommand{\mcM}{\mathcal{M}}
\newcommand{\mcH}{\mathcal{H}}
\newcommand{\mcI}{\mathcal{I}}

\newcommand{\calM}{\ensuremath{\mathcal{M}}}
\newcommand{\calC}{\ensuremath{\mathcal{C}}}

\newcommand{\rest}[1]{\llparenthesis\, #1 \,\rrparenthesis}
\newcommand{\restbis}[1]{\llparenthesis\, #1 \,\rrparenthesis^\sharp}

\newcommand{\ssk}{\smallskip}

\newenvironment{Dem}[1][\unskip]{%
    \begin{list}{\hspace{0.5cm}{\sf \textbf{Proof #1 --}}}{%
        \setlength{\topsep}{0pt}%
        \setlength{\leftmargin}{0pt}%
        \setlength{\rightmargin}{0pt}%
        \setlength{\listparindent}{0pt}%
        \setlength{\itemindent}{0pt}%
        \setlength{\parsep}{0pt}%
        \addtolength{\leftmargin}{20pt}%
        \addtolength{\rightmargin}{0pt}%
        } \item }{\hfill{\begin{flushright} $\rhd$ \end{flushright}}\end{list}\smallskip}

\titleformat{\section}[block]
  {\normalfont\sffamily}{\thesection}{.5em}{\titlerule\\[.8ex]\bfseries} 

\titleformat{\subsection}[wrap]{\sffamily\fontseries{b}\selectfont\filright}{\thesubsection.}{.5em}{}
                  \titlespacing{\subsection}{20pc}{1.5ex plus .1ex minus .2ex}{1pc}

\titleformat{\subsubsection}[wrap]{\sffamily\fontseries{b}\selectfont\filright}{\thesubsubsection.}{.5em}{}
                  \titlespacing{\subsubsection}{20pc}{1ex plus .1ex minus .2ex}{1pc}

\def\XXint#1#2#3{{\setbox0=\hbox{$#1{#2#3}{\int}$}
     \vcenter{\hbox{$#2#3$}}\kern-.5\wd0}}

\numberwithin{subsection}{section}
\numberwithin{subsubsection}{subsection}

\newtheoremstyle{mystyle}
{3pt}
{3pt}
{\sffamily}
{0em}
{\bfseries\sffamily}
{.}
{.5em}
{}
\theoremstyle{mystyle}
\newtheorem{thm}{Theorem}

\newtheorem{cor}[thm]{\hspace{-0.15cm}  {Corollary} }
\newtheorem{lem}[thm]{\hspace{-0.15cm}  {Lemma} }
\newtheorem{prop}[thm]{\hspace{-0.15cm} {Proposition}}
\newtheorem{defn}[thm]{ \hspace{-0.3cm} {Definition}}
\newtheorem*{defn*} {Definition}
\newtheorem*{lem*} {Lemma}
\newtheorem{rem}[thm]{\hspace{-0.15cm} {Remark}}

\numberwithin{equation}{section} 

\newcommand\restr[2]{{
  \left.\kern-\nulldelimiterspace 
  #1 
  \vphantom{\big|} 
  \right|_{#2} 
  }}

 \author[I. Bailleul]{Ismael Bailleul}
\address{Institut de Recherche Mathematiques de Rennes, 263 Avenue du General Leclerc, 35042 Rennes, France}
 \email{ismael.bailleul@univ-rennes1.fr}

\author[F. Bernicot]{Fr\'ed\'eric Bernicot}
\address{CNRS - Universit\'e de Nantes \\ Laboratoire de Math\'ematiques Jean Leray \\ 2, Rue de la Houssini\`ere 44322 Nantes Cedex 03, France}
\email{frederic.bernicot@univ-nantes.fr}
 
\author{Dorothee Frey}
\address{CNRS - Universit\'e Paris-Sud, Laboratoire de Math\'ematiques, UMR 8628, 91405  Orsay, France }
\curraddr{Delft Institute of Applied Mathematics, Delft University of Technology, P.O. Box 5031, 2600 GA Delft, The Netherlands}
\email{d.frey@tudelft.nl}

\date{\today}

\begin{document}
\date{\today}

\vspace*{3ex minus 1ex}
\begin{center}
\huge\sffamily{Space-time paraproducts for paracontrolled calculus, 3d-PAM and multiplicative Burgers equations \vspace{0.5cm}}
\end{center}
\vskip 5ex minus 1ex

\begin{center}
{\sf I. BAILLEUL\footnote{I.Bailleul thanks the U.B.O. for their hospitality, part of this work was written there.}, F. BERNICOT\footnote{F. Bernicot's research is partly supported by the ANR projects AFoMEN no. 2011-JS01-001-01 and HAB no. ANR-12-BS01-0013 and the ERC Project FAnFArE no. 637510.} and D. FREY\footnote{D. Frey's research is supported by the ANR project HAB no. ANR-12-BS01-0013. \vspace{0.1cm}

{\sf AMSClassification:} 60H15; 35R60; 35R01 \vspace{0.1cm}

}}
\end{center}

\vspace{1cm}

\begin{center}
\begin{minipage}{0.8\textwidth}
\renewcommand\baselinestretch{0.7} \rmfamily {\scriptsize {\bf \sc \noindent Abstract.} 
We sharpen in this work the tools of paracontrolled calculus in order to provide a complete analysis of the parabolic Anderson model equation and Burgers system with multiplicative noise, in a $3$-dimensional Riemannian setting, in either bounded or unbounded domains. With that aim in mind, we introduce a pair of intertwined space-time paraproducts on parabolic H\"older spaces, with good continuity, that happens to be pivotal and provides one of the building blocks of higher order paracontrolled calculus.  }
\end{minipage}
\end{center}

\bigskip

\setcounter{tocdepth}{2}
\begin{quote}
\footnotesize  \tableofcontents
\end{quote}

\section[\hspace{0.7cm} Introduction]{Introduction}
\label{SectionIntroduction}

It is probably understated to say that the work \cite{H} of Hairer has opened a new era in the study of stochastic singular parabolic partial differential equations. It provides a setting where one can make sense of a product of a distribution with parabolic non-positive H\"older regularity index, say $a$, with a function with non-negative regularity index, say $b$, even in the case where $a+b$ is non-positive, and where one can make sense of, and solve, a large class of parabolic stochastic singular partial differential equations by fixed point methods. The parabolic Anderson model equation (PAM) 
\begin{equation}
\label{EqPAMEq}
(\partial_t+L) u = u\zeta,
\end{equation}
studied in Section \ref{SectionPAM} in a $3$-dimensional unbounded background, is an example of such an equation, as it makes sense in that setting to work with a distribution $\zeta$ of H\"older exponent $\alpha-2$, for some $\alpha<\frac{1}{2}$, while one expects the solution $u$ to the equation to be of parabolic H\"older regularity $\alpha$, making the product $u\zeta$ ill-defined since $\alpha+(\alpha-2)\leq 0$.

\medskip

The way out of this quandary found by Hairer has its roots in Lyons' theory of rough paths, which already faced the same problem. Lyons' theory addresses the question of making sense of, and solving, controlled differential equations 
\begin{equation}
\label{EqControlledRDE}
dz_t = V_i(z_t)\,dX^i_t
\end{equation}
in $\R^d$ say, driven by an $\R^\ell$-valued $\frac{1}{p}$-H\"older control $X = \big(X^1,\dots,X^\ell\big)$, with $p\geq 2$, and where $V_i$ are sufficiently regular vector fields on $\R^d$. Typical realizations of a Brownian path are $\frac{1}{p}$-H\"older continuous, with $p>2$, for instance. One expects a solution path to equation \eqref{EqControlledRDE} to be $\frac{1}{p}$-H\"older continuous as well, in which case the product $ V_i(z_t)\,dX^i_t$, or the integral $\int_0^t  V_i(z_s)\,dX^i_s$, cannot be given an intrinsic meaning since $\frac{1}{p}+\big(\frac{1}{p}-1\big)\leq 0$. Lyons' deep insight was to realize that one can make sense of, and solve, equation \eqref{EqControlledRDE} if one assumes one is given an enriched version of the driving signal $X$ that formally consists of $X$ together with its non-existing iterated integrals. The theory of regularity structures rests on the same philosophy, and the idea that the enriched noise should be used to give a local description of the unknown $u$, in the same way as polynomials are used to define and describe locally $C^k$ functions.

\ssk

At the very same time that Hairer built his theory, Gubinelli, Imkeller and Perkowski proposed in \cite{GIP} another implementation of that philosophy building on a different notion of local description of a distribution, using paraproducts on the torus. The machinery of paracontrolled distributions introduced in \cite{GIP} rests on a first order Taylor expansion of a distribution that happened to be sufficient to deal with the stochastic parabolic Anderson equation \eqref{EqPAMEq} on the $2$-dimensional torus, the stochastic additive Burgers equation in one space dimension \cite{GIP}, the $\Phi^4_3$ equation on the $3$-dimensional torus \cite{CC, ZZ3} and the stochastic Navier-Stokes equation with additive noise \cite{ZZ1,ZZ2}. The KPZ equation can also be dealt with using this setting \cite{GP}. Following Bony's approach \cite{Bony}, the paraproduct used in \cite{GIP} is defined in terms of Fourier analysis and does not allow for the treatment of equations outside the flat background of the torus or the Euclidean space, if one is ready to work with weighted functional spaces. The geometric restriction on the background was greatly relaxed in our previous work \cite{BB} by building paraproducts from the heat semigroup associated with the operator $L$ in the semilinear equation. A theory of paracontrolled distributions can then be considered in doubling metric measure spaces where one has small time Gaussian estimates on the heat kernel and its 'gradient' -- see \cite{BB}. This setting already offers situations where the theory of regularity structures is not known to be working. The stochastic parabolic Anderson model equation in a $2$-dimensional doubling manifold was considered in \cite{BB} as an example. The first order 'Taylor expansion' approach of paracontrolled calculus seems however to restrict a priori its range of application, compared to the theory of regularity structures, and it seems clear that a kind of higher order paracontrolled calculus is needed to extend its scope. We tackle in the present work the first difficulty that shows off in this program, which is related to the crucial use of commutator estimates between the heat operator and a paraproduct, which is one of the three workhorses of the paracontrolled calculus method, together with Schauder estimates and another continuity result on some commutator. The development of a high order paracontrolled calculus is the object of another work \cite{BB16}.

\ssk

Working in unbounded spaces with weighted functional spaces requires a careful treatment which was not done so far. We shall illustrate the use of our machinery on two examples: The parabolic Anderson model (PAM) equation \eqref{EqPAMEq}  in a possibly unbounded $3$-dimensional Riemannian manifold, and Burgers equation with multiplicative noise in the $3$-dimensional closed Riemannian manifold. Hairer and Labb\'e have very recently studied the (PAM) equation  in $\R^3$ from the point of view of regularity structures \cite{HLR3} -- see also the work \cite{HP} of Hairer and Pardoux. They had to introduce some weights $\varpi$ to get a control on the growth of quantities of interest at spatial infinity. A non-trivial part of their work consists into tracking the time-behavior of their estimates, with respect to the time, which requires the use of time-dependent weights. For the same reason, we also need to use weighted spaces and working with the weights of \cite{HP,HLR3} happens to be convenient. Our treatment is however substantially easier, as we do not need to 
travel backwards in time such as required in the analysis of the reconstruction operator in the theory of regularity structures. As a matter of fact, our results on the (PAM) equation give an alternative approach, and  provide a non-trivial extension, of the results of \cite{HLR3} to a non-flat setting, with a possibly wider range of operators $L$ than can be treated presently in the theory of regularity structures. As for Burgers equation with multiplicative noise, it provides a description of the random evolution of a velocity field on the $3$-dimensional torus, subject to a random rough multiplicative forcing, and whose dynamics reads 
\begin{equation}
\label{EqBurgersEq}
(\partial_t+L)u + (u\cdot\nabla)u = \textrm{M}_\zeta u,
\end{equation}
where $\zeta$ is a $3$-dimensional white noise with independent coordinates, and 
$$
\textrm{M}_\zeta u := \big(\zeta^1 u^1, \zeta^2 u^2, \zeta^3 u^3\big),
$$ 
for the velocity field $u = \big(u^1, u^2, u^3\big) : M^3 \to \RR^3$. With zero noise $\zeta$, this $3$-dimensional Burgers system plays a very important role in the theory of PDEs coming from fluid mechanics, and later from condensed matter physics and statistical physics. It has been proposed by Burgers in the 30's as a simplified model of dynamics for Navier-Stokes equations. A change of variables, called after Cole and Hopf, can be used to reduce the deterministic quasilinear parabolic equation to the heat equation, thus allowing the derivation of exact solutions in closed form. Despite this fact, the study of Burgers system is still very fashionable as a benchmark model that can be used  to understand the basic features of the interaction between nonlinearity and dissipation. Motivated by the will to turn Burgers equation into a model for turbulence, stochastic variants have been the topic of numerous recent works \cite{BCJL, HV,HW, H, GIP, GP}, where a random forcing term is added in the equation, mainly in one space dimension, with an additive space-time white noise -- that is with a space-time white noise instead of $\textrm{M}_\zeta u$ with $\zeta$ space white noise. The Cole-Hopf transformation can formally be used again, and turns a solution to the $1$-dimensional stochastic Burgers equation with additive space-time noise to the heat equation with multiplicative space-time noise, with a very singular noise, such as detailed in \cite{GP}. A similar change of variable trick can be used for the study of the above multidimensional stochastic Burgers system with multiplicative noise \eqref{EqBurgersEq}; we shall analyse it in Section \ref{SubsectionDeterministicBurger}. Also, one can consider the study of this example as a first step to understanding the dynamics of the $3$-dimensional stochastic incompressible Navier-Stokes equation, with multiplicative noise, where the incompressibility brings the additional difficulty to deal with the Leray projector to keep the vanishing divergence property. In any case, equation \eqref{EqBurgersEq} seems not to have been studied so far, to the best of our knowledge.

\ssk

Contrary to the theory of regularity structures, whose introduction requires to set up a whole new algebraic-analytic setting, the analytic part of paracontrolled calculus is based only on classical ingredients, and its use in solving some singular stochastic partial differential equation involves an elementary reasoning. This machinery is described in simple terms in Section \ref{SectionParacontrolledNutshell}, which serves as a baseline for the study of the parabolic Anderson and Burgers equations in Section \ref{SectionPAM}. 

\medskip

The geometric and functional settings in which we lay down our study are described in Section \ref{SectionSetting}. In short, we work on a doubling metric measure manifold $(M,d,\mu)$, equipped with a Riemannian operator
$$
L = -\sum_{i= 1}^{\ell_0} V_i^2
$$ 
given by the finite sum of squares of vector fields. The heat semigroup of the operator $L$ is assumed to have a kernel that satisfies Gaussian pointwise bounds, together with its iterated derivatives; precise conditions are given in the beginning of Section \ref{SubsectionSubLaplacian}. Such a setting covers a number of interesting cases. One can use the semigroup to construct in an intrinsic way the scale of spatial H\"older spaces $C^\alpha(M)$ on $M$ and a scale of parabolic H\"older spaces $\calC^\alpha\big([0,T] \times M\big)$ in which the (PAM) and Burgers equations will eventually be solved. Some Schauder-type regularity estimates for the heat semigroup, proved in Section \ref{SubsectionSchauder}, will be instrumental for that purpose. We call resolution map of the heat semigroup the map that associates to a distribution $f$ the solution to the equation $(\partial_t+L) v = f$, with zero initial condition. One of our main contributions is the introduction of a pair of paraproducts built from the heat semigroup, intertwined via the resolution map, that are used to get exact formulas where formulas with a remainder were used previously \cite{GIP,GP,BB}. These two paraproducts share the same algebraic structure and the same analytic properties, most importantly a cancellation property that we introduce in Section \ref{SubsectionCancellations}. It allows in particular to set the stage in a more natural function space than previously done. They consist in some sense of space-time paraproducts in the parabolic variable, whose continuity properties together with Schauder estimates allow to obtain some crucial estimates in $L^\infty_T C^\alpha(M)$ spaces.

\ssk

The technical core of the paracontrolled calculus, such as defined by Gubinelli, Imkeller and Perkowski, is a continuity estimate for a corrector that allows to make sense of an a priori undefined term by compensating it by another potentially undefined term with a simpler structure, and to separate analytic from probabilistic considerations. We prove in Section \ref{SubsectionIteratedCommutators} that this result holds in our general setting as well. As a result, we are able to prove the following kind of results on the (PAM) in a $3$-dimensional possibly unbounded measured manifold $(M,d,\mu)$ that is Ahlfors regular, and working with a second order differential operator $L$ that satisfies some mild assumptions stated in Section \ref{SubsectionSubLaplacian}. We also study the multiplicative Burgers equations in a bounded ambiant space. In statements below, $\xi$ stands for a space white noise on $(M,\mu)$, and $\xi^\epsilon := \big(e^{-\epsilon L}\big)\xi$ stands for its regularization via the heat semigroup. Full details on the mathematical objects involved in the statements will be given along the way. The notion of solution to the (PAM) equation \eqref{EqPAMEq} depends on a notion of (PAM)-enhancement $\widehat\zeta$ of a distribution $\zeta\in C^{\alpha-2}(M)$. To every such enhancement of $\zeta$ is associated a Banach space $\mcD\big(\widehat\zeta\,\big)$ of distributions within which one can make sense of the equation and look for the solution to it -- this is the space of paracontrolled distributions; see Sections \ref{SectionParacontrolledNutshell} and \ref{SectionPAM}. We refer to Section \ref{SubsectionParabolicHolderSpaces} for the definition of the weighted spaces used below, and to Section \ref{SubsectionSubLaplacian} for the statement of Assumptions {\bf (A)} on the heat semigroup generated by $L$. Assumption {\bf (B)} is a statement about the, probabilistic, renormalisation process needed to make sense of some ill-defined terms; we take it for granted in the present work. It is fully spelled out in Section \ref{SectionPAM}, and hints about the problems involved in this operation are given in Section \ref{SectionRenormalisation}.

\ssk 

\begin{thm}   \label{thm:pam}   {\sf 
Let $(M,d,\mu)$ be a $3$-dimensional possibly unbounded metric measure manifold. Assume the heat semigroup satisfies Assumptions {\bf (A)} and that the vector fields $V_i$ are divergence-free. Let further work under assumption {\bf (B)}. Given $\alpha \in \big(\frac{2}{5},\frac{1}{2}\big)$, and a (PAM)-enhancement of a distribution $\zeta\in C^{\alpha-2}$, the parabolic Anderson model equation on $M$ has a unique paracontrolled solution in $\mcD\big(\widehat\zeta\,\big)$. Moreover, the space white noise $\xi$ has a natural (PAM)-enhancement, and there exists a sequence $\big(\lambda^\epsilon\big)_{0<\epsilon\leq 1}$  of {\it deterministic and time-independent functions} such that for every finite positive time horizon $T$ and every initial data $u_0 \in C^{4\alpha}_{w_0}(M)$, the solution $u^\epsilon$ of the renormalized equation 
   \begin{equation*} 
   \partial_t u^\epsilon + L u^\epsilon = u^\epsilon\,\big(\xi^\epsilon - \lambda^\epsilon\big),\qquad u^\epsilon(0)=u_0
   \end{equation*}
converges in probability to the solution $u\in \calC^\alpha_w\big([0,T]\times M\big)$ of the parabolic Anderson model equation on $M$ associated with the natural enhancement of $\xi$. The result holds with $w=1$ and $T=\infty$, if $\mu(M)$ is finite.  }
\end{thm}

\ssk 

Let emphasize that uniqueness has to be understood as uniqueness of a solution in a suitable class of paracontrolled distributions, in which the problem is formulated. Note also that we use weighted spatial and parabolic H\"older spaces to deal with the unbounded nature of the ambient space $M$. In $\RR^3$, one can typically work with the weights $w(x,\tau) = e^{\tau(1+|x|)}$ and $w_0(x) = w(x,0)$ a constant -- these weights were already used by Hairer and Labb\'e in \cite{HLR3}; see section \ref{SubsectionParabolicHolderSpaces}. Hairer and Labb\'e \cite{HLR3} are able to work in the range $-\frac{1}{2}<\alpha\leq 0$, in the setting of regularity structures; we do not know how to deal with such a situation in our setting. Note on the other hand that we described in the appendix of \cite{BB} how to extend the paracontrolled calculus to a Sobolev setting. Together with the present work, this allows to solve the (PAM) equation in Sobolev spaces $W^{\alpha,p}$ for a large enough finite positive exponent $p$. The above H\"older setting corresponds to working with $p=\infty$. The robustness of our framework in terms of the operator $L$ or the ambient geometry is useful, at least insofar as the tools of regularity structures have not been adapted so far in a non-flat setting.  Moreover, as explained before, it is easier to deal with the time-dependent weight through the current paracontrolled approach than via the regularity structures theory, as done in \cite{HLR3}.

\medskip

As we shall see, the computations needed to handle the (PAM) equation and multiplicative Burgers system involve almost the same quantities. As far as the latter is concerned, we can prove the following result, under the same conditions on the operator $L$ as above, in a \textit{bounded} geometry. We identify in the renormalized equation \eqref{eq:renor} below a symmetric matrix $d$ with its associated quadratic form. We can work in such a bounded domain with the weight $w : \; (x,\tau) \mapsto e^{\kappa \tau}$, for a large enough positive constant $\kappa$. Note here that the above mentioned notion of enhancement $\widehat\zeta$ of a distribution $\zeta\in C^{\alpha-2}$ depends on the equation under study, which is why we called it (PAM)-enhancement above. We denote by {\bf (B')} an assumption similar to assumption {\bf (B)}, about the renormalisation process of a number of ill-defined terms that appear in the analysis of the Burgers system. We work here in the $3$-dimensional torus with $\RR^3$-valued fields $u$.

\ssk

\begin{thm}   \label{thm:burgers}   {\sf 
Assume the heat semigroup satisfies Assumptions {\bf (A)} and that the vector fields $V_i$ are divergence-free. Let further work under assumption {\bf (B')}. Given $\alpha \in \big(\frac{2}{5},\frac{1}{2}\big)$, and a Burgers-enhancement of $\zeta\in C^{\alpha-2}$, the multiplicative Burgers equation \eqref{EqBurgersEq} on $M$ has a unique  local in time paracontrolled solution in $\mcD\big(\widehat\zeta\,\big)$. Moreover, the space white noise $\xi$ has a natural Burgers-enhancement, and there exists sequences of {\it time-independent and deterministic} ${\mathbb R}^3$-valued functions $\big(\lambda^\epsilon\big)_{0<\epsilon\leq 1}$ and $(3\times 3)$-symmetric-matrix-valued functions $\big(d^\epsilon\big)_{0<\epsilon\leq 1}$ on $M$, such that if one denotes by $u^\epsilon$ the solution of the renormalized equation
   \begin{equation} 
   \label{eq:renor}
   \partial_t u^\epsilon + L u^\epsilon + \big(u^\epsilon\cdot\nabla\big) u^\epsilon = \textrm{M}_{\xi^\epsilon-\lambda^\epsilon} u^\epsilon - d^\epsilon\big(u^\epsilon,u^\epsilon\big) \qquad u^\epsilon(0)=u_0  
   \end{equation}
with initial condition $u_0 \in \calC^{4\alpha}$, then $u^\epsilon$ converges in probability to the solution $u\in \calC^\alpha$ of the multiplicative Burgers equation, locally in time.    }
\end{thm}

\ssk

Details on Theorems \ref{thm:pam} and \ref{thm:burgers} can be found in Section \ref{SectionPAM}. These statements are two-sided, with the well-posedness of the paracontrolled version of the equations on the one hand, and the link between this notion of solution and the convergence of solutions to a renormalized regularized version of the initial equation on the other hand. Assumptions {\bf (B)} and \textbf{(B')} deal with the latter side of the stody. Note that after the very recent works \cite{BHZ, CH} of Bruned-Hairer-Zambotti and Chandra-Hairer on renormalisation within the regularity structure approach to singular PDEs, there is no doubt anymore that this probabilistic step should be doable in a paracontrolled setting as well, in some generality.

\bigskip

\noindent {\sf \textbf{Notations.}} Let us fix here some notations that will be used throughout the work. 
\begin{itemize}
   \item Given a metric measure space $(M,d,\mu)$, we shall denote its parabolic version by $(\mcM,\rho,\nu)$, wherethe parabolic space 
   $$
   \mcM:= M \times {\mathbb R}
   $$ 
   is equipped with the parabolic metric
$$ 
\rho\big((x,\tau), (y,\sigma)\big) = d(x,y) + \sqrt{|\tau-\sigma|}
$$
and the parabolic measure $\nu=\mu \otimes dt$. Note that for  $(x,\tau) \in \mcM$ and small radius $r>0$, the parabolic ball $B_\mcM\big((x,\tau),r\big)$
has volume 
$$ 
\nu \Big( B_\mcM \big( (x,\tau),r \big) \Big) \approx r^2\,\mu\big( B(x,r) \big).  \vspace{0.15cm}
$$
We shall denote by $e$ a generic element of the parabolic space $\mcM$.   \vspace{0.15cm}
   
   \item Given an unbounded linear operator $L$ on $L^2(M)$, we denote by $\mcD_2(L)$ its domain. We give here the definition of a distribution, as it is understood in the present work. The definition will always be associated with the operator $L$ described in Subsection \ref{SubsectionSubLaplacian} below. 

\ssk

Fix a point $o\in M$ and then define a Fr\'echet space $\mcS_o$ of test functions $f$ on $\mcM$ requiring that 
$$ 
N_a(f) := \Big\|\big(1+|\tau|\big)^{a_1}\big(1+d(o,\cdot)\big)^{a_2} \partial_\tau^{a_3} (L^*)^{a_4} f\Big\|_{2,d\nu} <\infty,
$$
for all tuples $a=(a_1,\dots,a_4)$ of integers; we equip $\mcS_o$ of the Fr\'echet space structure associated with the following family of semi-norms $N_a$. A \textsf{\textbf{distribution}} is a continuous linear functional on $\mcS_o$; we write $\mcS_o'$ for the set of all distributions.   \vspace{0.15cm}

 \item Spatial H\"older spaces $C^\gamma$ and parabolic Space-time H\"older spaces $\calC^\gamma$ will be rigorously defined in Section \ref{SubsectionParabolicHolderSpaces}, and the weights $\varpi$ and $p_a$ will be introduced in Section \ref{SubsectionSchauder}. To deal with remainder terms in some paracontrolled expansions, we shall use the following notation. For $\gamma\in {\mathbb R}$ and $c$ a non negative integer, we shall denote by $\rest{\gamma}_c$ an element of $\calC^\gamma_{\varpi p_{ca}}$, and by $\restbis{\gamma}_c$ an element of $L^\infty_T C^{\gamma}_{\varpi p_{ca}}$.   \vspace{0.15cm}
     
   \item As a last bit of notation, we shall always denote by $K_Q$ the kernel of an operator $Q$, and write $\lesssim_T$ for an inequality that holds up to a positive multiplicative constant that depends only on $T$. 
\end{itemize}
\bigskip

\section[\hspace{0.7cm} Paracontrolled calculus in a nutshell]{Paracontrolled calculus in a nutshell}
\label{SectionParacontrolledNutshell}

The theories of regularity structures and paracontrolled calculus aim at giving a framework for the study of a class of classically ill-posed stochastic parabolic partial differential equations (PDEs), insofar as they involve illicit operations on the objects at hand. This is typically the case in the above parabolic Anderson model and Burgers equations, where the products $u\zeta$ and $\textrm{M}_\zeta u$ are a priori meaningless, given the expected regularity properties of the solutions to the equations. So a regularization of the noise does not give a family of solutions to a regularized problem that converge in any reasonable functional space to a limit that could be defined as a solution to the original equation. To bypass this obstacle, both the theory of regularity structures and paracontrolled calculus adopt a point of view similar to the point of view of rough paths analysis, according to which a good notion of solution requires the enhancement of the notion of noise into a finite collection of objects/distributions, built by purely probabilistic means, and that a solution to the equation should locally be entirely described in terms of these objects. This collection of reference objects depends on the equation under study, and plays in the setting of regularity structures the role played by polynomials in the world of $C^k$ maps, where they provide local descriptions of a function in the form of a Taylor expansion. Something similar holds in paracontrolled calculus. In both approaches, the use of an ansatz for the solution space allows to make sense of the equation and get its well-posed character by deterministic fixed point methods, and provides as a consequence solutions that depend continuously on all the parameters in the problem.

To be more concrete, let us take as an introduction to these theories the example of the $2$-dimensional (PAM) equation, fully studied in \cite{H,GIP,HL,BB}. The space white noise $\zeta$ is in that case $(-1^-)$-H\"older continuous, and the intuition suggests that the solution $u$ to the (PAM) equation should be $(1^-)$-H\"older continuous, as a consequence of the regularizing effect of the heat semigroup. So at small time-space scales, $u$ should essentially be constant, as a first approximation. This could suggest to try a perturbative approach in which, if one denotes by $Z$ the solution to the equation $(\partial_t+\Delta)Z = \zeta$, with null initial condition, one looks for a distribution/function $v := u-Z$ with better regularity than the expected regularity of $u$, so as to get a well-posed equation for $v$. Such an attempt is bound to fail as $v$ needs to satisfy the same equation as $u$. The same trick invented by Da Prato-Debbusche in their study of the $2$-dimensional stochastic quantization equation \cite{DPD}, also fails in the study of $3$-dimensional scalar $\Phi^4_3$ equation, but a local 'version' of this idea is at the heart of the theory of regularity structures, while a tilted version of that point of view is also the starting point of paracontrolled calculus. Both make sense, with different tools, of the fact that a solution should locally ``look like'' $Z$. Whereas 'usual' Taylor expansions are used in the theory of regularity structures to compare a distribution to a linear combination of some given model distributions constructed by purely probabilistic means, such as the a priori undefined product $Z\zeta$, the paracontrolled approach uses paraproducts as a means of making sense of the sentence ``$u$ looks like $Z$ at small scales'', such as given in the definition below. For readers unfamiliar with paraproducts, recall that any distribution $f$ can be described as an infinite sum of smooth functions $f_i$ with the Fourier transform $\widehat{f_i}$ of $f_i$ essentially equal to the restriction of $\widehat{f}$ on a compact annulus depending on $i$. A product of two distributions $f$ and $g$ can thus always be written formally as
\begin{equation}
\begin{split}
fg &= \sum f_i g_j = \sum_{i\leq j-2} f_i g_j + \sum_{|i-j|\leq 1} f_i g_j +\sum_{j\leq i-2} f_i g_j \\
    &=: \Pi_f(g) + \Pi(f,g) + \Pi_g(f).
\end{split}
\end{equation}
The term $\Pi_f(g)$ is called the paraproduct of $f$ and $g$, and the term $\Pi(f,g)$ is called the resonant term. The paraproduct is always well-defined for $f$ and $g$ in H\"older spaces, with possibly negative indices $\alpha$ and $\beta$ respectively, while the resonant term only makes sense if $\alpha+\beta>0$. (The book \cite{BCD} provides a gentle introduction to paraproducts and their use in the study of some classes of PDEs.) This result of Bony on paraproducts \cite{Bony} already offers a setting that extends Schwartz operation of multiplication of a distribution by a smooth function; it is not sufficient however for our needs, even for the (PAM) equation in dimension $2$, as $u$ is expected there to be $1^-$-H\"older and $\zeta$ is $(-1^-)$-H\"older in that case. Needless to say, things are even worse in dimension 3 and for Burgers system. However, the point is that we do not want to multiply any two distributions but rather very special pairs of distributions. A reference distribution $Z$ in some parabolic H\"older space $\calC^\alpha$, defined later, is given here.

\ssk

\begin{defn*}
Let $\beta>0$ be given. A pair of\textbf{ distributions} $(f,g)\in \calC^\alpha\times\calC^\beta$ is said to be \textbf{paracontrolled by $Z$} if 
$$
(f,g)^\sharp := f - \Pi_g(Z)\in\calC^{\alpha+\beta}.
$$ 
\end{defn*}

\ssk

The distribution $g$ is called the \textit{derivative of $f$ with respect to $Z$.} The following elementary remark gives credit to this choice of name. It also partly explains why we shall solve the (PAM) equation in the way we do it here -- using some kind of Cole-Hopf transform. Assume $\alpha$ is positive, and write $(2\alpha)$ for a function in $\calC^{2\alpha}$ that may change from line to line. \textit{For a pair $(f,f)$ paracontrolled by $Z$, one can write $f = e^Zg$, for some function $g$ in $\calC^{2\alpha}$.} It suffices indeed to notice that Bony's decomposition gives
\begin{equation*}
\begin{split}
e^{-Z}f &= \Pi_{e^{-Z}}(f) + \Pi_f\big(e^{-Z}\big) + (2\alpha)   \\
            &= \Pi_{e^{-Z}}\big(\Pi_f(Z)\big) + \Pi_f\big(\Pi_{-e^{-Z}}(Z)\big) + (2\alpha)  \\
            &= \Pi_{e^{-Z}f}(Z) - \Pi_{e^{-Z}f}(Z) + (2\alpha) = (2\alpha).
\end{split}
\end{equation*}
We used in the second and third equalities two elementary results on paraproducts which are well-known in the classical setting, and proved below in the more general setting of the present work. 

\ssk

The twist offered by this definition, as far as the multiplication problem of $u$ by $\zeta$ is concerned, is the following. Take for $Z$ the solution to the equation $(\partial_t+L)Z = \zeta$, with null initial condition; the noise $\zeta$ is thus here $(\alpha-2)$-H\"older. From purely analytic data, the product $u\zeta$ is meaningful only if $\alpha+(\alpha-2)>0$, that is $\alpha>1$. For a distribution $(u,u')$ controlled by $Z$, with $\beta=\alpha$ say, the  formal manipulation 
\begin{equation*}
\begin{split}
u\zeta &= \Pi_u(\zeta) + \Pi_\zeta(u) + \Pi(u,\zeta) \\
           &= \Pi_u(\zeta) + \Pi_\zeta(u) + \Pi\big(\Pi_{u'}(Z),\zeta\big) + \Pi\big(\rest{2\alpha},\zeta\big) \\
           &=: \Pi_u(\zeta) + \Pi_\zeta(u) + \textsf{C}(Z,u',\zeta) + u'\,\Pi(Z,\zeta) + \Pi\big(\rest{2\alpha},\zeta\big),
\end{split}
\end{equation*}
gives a decomposition of $u\zeta$ where the first two terms are always well-defined, with known regularity, and where the last term makes sense provided $2\alpha + (\alpha-2) > 0$, that is $\alpha>\frac{2}{3}$. It happens that the corrector
$$
\textsf{C}(Z,u',\zeta) := \Pi\big(\Pi_{u'}(Z),\zeta\big) - u'\Pi(Z,\zeta)
$$
can be proved to define an $\big(\alpha+\alpha+(\alpha-2)\big)$-H\"older distribution if $\alpha>\frac{2}{3}$, although the resonant term $\Pi\big(\Pi_{u'}(Z),\zeta\big)$ is only well-defined on its own if $\alpha>1$. So we see that the only undefined term in the decomposition of $u\zeta$ is the product $u'\,\Pi(Z,\zeta)$, where the resonant term $\Pi(Z,\zeta)$ does not make sense so far. This is where probability comes into play, to show that one can define a random distribution $\Pi(Z,\zeta)$ as a limit in probability  of renormalized quantities of the form $\Pi(Z^\epsilon,\zeta^\epsilon) - c^\epsilon$, where $\zeta^\epsilon$ is a regularized noise, with associated $Z^\epsilon$, and $c^\epsilon$ is a deterministic function, a constant in some cases. The convergence can be proved to hold in $\calC^{\alpha+(\alpha-2)}$, so the product $u'\,\Pi(Z,\zeta)$ eventually makes perfect sense if $\alpha+(2\alpha-2)>0$, that is $\alpha>\frac{2}{3}$. This combination of analytic and probabilistic ingredients shows that one can define the product $u\zeta$, or more properly $(u,u')\,\zeta$, for $\alpha>\frac{2}{3}$, which is definitely beyond the scope of Bony's paradigm. Once the distribution $\zeta$ has been enhanced into a pair $\widehat\zeta := \big(\zeta,\Pi(Z,\zeta)\big)$ with good analytic properties, one can define the product $(u,u')\,\widehat\zeta$ as above for a generic distribution paracontrolled by $Z$, and reformulate a singular PDE such as the (PAM) equation in dimension $2$ as a fixed point problem in some space of paracontrolled distribution, and solve it uniquely by a fixed point method. Note that the very notion of product, and hence the meaning of the equation, depends on the choice of enhancement of $\zeta$ into $\widehat{\zeta}$.

\medskip
 
The above reasoning will not be sufficient, however, to deal with the (PAM) and multiplicative Burgers equations in dimension $3$, for which $\alpha<\frac{1}{2}$, and one needs first to reformulate the equation differently to make it accessible to this first order expansion calculus. In analogy with Lyons' rough paths theory, and in parallel with the logical structure of the theory of regularity structures, one may also consider developing a \textsf{\textbf{higher order paracontrolled calculus}} where a collection of reference functions $(Z_1,..,Z_k)$, with increasing regularity (for example $Z_i$ of regularity $i\alpha$ for some $\alpha>0$), are given, and used to give some sort of Taylor expansion of a function $f\in\calC^\alpha$ of the form 
$$
(f,g_1,..,g_k)^\sharp := f - \left(\Pi_{g_1}(Z_1)+ ...+ \Pi_{g_k}(Z_k)\right) \in\calC^{k\alpha+\beta}.
$$ 
for some tuple $(g_1,..,g_k)$ of $\calC^\alpha$ functions with similar expansions at lower order. We shall develop this framework in a forthcoming work.

\bigskip

\section[\hspace{0.7cm} Geometric and functional setting]{Geometric and functional settings}
\label{SectionSetting}

We describe in this section the geometric and functional settings in which we shall construct our space-time paraproducts in Section \ref{SectionParaproducts}, and provide a number of tools. We shall work in a Riemannian setting under fairly general conditions; parabolic H\"older spaces are defined Section \ref{SubsectionParabolicHolderSpaces} purely in terms of the semigroup generated by $L$. In Section \ref{SubsectionSchauder} we prove some fundamental Schauder-type regularity estimates. The cancellation properties put forward in Section \ref{SubsectionCancellations} are fundamental for proving in Section \ref{SectionParaproducts} some continuity results for some iterated commutators and correctors.

\bigskip

\subsection[\hspace{-0.5cm} Riemannian framework]{Riemannian framework}
\label{SubsectionSubLaplacian}

Our basic setting in this work will be a complete volume doubling measured Riemannian manifold $(M,d,\mu)$; all kernels mentioned in the sequel are with respect to the fixed measure $\mu$. We are going to introduce in the sequel a number of tools to analyze singular partial differential equations involving a parabolic operator on $\R_+\times M$
$$
\mathscr{L} := \partial_t + L,
$$
with $L$ built from first order differential operators $(V_i)_{i=1..\ell_0}$ on $M$, defined as operators $V_i$ satisfying the Leibniz rule
\begin{equation} \label{eq:Leibniz} 
V_i(fg) = f V_i(g) + g V_i(f)
\end{equation} 
for all functions $f,g$ in the domain of $L$. Given a tuple $I = (i_1,\dots,i_k)$ in $\{1,\dots,\ell_0\}^k$, we set $|I| := k$ and 
$$
V_I := V_{i_k}\cdots V_{i_1}.
$$

\bigskip

\noindent {\bf Assumption (A)} We shall assume throughout that 
\begin{itemize}
\item {\sf the operator $L$ is a sectorial operator in $L^2(M)$, $L$ is injective on $L^2(M)$, or the quotient space of $L^2(M)$ by the space of constant functions if $\mu$ is finite, it has a bounded $H^\infty$-calculus on $L^2(M)$, and  $-L$ generates a holomorphic semigroup $(e^{-tL})_{t>0}$ on $L^2(M)$}, \vspace{0.2cm}

\item {\sf one has 
$$
L = -\sum_{i=1}^{\ell_0} V_i^2
$$
and
$$
\mcD(L) \subseteq \mcD(V_i^2) := \big\{f\in L^2, V_i^2(f)\in L^2(\mu)\big\},
$$ 
for some operators $V_i$ satisfies the Leibniz rule \eqref{eq:Leibniz} on $\mcD(L)$,} \vspace{0.2cm} 

  \item {\sf the heat semigroup is conservative, that is $\big(e^{-tL}\big)({\bf 1}_M) = {\bf 1}_M$ for every $t>0$, where $1_M$ stands for the constant function on $M$ -- or in a weak sense that $L({\bf 1}_M)=0$},  \vspace{0.2cm}
   
   \item {\sf the semigroup has regularity estimates at any order}, by which we mean that for every tuple $I$, the operators $\Big(t^\frac{| I |}{2}V_I\Big)e^{-tL}$ and $e^{-tL}\Big(t^\frac{| I |}{2}V_I\Big)$ have kernels $K_t(x,y)$ satisfying the Gaussian estimate
\begin{equation}
\label{EqConditionGaussianEstimate}
\Big| K_t(x,y) \Big| \lesssim  \frac{1}{\mu\big(B(x,\sqrt{t})\big)}\, e^{-c\,\frac{d(x,y)^2}{t}}
\end{equation}
and the following regularity estimate. For $d(x,z)\leq \sqrt{t}$
\begin{equation}
\label{finite}
\Big| K_t(x,y) - K_t(z,y)\Big| \lesssim \frac{d(x,z)}{\sqrt{t}}  \frac{1}{\mu\big(B(x,\sqrt{t})\big)} \, e^{-c\,\frac{d(x,y)^2}{t}},
\end{equation}
for some constants which may depend on $|I|$.  
\end{itemize} 

\medskip

\noindent  Let us point out that the regularity property \eqref{finite} for $|I|=k$ can be obtained from \eqref{EqConditionGaussianEstimate} with $k+1$ writing the ``finite-increments'' formula
$$ 
\big|K_t(x,y)-K_t(z,y)\big| \lesssim d(x,z) \ \sup_j \ \sup_{w \in (x,z)} \big|X_j K_t(w,y)\big|
$$
where $(x,y)$ stands for a geodesic joining $x$ to $z$ and of length $d(x,z)$, and $(X_j)$ stands for a local frame field near $(x,y)$, and it acts here as a first order differential on the first component of $K$. As a matter of fact, it suffices for the present work to assume that the semigroup has regularity estimates of large enough order. Observe that under Assumption \ref{assumption1}, the semigroup $e^{-tL}$ may be defined as acting on the distributions. This can be rigorously done by duality, since for every integer $N\geq 0$, $(L^{*})^N e^{-tL^*}$ has a kernel satisfying pointwise Gaussian estimates and so is acting on the test functions $\mcS_o$ -- we refer to \cite{BDY} for more details. One can keep in mind the following two examples.

\medskip

\begin{enumerate}
   \item \textbf{Euclidean domains.} In the particular case of the Euclidean space, all of the current work can be reformulated in terms of Fourier transform rather than in terms of the heat semigroup; which may make some reasoning a bit more familiar but does not really simplify anything. The case of a bounded domain with its Laplacian associated with Neumann boundary conditions fits our framework if the boundary is sufficiently regular. We may also consider other kind of second order operator, like $L=-\div(A \nabla)$ for some smooth enough matrix-valued map satisfying the ellipticity (or accretivity) condition. 
      \vspace{0.2cm}

   \item \textbf{Riemannian manifolds.} Smooth closed manifolds equipped with an operator $L$ of Hörmander type as above, with $V_i$ smooth, all satisfy assumptions {\bf (A)}. Here is a simple setting within which one can deal with unbounded spaces. Assume $M$ is a parallelizable $d$-dimensional manifold with a smooth global frame field $V = (V_1,\dots,V_d)$. One endows $M$ with a Riemannian structure by turning $V$ into orthonormal frames. The above assumption on the heat kernel holds true if $M$ has bounded geometry, that is if \vspace{0.1cm}
   \begin{itemize}
      \item[\textsf{(i)}] the curvature tensor and all its covariant derivatives are bounded in the frame field $V$,  \vspace{0.1cm}
      \item[\textsf{(ii)}] Ricci curvature is bounded from below, \vspace{0.1cm}
      \item[\textsf{(iii)}] and $M$ has a positive injectivity radius;  \vspace{0.1cm}
\end{itemize}
see for instance \cite{CRT} or \cite{Triebel}. One can include the Laplace operator in this setting by working with its canonical lift to the orthonormal frame bundle, given by $\frac{1}{2}\sum_{i=1}^d H_i^2 + \frac{1}{2}\sum_{1\leq j<k\leq d} V_{jk}$, where the $H_i$ are the canonical horizontal vector fields of the Levi-Civita connection, and the $V_{jk}$ are the canonical vertical vector fields on the orthonormal frame bundle, inherited from its $SO(\R^d)$-principal bundle structure. The bundle $OM$ is parallelizable and satisfies assumptions {\bf (A)} if the base Riemannian manifold $M$ satisfies the conditions {\sf (i--iii)}. This example shows that the assumption that $M$ is parallelizable is essentially done here for convenience.
\end{enumerate}

\bigskip

\subsection[\hspace{-0.5cm} Approximation operators and cancellation property]{Approximation operators and cancellation property}
\label{SubsectionCancellations}

We introduce in this section a notion of approximation operators that will be the building blocks for the definition and study of the paraproducts, commutators and correctors, used in our analysis of singular PDEs. Some of them enjoy some kind of orthogonality, or cancellation, property quantified by condition \eqref{eq:mcOa} below. This property is to be thought of as a quantitative replacement for the property of frequency localisation of Fourier multipliers in the Littlewood-Paley decomposition heavily used in the classical Fourier definition of paraproducts; see \cite{BCD}. Indicators of annulii will somehow be replaced in our setting by continuous functions of order $1$ on such annulii, with exponential decrease at $0$ and $\infty$. Note that we shall be working in a parabolic setting with mixed cancellation effects in time and space.

\ssk

All computations below make sense for a choice of large enough integers $b,\ell_1$ that will definitely be fixed at the end of Section \ref{SubsectionIntertwined} to ensure some continuity properties for some useful operators. Recall that generic elements of the parabolic space $\calM=\RR\times M$ are denoted by $e=(x,\tau)$ or $e'=(y,\sigma)$, and that $t$ stands for a scaling parameter. The following parabolic Gaussian-like kernels $(\G_t)_{0<t\leq 1}$ will be used as reference kernels in this work. For $0<t\leq 1$ and $\sigma\leq\tau$, if $d(x,y)\leq 1$, set 
$$
\G_t\big((x,\tau),(y,\sigma)\big):= \frac{1}{\nu\Big(B_\mcM \big((x,\tau),\sqrt{t}\big) \Big)}\,\left(1+\frac{\rho \big((x,\tau), (y,\sigma)\big)^2}{t} \right)^{-\ell_1},  
$$
otherwise we set 
\begin{equation*}
\begin{split}
\G_t\big((x,\tau)&,(y,\sigma)\big) :=  \\
&\frac{1}{\nu\Big(B_\mcM \big((x,\tau),1\big) \Big)} \left(1+\frac{|\tau-\sigma|}{t} \right)^{-\ell_1} \left(1+\frac{d(x,y)^2}{t} \right)^{-\ell_1}\,\exp\left(-c \, \frac{d(x,y)^2}{t}\right)
\end{split}
\end{equation*}
for $d(x,y)\geq 1$, and $\G_t \equiv 0$ if $\tau\leq \sigma$. We do not emphasize the dependence of $\G$ on the positive constant $c$ in the notation for the 'Gaussian' kernel, and we shall allow ourselves to abuse notations and write $\G_t$ for two functions corresponding to two different values of that constant. This will in particular be the case in the proof of Lemma \ref{lem:weight}. We have for instance, for two scaling parameters $s,t\in (0,1)$, the estimate
\begin{equation}
\label{EqIteratedG}
\int_\mcM \G_{t}\big(e,e'\big) \, \G_{s}\big(e',e''\big) \, \nu(de') \lesssim \G_{t+s}\big(e,e''\big).
\end{equation}
(Indeed, the space variables and the time variables are separated in the kernel $\G_t$. Then both in space and time variables, the previous inequality comes from classical estimates for convolution of functions with fast decay at infinity, such as done in \cite[Lemma A.5]{BB} for example.) This somewhat unnatural definition of a Gaussian-like kernel is justified by the fact that we shall mainly be interested in local regularity matters; the definition of $\G_t$ in the domain $\big\{d(x,y)\geq 1\big\}$ is only technical and will allow us to obtain global estimates with weights. Presently, note that a large enough choice of constant $\ell_1$ ensures that we have
\begin{equation} 
\sup_{t\in(0,1]} \sup_{e\in \mcM} \ \int_\mcM \G_{t}(e,e')\, \nu(de') <\infty,
 \label{eq:Gt} \end{equation}
so any linear operator on a function space over $\mcM$, with a kernel pointwisely bounded by some $\G_t$ is bounded in $L^{p}(\nu)$ for every $p\in[1,\infty]$.

\medskip

\begin{defn*}
We shall denote throughout by ${\sf G}$ the set of families $(\mcP_t)_{0<t\leq 1}$ of linear operators on $\mcM$ with kernels pointwisely bounded by
$$ 
\Big| K_{\mcP_t}(e,e') \Big| \lesssim  \G_{t}(e,e').
$$
\end{defn*}

The letter {\sf G} is chosen for 'Gaussian'. A last bit of notation is needed before we introduce the cancellation property for a family of operators in a parabolic setting. Given a real-valued integrable function $m$ on $\R$, define its rescaled version as
$$
m_t(\cdot) := \frac{1}{t}\, m\Big(\frac{\cdot}{t}\Big);
$$
the family $(m_t)_{0<t\leq 1}$ is uniformly bounded in $L^1(\R)$. We also define the ``convolution'' operator $m^\star$ associated with $m$ via the formula
$$
m^\star(f)(\tau) := \int_{0}^\infty m(\tau-\sigma)f(\sigma)d\sigma.
$$
Note that if $m$ has support in $\R_+$, then the operator $m^\star$ has a kernel supported on the same set $\big\{(\sigma,\tau)\,; \sigma\leq \tau\big\}$ as our Gaussian-like kernel. Moreover, we let the reader check that if $m_1,m_2$ are two $L^1$-functions with $m_2$ supported on $[0,\infty)$, with convolution $m_1 \ast m_2$, then we have
$$
\big(m_1 \ast m_2\big)^\star = m_1^\star \circ m_2^\star.
$$
Given an integer $b\geq 1$, we define a special family of operators on $L^2(M)$ setting $\gamma_b := (b-1)!$ and 
$$ 
Q_t^{(b)} := \gamma_b^{-1} (tL)^b e^{-tL} \qquad \textrm{and} \qquad  -t\partial_t P_t^{(b)} := Q_t^{(b)},
$$
with $P^{(b)}_0 = \textrm{Id}$, so $P_t^{(b)}$ is an operator of the form $p_b(tL)e^{-tL}$, for some polynomial $p_b$ of degree $(b-1)$, with value $1$ in $0$. Under Assumption \ref{assumption1}, the operators $P_t^{(b)}$ and $Q_t^{(b)}$ both satisfy the Gaussian regularity estimates \eqref{EqConditionGaussianEstimate} at any order
\begin{equation}
\label{EqGaussianEstimateKQt}
\left| K_{t^\frac{| I |}{2}V_IR} (x,y) \right| \vee \left| K_{t^\frac{| I |}{2}RV_I} (x,y) \right| \lesssim  \frac{1}{\mu\big(B(x,\sqrt{t})\big)} \, e^{-c\,\frac{d(x,y)^2}{t}},
\end{equation}
with $R$ standing here for $P_t^{(b)}$ or $Q_t^{(b)}$.  \\

The parameters $b$ and $\ell_1$ will be chosen large enough, and fixed throughout the paper. See Proposition \ref{prop:modifiedparaproduit} and the remark after Proposition \ref{PropContinuityModifParapdct} for the precise choice of $b$ and $\ell_1$. 

\medskip

\begin{defn*} {\sf  
Let an integer $a\in\llbracket 0,2b \rrbracket$ be given. The following collection of families of operators is called the \textbf{standard collection of operators with cancellation of order $a$}, denoted by ${\sf StGC}^a$. It is made up of all the space-time operators
$$ 
\Big( \big(t^\frac{|J|}{2} V_J\big) (tL)^\frac{a-|J|-2k}{2} P_t^{(c)}  \otimes m^\star_t \Big)_{0<t\leq 1} 
$$
where $k$ is an integer with $2k+|J| \leq a$, and $c\in \llbracket 1,b\rrbracket$, and $m$ is any smooth function supported on $\big[\frac{1}{2},2\big]$ such that
\begin{equation} 
\label{eq:moment} 
\int \tau^i m(\tau) \, d\tau=0,  
\end{equation}
for all $0\leq i\leq k-1$, with the first $b$ derivatives bounded by $1$. These operators are uniformly bounded in $L^p(\mcM)$ for every $p\in [1,\infty]$, as functions of the scaling parameter $t$. So a standard collection of operators $\mcQ$ can be seen as a bounded map $\mcQ: t\rightarrow \mcQ_t$ from $(0,1]$ to the set ${\mathcal B}(L^p)$ of bounded linear operators on $L^p(\mcM)$.  We also set 
$$
{\sf StGC}^{[0,2b]} := \bigcup_{0\leq a \leq 2b} {\sf StGC}^a.
$$   }
\end{defn*}

\ssk

The cancellation effect of such operators is quantified in Proposition \ref{prop:SO} below; note here that it makes sense at an intuitive level to say that $L^\frac{a-|J|-2k}{2}$ encodes cancellation in the space-variable of order $a-|J|-2k$, that $V_J$ encodes a cancellation in space of order $|J|$  and that the moment condition \eqref{eq:moment} encodes a cancellation property in the time-variable of order $k$ for the convolution operator $m^\star_t$. Since we are in the parabolic scaling, a cancellation of order $k$ in time corresponds to a cancellation of order $2k$ in space, so that $V_J L^\frac{a-|J|-2k}{2} P_t^{(c)}  \otimes m^\star_t$ has a space-time cancellation property of order $a$. We invite the reader to check that each operator $\big(t^\frac{|J|}{2} V_J\big) (tL)^\frac{a-|J|-2k}{2} P_t^{(c)}  \otimes m^\star_t$ in the standard collection has a kernel pointwisely bounded from above by some $\G_t$. This justifies the choice of name ${\sf StGC}^a$ for this space, where {\sf St} stands for 'standard', {\sf G} for 'Gaussian' and {\sf C} for 'cancellation'. The paracontrolled analysis, that we are going to explain, is based on these specific operators. We emphasize that because of the Gaussian kernel $\G_t$ and the function $m$, all of these operators have a support in time included in
$$ 
\{(\tau,\sigma),\ \tau \geq \sigma\}.
$$
In particular, that means that we never travel  backwards in time through these operators. This fact will be very important, to deal further with the weight $\varpi$, which will depend on time. We give one more definition before stating the cancellation property.

\ssk

\begin{defn*} {\sf 
Given an operator $Q := V_I\,\phi(L)$, with $|I|\geq 1$, defined by functional calculus from some appropriate function $\phi$, we write $Q^\bullet$ for the \textbf{formal dual operator}
$$ 
Q^{\bullet} :=  \phi(L) V_I.
$$
For $I= \emptyset$, and $Q = \phi(L)$, we set $Q^\bullet := Q$. For an operator $Q$ as above we set
$$
\big(Q\otimes m^\star\big)^\bullet := Q^\bullet\otimes m^\star.
$$  }
\end{defn*}

\ssk

Note that the above definition is \textit{not} related to any classical notion of duality, and let us emphasize that we do \textit{not assume} that $L$ is self-adjoint in $L^2(\mu)$. This notation is only used to indicate that an operator $Q$, resp. $Q^\bullet$, can be composed on the right, resp. on the left, by another operator $\psi(L)$, for a suitable function $\psi$, due to the functional calculus on $L$. In the setting of analysis on a finite dimensional torus, the operators $Q^{(b)}_t$ are given in Fourier coordinates $\lambda$, as the multiplication operators by $(t|\lambda|^2)^be^{-t|\lambda|^2}$; as this function is almost localized in an annulus $|\lambda|\sim t^{-\frac{1}{2}}$, the operators $Q_s^{(b)}$ and $Q_t^{(b)}$ are almost orthogonal if $\frac{s}{t}$ is either very small or very big. This is encoded in the elementary estimate 
\begin{equation}
\label{EqElementaryCancellation}
\left| K_{Q^{(b)}_{s} \circ Q^{(b)}_{t}} (x,y) \right| \lesssim  \left(\frac{ts}{(s+t)^2}\right)^\frac{b}{2} \frac{1}{\mu\big(B(x,\sqrt{s+t})\big)}\,\exp\left(-c\,\frac{d^2(x,y)}{t+s}\right).
\end{equation}  
The frequency analysis of the operators $\mcQ_s^{(b)}$ is not very relevant in the non-homogeneous parabolic space $\mcM$. We keep however from the preceeding analysis the idea that relation \eqref{EqElementaryCancellation} encodes some kind of orthogonality, or cancellation effect.

\medskip

\begin{prop} \label{prop:SO} {\sf 
Consider $\mcQ^1\in {\sf StGC}^{a_1}$ and $\mcQ^2\in {\sf StGC}^{a_2}$ two standard collections with cancellation, and set $a:=\min(a_1,a_2)$. Then for every $s,t\in(0,1]$, the composition $\mcQ^1_{s} \circ \mcQ^{2\bullet}_{t}$ has a kernel pointwisely bounded by
\begin{equation}
\label{prop:cancellation}
\left| K_{\mcQ^1_{s}\circ\mcQ^{2\bullet}_{t}} (e,e') \right| \lesssim  \left(\frac{ts}{(s+t)^2}\right)^\frac{a}{2} \G_{t+s}(e,e'). 
\end{equation}   }
\end{prop}

\medskip

\begin{Dem} 
Given
$$ 
\mcQ^1_{s} =  s^\frac{j_1}{2} V_{J_1}(sL)^\frac{a_1-j_1-2k_1}{2} P_s^{(c_1)}  \otimes m^{1\star}_s     \quad \textrm{and}   \quad  \mcQ^{2\bullet}_{t} =  (t L)^\frac{a_2-j_2-2k_2}{2} P_t^{(c_2)} t^\frac{j_2}{2} V_{J_2} \otimes \,m^{2\star}_t
$$
a standard operator and the dual of another, we have
$$ 
\mcQ^1_{s} \circ\mcQ^{2\bullet}_{t} = s^\frac{a_1-2k_1}{2} t^\frac{a_2-2k_2}{2}  V_{J_1}L^\frac{a_1-j_1-2k_1+a_2-j_2-2k_2}{2} P_s^{(c_1)}P_t^{(c_2)} V_{J_2} \otimes \big(m^{1}_s \ast m^{2}_t\big)^\star.
$$
Assume, without loss of generality, that $0<s \leq t$. Then the kernel of the time-convolution operator $m^{(1)}_s \ast m^{(2)}_t$ is given by
$$ 
K_{m^{1}_s * m^{2}_t}(\tau-\sigma)=\int m^{1}\left(\frac{\tau-\lambda}{s}\right) m^{2} \left(\frac{\lambda-\sigma}{t}\right) \, \frac{d\lambda}{st}.
$$
Since $m^{1}$ has vanishing $k_1$ first moments, we can perform $k_1$ integration by parts and obtain that
$$ 
\left|K_{m^1_s \ast m^2_t}(\tau,\sigma)\right| \lesssim \left(\frac{s}{t}\right)^{k_1} \int \partial^{-k_1} m^{1}\left(\frac{\tau-\lambda}{s}\right)  \partial^{k_1} m^{2} \left(\frac{\lambda-\sigma}{t}\right) \, \frac{d\lambda}{st},
$$
where we slightly abuse notations and write $\partial^{-k_1}m^{1}$ for the $k_1^\textrm{th}$ primitive of $m^{1}$ null at $0$. Then we get
\begin{align*}
 \left|K_{m^{1}_s\ast m^{2}_t}(\tau,\sigma)\right| & \lesssim \left(\frac{s}{t}\right)^{k_1} \int \left(1+\frac{|\tau-\lambda|}{s}\right)^{-\ell_1+2}  \left(1+\frac{\lambda-\sigma}{t}\right)^{-\ell_1+2} \, \frac{d\lambda}{st} \\
 & \lesssim \left(\frac{s}{t}\right)^{k_1}  \left(1+\frac{|\tau-\sigma|}{s+t}\right)^{-\ell_1} (s+t)^{-1}. 
\end{align*}
In the space variable, the kernel of $V_{J_1}L^\frac{a_1-j_1-2k_1+a_2-j_2-2k_2}{2} P_s^{(c_1)}P_t^{(c_2)}V_{J_2}$ is bounded above by
$$
(s+t)^\frac{-a_1+2k_1-a_2+2k_2}{2} \, \mu\Big(B(x,\sqrt{s+t})\Big)^{-1} \exp\left(-c \frac{d(x,y)^2}{s+t}\right),
$$
as a consequence of the property \eqref{EqGaussianEstimateKQt}. Altogether, this gives
\begin{align*}
\left| K_{\mcQ^1_{s} \circ \mcQ^{2\bullet}_{t}} (e,e') \right| & \lesssim \left(\frac{s}{t}\right)^{k_1} s^\frac{a_1-2k_1}{2} t^\frac{a_2-2k_2}{2} (s+t)^\frac{-a_1+2k_1-a_2+2k_2}{2} \G_{t+s} (e,e') \\
 & \lesssim \left(\frac{s}{t}\right)^\frac{a_1}{2}\G_{t+s} (e,e') \\
 & \lesssim \left(\frac{s}{t}\right)^\frac{a}{2} \G_{t+s} (e,e'),
\end{align*}
where we used that $s\leq t$ and $a\leq a_1$.
\end{Dem}

\medskip

\begin{defn*} {\sf
Let $0\leq a\leq 2b$ be an integer. We define the subset ${\sf GC}^a$ of ${\sf G}$ of \textbf{families of operators with the cancellation property of order $a$} as the set of elements $\mcQ$ of ${\sf G}$ with the following cancellation property. For every $0<s,t\leq 1$ and every standard family $\mcS\in {\sf StGC}^{a'}$, with $a'\in\llbracket a,2b\rrbracket$, the operator $\mcQ_t\circ\mcS_s^\bullet$ has a kernel pointwisely bounded by
\begin{equation} 
\label{eq:mcOa}
\Big| K_{\mcQ_t \circ\mcS_s^\bullet}(e,e') \Big| \lesssim \left(\frac{st}{(s+t)^2}\right)^\frac{a}{2} \G_{t+s}(e,e'). 
\end{equation}   }
\end{defn*}

\medskip

Here are a few examples. Consider a smooth function $m$ with compact support in $[2^{-1},2]$, an integer $c\geq 1$, and a tuple $I$ of indices.
\begin{itemize}
    \item The families $\Big(Q_t^{(\frac{a}{2})} \otimes m_t^\star\Big)_{0<t\leq 1}$ and $\Big(t^\frac{|I|}{2} V_I P_t^{(c)} \otimes m_t^\star\Big)_{0<t\leq 1}$ belong to ${\sf GC}^a$ if $|I|\geq a$;   \vspace{0.15cm}
 
   \item If $\int \tau^k m(\tau) \, d\tau=0$ for all integer $k=0,...,a-1$, then we can see by integration by parts along the time-variable that $\big(P_t^{(c)} \otimes m_t^\star\big)_{0<t\leq 1} \in {\sf GC}^a$.   \vspace{0.15cm}
 
   \item If $\int \tau^k m(\tau) \, d\tau=0$ for all integer $k=0,...,a_2$ with $a_1+a_2=a$, then the families $\Big(Q_t^\frac{a_1}{2} \otimes m_t^\star\Big)_{0<t\leq 1}$ and $\Big(t^\frac{|I|}{2} V_I P_t^{(c)} \otimes m_t^\star\Big)_{0<t\leq 1}$, where $|I|\geq a_1$, both belong to $ {\sf GC}^a$.   \vspace{0.15cm}
\end{itemize}

We see on these examples that cancellation in the parabolic setting can encode some cancellations in the space variable, the time-variable or both at a time.

\medskip

We introduced above the operators $Q^{(b)}_t$ and $P^{(b)}_t$ acting on the base manifold $M$. We end this section by introducing their parabolic counterpart. Choose arbitrarily a smooth real-valued function $\varphi$ on $\R$, with support in $\big[\frac{1}{2},2\big]$, unit integral and such that for every integer $k=1,\dots,b$, we have
$$ 
\int \tau^k \varphi(\tau) \, d\tau =0.
$$ 
Set 
$$
\mcP_t^{(b)}:= P^{(b)}_t \otimes \varphi_t^\star \qquad \textrm{and} \qquad \mcQ_t^{(b)}:= -t \partial_t \mcP^{(b)}_t,
$$
Denote by ${\sf M}_\tau$ the multiplication operator in $\R$ by $\tau$. An easy computation yields that
$$ 
\mcQ_t^{(b)} = Q_t^{(b)} \otimes \varphi_t^\star + P_t^{(b)} \otimes \psi_t
$$
where $\psi(\sigma) := \varphi(\sigma) + \sigma\varphi'(\sigma)$. (For an extension of the present theory to the setting of Sobolev spaces, such as done in the appendix B of \cite{BB}, it would be convenient to work with $\varphi\ast\varphi$ rather than $\varphi$.) Note that, from its very definition, a parabolic operator $\mcQ_t^{(b)}$ belongs at least to ${\sf GC}^{2}$, for $b\geq 2$. Remark that if $\zeta$ is a time-independent distribution then $\mcQ^{(b)}_t\zeta = Q_t^{(b)}\zeta$. Note also that due to the normalization of $\varphi$, then for every $f\in L^p ({\mathbb R})$ supported on $[0,\infty)$ then
$$
\varphi^\star_t(f) \xrightarrow[t\to 0]{} f \qquad \textrm{in $L^p$}.
$$
So, the operators $\mcP_t$ tend to the identity as $t$ goes to $0$, on the set of functions $f\in L^p(\mcM)$ with time-support included in $[0,\infty)$, whenever $p \in [1,\infty)$, and on the set of functions $f\in C^0(\mcM)$ with time-support included in $[0,\infty)$. The following {\bf Calder\'on reproducing formula} follows as a consequence. For every continuous function $f\in L^\infty(\mcM)$ with time-support in $[0,\infty)$, we have 
\begin{equation}
\label{eq:calderon}
 f = \int_0^1 \mcQ^{(b)}_t (f) \, \frac{dt}{t} + \mcP^{(b)}_1(f).
\end{equation}

This formula will play a fundamental role for us. Noting that the measure $\frac{dt}{t}$ gives unit mass to intervals of the form $\big[2^{-(i+1)},2^{-i}\big]$, and considering the operator $\mcQ^{(b)}_t $ as a kind of multiplier roughly localized at frequencies of size $t^{-\frac{1}{2}}$, Calder\'on's formula appears as nothing else than a continuous time analogue of the Paley-Littlewood decomposition of $f$, with  $\frac{dt}{t}$ in the role of the counting measure.

\bigskip

\subsection[\hspace{-0.5cm} Parabolic H\"older spaces]{Parabolic H\"older spaces}
\label{SubsectionParabolicHolderSpaces}

We define in this section space and space-time weighted H\"older spaces, with possibly negative regularity index, and give a few basic facts about them. The setting of weighted function spaces is needed for the applications to the parabolic Anderson model and multiplicative Burgers equations on unbounded domains studied in Section \ref{SectionPAM}. The weights we use were first introduced in \cite{HLR3}.

\medskip

Let us start recalling the following well-known facts about H\"older spaces on $M$, and single out a good class of weights on $M$. A function $w : M \rightarrow [1,\infty)$ will be called a \textbf{spatial weight} if one can associate to any positive constant $c_1$ a positive constant $c_2$ such that one has
\begin{equation}  
\label{eq:weight}
\begin{split}
w(x) \, e^{-c_1 d(x,y)} \leq c_2 \, w(y),
\end{split}
\end{equation}
for all $x,y$ in $M$. Given $0<\alpha\leq 1$, the classical metric H\"older space $H^\alpha_w$ is defined as the set of real-valued functions $f$ on $M$ with finite $H^\alpha_w$-norm, defined by the formula
$$
\|f\|_{H_\omega^\alpha} := \big\|w^{-1} f\big\|_{L^\infty(M)} + \sup_{0<d(x,y)\leq 1} \frac{\big| f(x) - f(y) \big|}{w(x)\, d(x,y)^\alpha} <\infty.
$$
Distributions on $M$ were defined in \cite{BB} using a very similar definition as in the end of Section \ref{SectionIntroduction}, where their parabolic counterpart is defined.

\ssk

\begin{defn*} \label{def:C} {\sf 
For $\alpha\in (-3,3)$ and $w$ a spatial weight, define $C^\alpha_w := C^\alpha_w(M)$ as the set of distributions on $M$ with finite $C^\alpha_w$-norm, defined by the formula
$$ 
\|f\|_{C^\alpha_\omega} := \Big\|w^{-1} e^{-L}f\Big\|_{L^\infty(M)} + \sup_{0<t\leq 1} t^{-\frac{\alpha}{2}} \Big\|w^{-1} Q_t^{(a)}f \Big\|_{L^\infty(M)},
$$
and equip that space with the induced norm. The latter does not depend on the integer $a> \frac{|\alpha|}{2}$, and one can prove that the two spaces $H^\alpha_w$ and $C^\alpha_w$ coincide and have equivalent norms when $0<\alpha<1$ -- see \cite{BB}.  }
\end{defn*}

\medskip

These notions have parabolic counterparts which we now introduce. A \textbf{space-time weight} is a function $\omega : \mcM \rightarrow [1,\infty)$ with $\omega(x,\cdot)$ non-decreasing function of time, for every $x\in M$, and such that there exists two constants $c_1$ and $c_2$ with
\begin{equation} \label{eq:stweight}
\begin{split}
\omega(x,\tau)\, e^{-c_1 d(x,y)} \leq c_2\,\omega(y,\tau),
\end{split}
\end{equation}
for all pairs of points of $\mcM$ of the form $\big((x,\tau), (y,\tau)\big)$. The function $w_\tau := \omega(\cdot,\tau)$ is in particular a spatial weight for every time $\tau$. For $0<\alpha\leq 1$ and a space-time weight $\omega$, the metric parabolic H\"older space $\mcH^\alpha_\omega = \mcH^\alpha_\omega(\mcM)$ is defined as the set of all functions on $\mcM$ with finite $\mcH^\alpha_\omega$-norm, defined by the formula
$$ 
\|f\|_{\mcH^\alpha_\omega} := \big\|\omega^{-1} f\big\|_{L^\infty(\mcM)} + \sup_{0<\rho\big((x,\tau),(y,\sigma)\big)\leq 1;\,\tau\geq \sigma} \;\; \frac{|f(x,\tau)-f(y,\sigma)|}{\omega(x,\tau) \, \rho\big((x,\tau),(y,\sigma)\big)^\alpha}.
$$
As in the above spatial setting, one can recast this definition in a functional setting, using the parabolic standard operators. This requires the use of the following elementary result. Recall that the kernels $\G_t$ depend implicitly on a constant $c$ that may take different values with no further mention of it. We make this little abuse of notation in the proof of this statement.

\ssk

\begin{lem} \label{lem:weight} {\sf 
Let $A$ be a linear operator on $\mcM$ with a kernel $K_A$ pointwisely bounded by a Gaussian kernel $\G_t$, for some $t\in(0,1]$. Then for every space-time weight $\omega$, we have
$$ 
\big\| \omega^{-1} Af\big\|_{L^\infty(\mcM)} \lesssim \big\| \omega^{-1} f\big\|_{L^\infty(\mcM)}.
$$   }
\end{lem}

\ssk

\begin{Dem}
Indeed, for every $(x,\tau)\in \mcM$ we have 
\begin{align*}
\frac{1}{ \omega(x,\tau)}\, \big|(Af)(x,\tau)\big| & \lesssim \int_\mcM \G_t\big((x,\tau),(y,\sigma)\big) \frac{\omega(y,\sigma)}{\omega(x,\tau)}\,  \frac{\big| f(y,\sigma)\big|}{\omega(y,\sigma)} \, \nu(dy d\sigma) \\
 & \lesssim \int_\mcM \G_t\big((x,\tau),(y,\sigma)\big) \frac{\omega(y,\sigma)}{\omega(x,\sigma)} \, \frac{\big| f(y,\sigma)\big|}{\omega(y,\sigma)} \, \nu(dy d\sigma) \\
  & \lesssim \int_\mcM  \G_t\big((x,\tau),(y,\sigma)\big) \, \frac{\big| f(y,\sigma)\big|}{\omega(y,\sigma)} \, \nu(dyd\sigma) \\
  & \lesssim \big\|\omega^{-1} f\big\|_\infty,
  \end{align*}
where 
\begin{itemize}
   \item we used in the second inequality the fact that the function $\omega(x,\cdot)$ of time is non-decreasing, and $\G_t$ is null if $\sigma\geq \tau$, \vspace{0.1cm}
   
   \item the implicit constant in $\G_t$ was changed in the right hand side of the third inequality, and we used the growth condition \eqref{eq:stweight} on $\omega$ as a function of its first argument here, \vspace{0.1cm}
   
   \item we used the uniform bound \eqref{eq:Gt} on a Gaussian integral in the last line.
\end{itemize}
\end{Dem}

\medskip

Recall that distributions were introduced in the end of Section \ref{SectionIntroduction}.

\medskip

\begin{defn*} \label{def:calC} {\sf 
For $\alpha\in (-3,3)$ and a space-time weight $\omega$, we define the parabolic H\"older space $\calC^\alpha_\omega := \calC^\alpha_\omega(\mcM)$ as the set of distributions with finite $\calC^\alpha_\omega$-norm, defined by 
$$ 
\|f\|_{\calC_\omega^\alpha} := \sup_{\genfrac{}{}{0pt}{}{\mcQ \in {\sf StGC}^{k}}{0\leq k\leq 2b}} \big\|\omega^{-1} \mcQ_1(f)\big\|_{L^\infty(\mcM)} +  \sup_{\genfrac{}{}{0pt}{}{\mcQ \in {\sf StGC}^{k}}{ |\alpha|< k\leq 2b}} \ \sup_{0<t\leq 1} t^{-\frac{\alpha}{2}} \big\|\omega^{-1} \mcQ_t(f) \big\|_{L^\infty(\mcM)},
$$   
equipped with the induced norm.  }
\end{defn*}

\ssk

The restriction $\alpha\in(-3,3)$ is irrelevant and will be sufficient for our purpose in this work; taking $b$ large enough we can allow regularity of as large an order as we want. Building on Calder\'on's formula \eqref{eq:calderon}, one can prove as in \cite{BB} that the two spaces $\mcH^\alpha_\omega$ and $\calC^\alpha_\omega$ coincide and have equivalent norms, when $0<\alpha<1$. 

\ssk

\begin{prop} {\sf 
For $\alpha\in(0,1)$ and every space-time weight $\omega$,  the two spaces $\mcH^\alpha_\omega$ and $\calC^\alpha_\omega$ coincide and have equivalent norms.   }
\end{prop}

\medskip

\begin{Dem}
We first check that $\mcH^\alpha_\omega$ is continuously embedded into $\calC^\alpha_\omega$. So fix a function $f\in \mcH^\alpha_\omega$, then by Lemma \ref{lem:weight} we easily deduce that
$$
\sup_{\genfrac{}{}{0pt}{}{\mcQ \in {\sf StGC}^{k}}{0\leq k\leq 2b}} \big\|\omega^{-1} \mcQ_1(f)\big\|_{L^\infty(\mcM)} \lesssim \big\|\omega^{-1} f\big\|_{L^\infty(\mcM)}.
$$
For the high frequency part, we consider $t\in(0,1]$ and $\mcQ \in {\sf StGC}^{k}$ with $\alpha <k\leq 2b$. Then $\mcQ_t$ has \textit{at least a cancellation of order $1$}, hence
\begin{align*} 
\mcQ_t\big(f\big)(e) & = \mcQ_t\big(f-f(e)\big)(e) \\
& = \int K_{\mcQ_t}\big(e,e'\big) \big(f(e')-f(e)\big) \, \nu(de').
\end{align*}
Due to the kernel support of $\mcQ_t$, the integrated quantity is non-vanishing (and so relevant) only for $\tau \geq \sigma$, with $e=(x,\tau)$ and $e'=(y,\sigma)$. If $ \rho(e,e')\leq 1$, then by definition
$$ 
\big|f(e')-f(e)\big| \leq \omega(e) \rho(e,e')^\alpha \|f\|_{\mcH^\alpha_\omega}
$$
and if $ \rho(e',e)\geq 1$, then by the property of the weight we have
$$ 
\big|f(e')-f(e)\big| \leq \big(\omega(e)+ \omega(e')\big) \big\|\omega^{-1}f\big\|_{L^\infty(\mcM)}. 
$$
Hence
\begin{align*}
\big| \mcQ_t(f)(e)\big| & \lesssim \omega(e) \left\{ \int_{\rho \leq 1} \G_t(e,e') \rho(e,e')^\alpha \, \nu(de') \right. \\
 &+ \qquad \qquad  \left. \int_{\rho \geq 1} \G_t(e,e') \left(1+ \frac{\omega(e')}{\omega(e)}\right) \, \nu(de') \right\} \|f\|_{\mcH^\alpha_\omega} \\
& \lesssim \omega(e) t^\frac{\alpha}{2} \|f\|_{\mcH^\alpha_\omega},
\end{align*}
uniformly in $e \in \mcM$ and $t\in(0,1)$; this concludes the proof of the continuous embedding of $\mcH^\alpha_\omega$ into $\calC^\alpha_\omega$.

\medskip

To prove the converse embedding, let us start by fixing a function $f\in \calC^\alpha_\omega$. The low frequency part of $f$ is easily bounded, using  Calder\'on's reproducing formula
\begin{align*}
\big\|\omega^{-1} f\big\|_{L^\infty(\mcM)} & \lesssim \Big\|\omega^{-1} \mcP^{(1)}_1 f\Big\|_{L^\infty(\mcM)} + \int_0^1 \Big\|\omega^{-1} \mcQ^{(1)}_t f\Big\|_{L^\infty(\mcM)} \, \frac{dt}{t} \\
& \lesssim \|f\|_{\calC^\alpha_\omega},
\end{align*}
since $\alpha>0$. Now fix $e=(x,\tau)$ and $e'=(y,\sigma)$ in $\mcM$, with $\rho := \rho (e,e')\leq 1$ and $\tau \geq \sigma$. We again decompose
$$ 
f = \mcP_1^{(1)}f + \int_0^1 \mcQ_t^{(1)}f \, \frac{dt}{t}.
$$
For $t<\rho^2$, we have
$$ 
\left| \mcQ_t^{(1)}f (e)\right| \lesssim t^\frac{\alpha}{2} \omega(e) \|f\|_{\calC^\alpha_\omega}
$$
and
$$ 
\left| \mcQ_t^{(1)}f (e')\right| \lesssim t^\frac{\alpha}{2} \omega(e') \|f\|_{\calC^\alpha_\omega} \lesssim t^\frac{\alpha}{2} \omega(e) \|f\|_{\calC^\alpha_\omega}  
$$
where we used that the weight is increasing in time and then that $d(x,y)\leq \rho\leq 1$ with the property of the weight. So we may integrate over $t<\rho^2$ and we have
\begin{align*}
\int_0^{\rho^2} \left| \mcQ_t^{(1)}f (e) - \mcQ_t^{(1)}f (e') \right| \, \frac{dt}{t} & \lesssim \left(\int_0^{\rho^2} t^\frac{\alpha}{2} \, \frac{dt}{t} \right) \omega(e) \|f\|_{\calC^\alpha_\omega} \\
& \lesssim \rho^\alpha \omega(e) \|f\|_{\calC^\alpha_\omega}.
\end{align*}
For the low frequency parts, $\mcQ_t^{(1)}$ with $\rho^2 \leq t \leq 1$ or $\mcP^{(1)}_1$, we use that
\begin{align*}
\left| \mcQ_t^{(1)} f(x,\tau) - \mcQ_t^{(1)} f(x,\sigma) \right| 
& \lesssim |\tau-\sigma|^\frac{1}{2} \left(\sup_{\varsigma\in (\sigma,\tau)} \left| \partial_\tau \mcQ_t^{(1)}f(x,\varsigma) \right|\right) \left(\sup_{\varsigma\in (\sigma,\tau)} \left| \mcQ_t^{(1)}f(x,\varsigma) \right|\right) \\
& \lesssim \rho \, \omega(x,\tau) \, t^\frac{\alpha-1}{2} \, \|f\| _{\calC^\alpha_\omega}
\end{align*}
where we used that $\rho \leq 1$ with the fact that the two collections of operators $\big(t \partial_\tau \mcQ_t^{(1)}\big)_{0<t\leq 1}$ and $\big(\mcQ_t^{(1)}\big)_{0<t\leq 1}$ are of type ${\sf StGC}^1$, that is have cancellation of order at least $1$, and that the weight is non-decreasing in time. Similarly we can estimate the variation in space with the assumed finite-increment representation \eqref{finite}, where one considers a local frame field $(X_j)$ in a neighbourhood of a geodesic $(x,y)$ from $x$ to $y$. This gives
\begin{align*}
\left| \mcQ_t^{(1)} f(x,\sigma) - \mcQ_t^{(1)} f(y,\sigma) \right| 
& \lesssim d(x,y)  \sup_{\ z \in (x,y)\,; j}  \left|X_j \mcQ_t^{(1)} f(z,\varsigma)\right| \\
& \lesssim \rho \,\omega(x,\tau) \, t^\frac{\alpha-1}{2} \, \|f\| _{\calC^\alpha_\omega}.
\end{align*}
So we get
\begin{align*}
\int_{\rho^2}^1 \left| \mcQ_t^{(1)}f (e) - \mcQ_t^{(1)}f (e') \right| \, \frac{dt}{t} & \lesssim \rho \left(\int_{\rho^2}^1 t^\frac{\alpha-1}{2} \, \frac{dt}{t} \right) \omega(e) \|f\|_{\calC^\alpha_\omega} \\
& \lesssim \rho^\alpha \, \omega(e) \, \|f\|_{\calC^\alpha_\omega},
\end{align*}
because $\alpha<1$. A similar estimate for $\mcP_1^{(1)}$ ends the proof of continuous embedding of $\calC^\alpha_\omega$ into $\mcH^\alpha_\omega$.
\end{Dem}

\medskip

The next proposition introduces an intermediate space whose unweighted version was first introduced in the setting of paracontrolled calculus in \cite{GIP}, and used in \cite{BB}. To fix notations, and given a space-time weight $\omega$, we denote by $\Big(C^\frac{\alpha}{2}_\tau L^\infty_x\Big)(\omega) = \Big(L^\infty_x C^\frac{\alpha}{2}_\tau\Big)(\omega)$ the set of parabolic distributions such that 
$$
\sup_{x\in M}\;\big\|f(x,\cdot)\big\|_{C^\frac{\alpha}{2}_{\omega(x,\cdot)}(\RR_+)} < \infty.
$$
Also $\Big(L^\infty_\tau C^\alpha_x\Big)(\omega)$ stands for the set of parabolic distributions such that 
$$
\underset{\tau}{\sup}\; \big\|f(\cdot,\tau)\big\|_{C^\alpha_{\omega_\tau}(M)} < \infty.
$$

\begin{prop} \label{prop:LLalpha} {\sf 
Given $\alpha\in (0,2)$ and a space-time weight $\omega$, set
$$ 
\mcE^\alpha_\omega := \Big(C_\tau^\frac{\alpha}{2} L^\infty_x\Big)(\omega) \cap \Big(L^\infty_\tau C^\alpha_x\Big)(\omega).
$$
Then $\mcE^\alpha_\omega$ is continuously embedded into $\calC^\alpha_\omega$. Furthermore, if $\alpha \in (0,1)$, the spaces $\mcE^\alpha_\omega, \calC^\alpha_\omega$ and $\mcH^\alpha_\omega$ are equal, with equivalent norms.  }
\end{prop}

\ssk

\begin{Dem} 
We first  check that $\mcE^\alpha_\omega$ is continuously embedded into $\calC^\alpha_\omega$, and fix for that purpose a function $f\in \mcE^\alpha_\omega$. As done in \cite[Proposition 2.12]{BB}, we know that for all integers $k,j$ with  $k+\frac{j}{2}>\frac{\alpha}{2}$ and every space function $g\in C^\alpha(M)$, we have
$$ 
\Big\| t^\frac{j}{2} V_{J}(tL)^k e^{-tL} g \Big\|_{L^\infty(M)} \lesssim t^\frac{\alpha}{2} \|g\|_{C^\alpha(M)},
$$
for any subset of indices $J$ with $|J|=j$. So consider a generic {\it standard} family $\left( t^\frac{j}{2} V_j(tL)^{\frac{a-j}{2}-k} P_t^{(c)}  \otimes m_t^\star \right)_{0<t\leq 1}$ in ${\sf StGC}^{a}$, with $3\leq a\leq b$, and a smooth function $m$ with vanishing first $k$ moments. If $k=0$ we have seen that we have
$$ 
\Big\| \omega_\tau^{-1} t^\frac{j}{2} V_J (tL)^\frac{a-j}{2} P_t^{(c)} f (\cdot,\tau) \Big\|_{L^\infty(M)} \lesssim t^\frac{\alpha}{2} \big\|f (\tau) \big\|_{C^\alpha_{\omega_\tau}} 
$$
for every $\tau$, so 
$$ 
\Big\| \omega^{-1} t^\frac{j}{2} V_J (tL)^\frac{a-j}{2}P_t^{(c)}  \otimes m_t^\star (f) \Big\|_{L^\infty(\mcM)} \lesssim t^\frac{\alpha}{2} \|f\|_{L^\infty_\tau C^\alpha_x(\omega)}
$$
since $m_t^\star$ is a $L^\infty(\R)$-bounded operator as a convolution with an $L^1$-normalized function.

\ssk

If $k=1$ (or $k\geq 1$), the same reasoning shows that we have 
$$  
\Big\| \omega(x,\cdot)^{-1}  m_t^\star (f) (x,\cdot) \Big\|_{L^\infty(\R_{+})} \lesssim t^\frac{\alpha}{2} \, \big\|f(x,\cdot) \big\|_{C^\frac{\alpha}{2}_{\omega(x,\cdot)}(\R_+)},
$$
for every $x\in M$, since $\frac{\alpha}{2}\in (0,1)$, and $m$ encodes a cancellation at order $1$ in time as it has a vanishing first moment. Hence
$$ 
\Big\| \omega^{-1} t^\frac{j}{2} V_J (tL)^\frac{a-j}{2}P_t^{(c)}  \otimes m_t^\star (f) \Big\|_{L^\infty(\mcM)} \lesssim t^\frac{\alpha}{2} \|f\|_{C_\tau^\frac{\alpha}{2} L^\infty_x(\omega)},
$$
which concludes the proof of the embedding $\mcE^\alpha_\omega \hookrightarrow \calC^\alpha_\omega$. The remainder of the statement is elementary since $\calC^\alpha_\omega = \mcH^\alpha_\omega$ is embedded in $\mcE^\alpha_\omega$. 
\end{Dem}

\ssk

Before turning to the definition of an intertwined pair of parabolic paraproducts, we close this section with two other useful continuity properties involving the H\"older spaces $\calC^\sigma_\omega$.

\begin{prop}  \label{prop:hol}   {\sf 
Given $\alpha\in(0,1)$, a space-time weight $\omega$, some integer $a\geq 0$ and a standard family $\mcP\in {\sf StGC}^a$, there exists a constant $c$ depending only on the weight $\omega$, such that 
$$
\omega(e)^{-1} \Big|\big(\mcP_t f\big)(e) - \big(\mcP_sf\big)(e') \Big| \lesssim \left(s+t+ \rho(e,e')^2\right)^\frac{\alpha}{2}\, e^{c d(x,y)} \big\|f\big\|_{\calC^\alpha_\omega},
$$
uniformly in $s,t\in(0,1]$ and $e= (x,\tau)$ and $e'=(y,\sigma)\in\mcM$, with $\tau\geq \sigma$.  }
\end{prop}

\ssk

\begin{Dem}
We explain in detail the most difficult case corresponding to $\mcP\in {\sf StGC}^0$, so $\mcP$ encodes a priori no cancellation. Then $\mcP_t$ takes the form
$$ 
\mcP_t = P_t^{(c)}  \otimes m_t^\star
$$
for some integer $c\geq 1$ and some smooth function $m$. There is no loss of generality in assuming that $\int m(\tau) \, d\tau$ is equal to $1$, as $\mcP$ is actually an element of ${\sf StGC}^1$ if $m$ has zero mean -- this case is treated at the end of the proof. 

\medskip

In this setting, since $f$ is bounded and continuous, we have the pointwise identity
$$ 
f = \lim_{t\to 0} \mcP_t(f).
$$

\ssk

\textbf{(i)} Consider first the case where $\rho(e,e')\leq 1$,, with $e=(x,\tau)$ and $e'=(y,\sigma)$. Decompose
\begin{align*}
&\omega(e)^{-1} \Big|\big(\mcP_t f\big)(e) - \big(\mcP_sf\big)(e') \Big| \\
&\qquad\leq \omega(e)^{-1} \big|f(e) - f(e') \big| + \omega(e)^{-1} \Big|\big(\mcP_t f\big)(e) - f(e) \Big| \\
&\qquad\qquad+ \omega(e)^{-1} \Big|\big(\mcP_t f\big)(e') - f(e') \Big| \\
&\qquad \lesssim \omega(e)^{-1}\big|f(e) - f(e') \big| + \big\|\omega^{-1} \left(\mcP_t f -f\right) \big\|_{L^\infty(\mcM)} + \big\|\omega^{-1}\left(\mcP_s f -f\right) \big\|_{L^\infty(\mcM)}.
\end{align*}
We have
$$ 
\omega(e)^{-1}\big|f(e) - f(e') \big| \leq \rho(e,e')^{\alpha} \|f\|_{\mcH_\omega^\alpha} \lesssim \rho(e,e')^{\alpha} \|f\|_{\calC_\omega^\alpha}.
$$
For the two other terms, we use that
$$ 
\big\|\omega^{-1} \big(\mcP_t f -f\big) \big\|_{L^\infty(\mcM)} \leq \int_0^t \Big\| \omega^{-1}\, u\,\partial_u \mcP_u f \Big\|_{L^\infty(\mcM)} \, \frac{du}{u},
$$
and note that
$$ 
u\partial_u \mcP_u = Q_u^{(c)} \otimes m_u + P_u^{(c)}  \otimes k_u
$$
with $k(\tau) = \partial_\tau\big(\tau m(\tau)\big)$, is actually the sum of two terms in ${\sf StGC}^{\geq 1}$ since it is clear for the first one and the function $k$ has a vanishing first moment. It follows by definition of the H\"older spaces with $\alpha<1$, that we have
$$  
\big\|\omega^{-1} \left(\mcP_t f -f\right) \big\|_{L^\infty(\mcM)} \lesssim \left(\int_0^t u^\frac{\alpha}{2} \, \frac{du}{u} \right) \|f\|_{\calC_\omega^\alpha} \lesssim t^\frac{\alpha}{2} \|f\|_{\calC_\omega^\alpha}.
$$
A similar estimate holds by replacing $t$ by $s$, which then concludes the proof in this case.

\ssk

\textbf{(ii)} In the case where $\rho(e,e')\geq 1$, we do not use the difference and use condition \eqref{eq:stweight} on the weight $\omega$ to write
$$ 
\omega(x,\tau)^{-1} \leq \omega(x,\sigma)^{-1} \lesssim \omega(y,\sigma)^{-1}.
$$
and obtain as a consequence the estimate
\begin{align*}
\omega(e)^{-1} \Big|\big(\mcP_t f\big)(e) - \big(\mcP_sf\big)(e') \Big| & \leq \omega(e)^{-1} \Big|\big(\mcP_tf\big)(e)\Big| + \omega(e)^{-1} \Big|\big(\mcP_s f\big)(e') \Big| \\
& \leq \big\| \omega^{-1} \mcP_t f  \big\|_{L^\infty(\mcM)} + \omega(x,\sigma)^{-1} \, \Big|\big(\mcP_t f\big)(e') \Big| \\
& \lesssim \big\| \omega^{-1} \mcP_t f  \big\|_{L^\infty(\mcM)} + e^{c d(x,y)} \big\| \omega^{-1} \mcP_s f  \big\|_{L^\infty(\mcM)},
\end{align*}
for some positive constant $c$. Since we know by Lemma \ref{lem:weight} that $\mcP_t$ and $\mcP_s$ are bounded in $L^\infty(\omega)$, we deduce that
\begin{align*}
\omega(e)^{-1} \Big|\big(\mcP_t f\big)(e) - \big(\mcP_sf\big)(e') \Big| &  \lesssim e^{c d(x,y)} \big\| \omega^{-1} f  \big\|_{L^\infty(\mcM)} \\
&  \lesssim e^{c d(x,y)} \|f\|_{\calC^\alpha_\omega},
\end{align*} 
since $\calC^\alpha_\omega \subset L^\infty_\omega$, given that $\alpha>0$. The expected estimate follows from that point.

\medskip

$\bullet$ In the easier situation where $\mcP\in {\sf StGC}^a$ for some integer $a\geq 1$, we can perform the same reasoning and use in addition the fact that 
$$ 
\lim_{t\to 0} \, \mcP_t(f) = 0,
$$ 
which makes the case easier since we do not have to deal with the first term $f(e)-f(e')$.
\end{Dem}

\medskip

With an analogous reasoning (indeed simpler) we may prove the following.

\begin{prop}  \label{prop:hol2}   {\sf 
Given $\alpha\in(-3,0)$, a space-time weight $\omega$ and a standard family $\mcP\in {\sf StGC}^0$, one has 
$$
\| \mcP_t f \|_{L^\infty(\mcM)} \lesssim t^\frac{\alpha}{2} \|f\|_{\calC^\alpha_\omega},
$$
uniformly in $t\in(0,1]$.  }
\end{prop}

\ssk

\begin{Dem} 
The proof follows the same idea as the the proof of Proposition \ref{prop:hol}. Indeed, we use the fact that since $\mcP$ is a standard family then
$$ 
\mcP_tf = \int_t^1 (-s \partial_s \mcP_s)f \, \frac{ds}{s} + \mcP_1 f.
$$
The key point is that $(-s \partial_s \mcP_s)_s$ can be split into a finite sum of families of ${\sf StGC}^{\geq 1}$, which allows us to conclude as previously. 
\end{Dem}

\bigskip

\subsection[\hspace{-0.5cm} Schauder estimates]{Schauder estimates}
\label{SubsectionSchauder}

We provide in this subsection a Schauder estimate for the heat semigroup in the scale of weighted parabolic H\"older spaces. This quantitative regularization effect of the heat semigroup will be instrumental in the proof of the well-posedness of the parabolic Anderson model (PAM) and multiplicative Burgers equations studied in Section \ref{SectionPAM}. Define here formally the linear resolution operator for the heat equation by the formula
\begin{equation}
\label{EqDefnResolution}
\mathscr{L}^{-1}(v)_\tau := \int_0^\tau e^{-(\tau-\sigma)L}v_\sigma \, d\sigma.
\end{equation}
We fix in this section a finite positive time horizon $T$ and consider the space 
$$
\mcM_T := M\times [0,T],
$$
equipped with its parabolic structure. Denote by $L^\infty_T$ the corresponding function space over $[0,T]$. We first state a Schauder estimate that was more or less proved in the unweighted case in \cite{GIP,BB} -- see Lemma A.9 in \cite{GIP} and Proposition 3.10 in \cite{BB}.

\ssk

\begin{prop}
\label{PropFirstSchauder}
Given $\beta\in(-2,0)$ and a space-time weight $\omega$, we have
$$
\big\|\mathscr{L}^{-1}(v)\big\|_{\calC^{\beta+2}_\omega} \lesssim_T \|v\|_{\big(L^\infty_TC^\beta_x\big)(\omega)}.
$$
\end{prop}

\ssk

We shall actually prove a refinement of this continuity estimate in the specific case where $\omega$ has a special structure motivated by the study of the (PAM) and multiplicative Burgers equations done in Section \ref{SectionPAM}. These special weights were first introduced by Hairer and Labb\'e in their study of the (PAM) equation in $\R^2$ and $\R^3$, via regularity structures \cite{HL, HLR3}. Let $o=o_\textrm{ref}$ be the reference point in $M$ fixed and used in the definition of $\calS_o$ at the end of Section \ref{SectionIntroduction}, and set 
\begin{equation}
p_a(x) := \big(1+d(o_\textrm{ref},x)\big)^a, \qquad \varpi(x,\tau) := e^{\kappa \tau} e^{(1+\tau) \big(1+d(o_\textrm{ref},x)\big)},
\label{def:varpi}
\end{equation}
for $0<a<1$ and a positive constant $\kappa$. (The introduction of an extra exponential factor $e^{\kappa \tau}$ in our space-time weight $\varpi$ will allow us to get around  an iterative step in the forthcoming application of the fixed point theorem used to solve the (PAM) and multiplicative Burgers equations, as done in \cite{HL,HLR3}.) For $\tau\geq 0$, we use the notation 
$$\varpi_\tau: \; x\in M \mapsto \varpi(x,\tau)$$
for  the spatial weight.
The space-time weight $\varpi$ satisfies condition \eqref{eq:stweight} on $[0,T]\times M$, uniformly with respect to $\kappa>0$. The above special weights satisfy in addition the following crucial property, already used in \cite{HL,HLR3}. We have
\begin{equation}
p_a(x) \varpi(x,\sigma) \lesssim \kappa^{-\epsilon} (\tau-\sigma)^{-a-\epsilon}\,\varpi(x,\tau),
\label{eq:varpi}
\end{equation}
for every non-negative real number $\epsilon$ small enough, uniformly with respect to $x\in M,\kappa>0$ and $0<\sigma<\tau\leq T$. The next improved Schauder-type continuity estimate shows how one can use the above inequality for the specific weights to compensate a gain on the weight by a loss of regularity.

\ssk

\begin{prop}  \label{prop:schauder}    {\sf
Given $\beta\in{\mathbb R}$, $a\in(0,1)$ and $\epsilon\in[0,1)$ small enough such that $a+e<1$, we have the continuity estimate
$$ 
\big\|\mathscr{L}^{-1} (v)\big\|_{\big(L^\infty_TC^{\beta+2(1-a-\epsilon)}_x\big)(\varpi)} \lesssim \kappa^{-\epsilon} \big\|v\big\|_{\big(L^\infty_TC^\beta_x\big)(\varpi p_a)}.
$$
Moreover if $-2+2(a+\epsilon) <\beta<0$, then 
$$ 
\| \mathscr{L}^{-1} (v) \|_{\calC^{\beta+2-2a-2\epsilon}_\varpi} \lesssim \kappa^{-\epsilon} \big\|v\big\|_{\big(L^\infty_TC^\beta_x\big)(\varpi p_a)}.
$$   }
\end{prop}

\ssk

\begin{Dem} 
Let us first check the regularity in space. So consider an integer $c \geq \frac{|\beta|}{2}+1$ and a parameter $r\in(0,1]$. Then for every fixed time $\tau \in [0,T]$ we have
$$ 
Q_r^{(c)}\big(\mathscr{L}^{-1}(v)_\tau\big) = \int_0^\tau  Q_r^{(c)} e^{-(\tau-\sigma)L}v_\sigma \, d\sigma.
$$
By using the specific property \eqref{eq:varpi} of the weights $p_a$ and $\varpi$, one has
\begin{align*}
\Big\|\varpi_\tau^{-1} Q_r^{(c)} e^{-(\tau-\sigma)L}&v_\sigma\Big\|_{L^\infty(M)} \lesssim \left(\frac{r}{r + \tau-\sigma} \right)^{c} \big\|\varpi_\tau^{-1} Q_{r+\tau-\sigma}^{(c)} v_\sigma\big\|_{L^\infty(M)} \\
& \lesssim \kappa^{-\epsilon} \left(\frac{r}{r + \tau-\sigma}\right)^{c} (r + \tau-\sigma)^{\frac{\beta}{2}} \left(\tau-\sigma\right)^{-a-\epsilon} \big\|v_\sigma\big\|_{C^\beta_{p_a\varpi_\sigma}}.
\end{align*}
So by integrating and using that $c$ is taken large enough, we see that
\begin{align*}
\Big\| \varpi_\tau^{-1} Q_r^{(c)}&\big(\mathscr{L}^{-1}(v)_\tau\big)\Big\|_{L^\infty(M)} \\
&\lesssim \kappa^{-\epsilon} \left\{\int_0^\tau \left(\frac{r}{r+\tau-\sigma}\right)^{c} (r + \tau-\sigma)^{\frac{\beta}{2}} (\tau-\sigma)^{-a-\epsilon}\, d\sigma\right\} \big\|v\big\|_{\big(L^\infty_T C_x^\beta\big)(p_a \varpi)} \\
& \lesssim \kappa^{-\epsilon} \tau^{\frac{\beta}{2}+1-a-\epsilon} \big\|v\big\|_{\big(L^\infty_T C_x^\beta\big)(p_a \varpi)}.
\end{align*}
This holds uniformly in $r\in(0,1]$ and $\tau\in[0,T]$ and so one concludes the proof of the first statement with the global inequality
\begin{align*}
\Big\|\varpi_\tau^{-1} \mathscr{L}^{-1}(v)_\tau\Big\|_{L^\infty(M)} &\lesssim \kappa^{-\epsilon} \left\{\int_0^\tau  (\tau-\sigma)^{-a-\epsilon} d\sigma\right\} \|v\|_{\big(L^\infty_TC_x^\beta\big)(p_a \varpi)} \\
&\lesssim \kappa^{-\epsilon} \tau^ {1-a-\epsilon} \|v\|_{\big(L^\infty_T C_x^\beta\big)(p_a \varpi)}.
\end{align*}

\ssk

For the second statement, we note that for $0\leq \sigma<\tau\leq T$ we have
\begin{align*}
\mathscr{L}^{-1}(v)_\tau - \mathscr{L}^{-1}(v)_\sigma  & = \Big(e^{-(\tau-\sigma)L}-\textrm{Id}\Big)\mathscr{L}^{-1}(v)_\sigma  + \int_\sigma^\tau e^{-(\tau-r)L}v_r \, dr \\ 
                                                     & = \int_0^{\tau-\sigma} Q^{(1)}_{r} \mathscr{L}^{-1}(v)_\sigma \, \frac{dr}{r} + \int_\sigma^\tau e^{-(\tau-r)L}v_r \, dr.
\end{align*}
We have by the previous estimate
\begin{align*} 
\left\|\varpi_\tau^{-1} \int_0^{\tau-\sigma} Q^{(1)}_{r} \mathscr{L}^{-1}(v)_\sigma \, \frac{dr}{r} \right\|_{L^\infty(M)} & \lesssim \kappa^{-\epsilon} \left(\int_0^{\tau-\sigma} r^{\frac{\beta}{2}+1-a-\epsilon} \, \frac{dr}{r} \right) \big\| \mathscr{L}^{-1}(v)_\sigma\big\|_{C_{\varpi_\sigma}^{\beta+2-2a-2\epsilon}} \\ 
 & \lesssim \kappa^{-\epsilon}  (\tau-\sigma)^{\frac{\beta}{2}+1-a-\epsilon} \big\|v\big\|_{\big(L ^\infty_TC_x^{\beta}\big)(\varpi)}
\end{align*}
where we used that $\varpi_\tau \geq \varpi_\sigma$ for $\sigma\leq \tau$. Moreover, since $\beta$ is negative, we also have
{\small
\begin{align*}
 & \left\| \varpi_\tau^{-1} \int_\sigma^\tau e^{-(\tau-r)L}v_r \, dr \right\|_{L^\infty(M)} \\
 & \qquad \lesssim \kappa^{-\epsilon} \int_\sigma^\tau \left(\int_{\tau-r}^1 (\tau-r)^{-a-\epsilon} \Big\| p^{-a}\varpi_r^{-1} Q_{s}^{(1)} v_r\Big\|_{L^\infty(M)} \, \frac{ds}{s} + \Big\| \varpi_r^{-1} e^{-L}\big(v_r\big)\Big\|_{L^\infty(M)}\right) \, dr \\
& \qquad \kappa^{-\epsilon} \lesssim \int_\sigma^\tau  \left(\big\|v_r\big\|_{C_{p_a \varpi_r}^\beta} (\tau-r)^{-a-\epsilon} \int_{\tau-r}^1 s^{\frac{\beta}{2}} \, \frac{ds}{s} + (\tau-r)^{-a-\epsilon} \Big\|e^{-L}\big(v_r\big)\Big\|_{C_{p_a \varpi_r}^\beta}\right) dr \\
& \lesssim \kappa^{-\epsilon} (\tau-\sigma)^{\frac{\beta}{2}+1-a-\epsilon} \big\|v\big\|_{\big(L^\infty_TC_x^\beta\big)(p_a \varpi)},
\end{align*}  }
where we used \eqref{eq:varpi} and $\frac{\beta}{2}+1-a-\epsilon >0$.
\end{Dem}

\medskip

The following result comes as a consequence of the proof, combined with Lemma \ref{lem:weight}; we single it out here for future reference. 

\ssk

\begin{lem} \label{lemm} {\sf 
Let $A$ be a linear operator on $\mcM$ with a kernel pointwisely bounded by $\G_t$ for some $t\in(0,1]$. Then for every $a+\epsilon \in(0,1)$, we have
$$ 
\| A\|_{L^\infty_{\varpi p_a}(\mcM) \to L^\infty_{\varpi}(\mcM)} \lesssim \kappa^{-\epsilon} t^{-a-\epsilon}.
$$   }
\end{lem}

\ssk

Schauder estimates can also be extended to spaces of positive regularity.

\ssk

\begin{prop}  \label{prop:schauder-bis}    {\sf
Given $\beta\in(0,2)$, $a\in(0,1)$ and $\epsilon\in[0,1)$ small enough such that $a+e<1$, we have the continuity estimate
$$ 
\| \mathscr{L}^{-1} (v) \|_{\calC^{\beta+2-2a-2\epsilon}_\varpi} \lesssim \kappa^{-\epsilon} \big\|v\big\|_{\calC^\beta_{\varpi p_a}}.
$$   }
\end{prop}

\begin{Dem}
This follows from Proposition \ref{prop:schauder}. For $v \in \calC^\beta_{\varpi p_a} \subset \big(L^\infty_T C^{\beta}_x)(\varpi p_a\big)$, it is known that $Lv \in \big(L^\infty_T C^{\beta-2}_x\big)(\varpi p_a)$, to which Proposition \ref{prop:schauder} can be applied since $\beta-2<0$. Now use that $\mathscr{L}^{-1}$ and $L$ commute to deduce that $L(\mathscr{L}^{-1} v) \in \calC^{\beta-2a-2\epsilon}_\varpi$, hence $\mathscr{L}^{-1} v \in \big(L^\infty_T C^{\beta+2-2a-2\epsilon}_x\big)(\varpi)$. On the other hand, $\partial_t(\mathscr{L}^{-1} v) = v - L\mathscr{L}^{-1} v$, from which follows that $\partial_t \mathscr{L}^{-1} v \in C^{\beta/2}_T L^\infty_x$, and consequently $\mathscr{L}^{-1} v \in C^{\beta/2+1}_T L^\infty_x$.
\end{Dem}

\medskip

The constraint $\beta<2$ is not relevant. Indeed, by iteration the previous Schauder estimates can be proved for an arbitrary exponent $\beta>0$. 

\bigskip

\section[\hspace{0.7cm}Time-space paraproducts]{Time-space paraproducts}
\label{SectionParaproducts}

We introduce in this section the machinery of paraproducts which we shall use in our analysis of the singular PDEs of Anderson \eqref{EqPAMEq} and Burgers \eqref{EqBurgersEq}. In the classical setting of analysis on the torus, the elementary definition of a paraproduct given in Section \ref{SectionParacontrolledNutshell} in terms of Fourier analysis should make convincing, for those who are not familiar with this tool, the fact that $\Pi_f(g)$ is a kind of "modulation" of $g$, insofar as each mode $g_j$ of $g$, in its Paley-Littlewood decomposition, is modulated by a signal which oscillates at frequencies much smaller -- the finite sum $\sum_{0\leq i\leq j-2} f_i$. So it makes sense to talk of a distribution/function of the form $\Pi_f(g)$ as a distribution/function that "locally looks like" $g$. This is exactly how we shall use paraproducts, as a tool that can be used to provide some kind of Taylor expansion of a distribution/function, in terms of some other 'model' distributions/functions. This will be used crucially to bypass the ill-posed character of some operations involved in the (PAM) and Burgers equations, along the line of what was written in Section \ref{SectionParacontrolledNutshell}. 

\ssk

Working in a geometric setting where Fourier analysis does not make sense, we shall define our paraproduct entirely in terms of the semigroup generated by the operator $\mathscr{L} = \partial_t + L$ on the parabolic space. The definition of a paraproduct comes  together with the definition of a resonant operator $\Pi(\cdot,\cdot)$, tailor-made to provide the decomposition
$$
fg = \Pi_f(g) + \Pi(f,g) + \Pi_g(f)
$$
of the product operation, and with $\Pi_f(g)$ and $\Pi(f,g)$ with good continuity properties in terms of $f$ and $g$ in the scales of H\"older spaces. Such a construction was already done in our previous work \cite{BB}, where the generic form of the operator $L$, given by its first order carr\'e du champ operator, imposed some restrictions on the range of the method and allowed only a first order machinery to be set up. The fact that we work here with an operator $L$ in H\"ormander form will allow us to set up a higher order expansion setting. We will use this for the description of the space in which to make sense of the two singular PDEs we want to analyse. However, this a priori useful setting is in direct conflict with one of the main technical tools introduced by Gubinelli, Imkeller and Perkowski in their seminal work \cite{GIP}.

\ssk

The case is easier to explain on the example of the (PAM) equation. A solution to that equation is formally given as a fixed point of the map 
$$
\Phi : u\mapsto e^{-\cdot \,L} u_0 + \mathscr{L}^{-1}(u \zeta),
$$
for which we shall need $u$ to be a priori controlled by $Z := \mathscr{L}^{-1}(\zeta)$, to make sense of the product $u \zeta$ -- more will actually be required, but let us stick to this simplified picture here; so the map $\Phi$ will eventually be defined on a space of distributions controlled by $Z$, such as defined in Section \ref{SectionParacontrolledNutshell}, where it will be shown to be a contraction. At a heuristic level, for a distribution $(u,u')$ controlled by $Z$, the product $u \zeta$ will be given by a formula of the form
$$
u\zeta = \Pi_{u}(\zeta) + (\cdots).
$$
To analyse the term $\mathscr{L}^{-1}(u \zeta)$, and recalling that $Z := \mathscr{L}^{-1}(\zeta)$, it is thus very tempting to write
\begin{equation*}
\mathscr{L}^{-1}\big(\Pi_{u}(\zeta)\big) = \Pi_{u}\big(Z\big) + \big[\mathscr{L}^{-1},\Pi_{u}\big](\zeta) + (\cdots)
\end{equation*}
and work with the commutator $\big[\mathscr{L}^{-1},\Pi_{u}\big]$. This is what was done in \cite{GIP,BB} to study the $2$-dimensional (PAM) equation on the torus and more general settings; and it somehow leads to a non-natural choice of function space for the remainder $f^\sharp$ of a paracontrolled distribution in a space-time setting. Unfortunately, we have little information on this commutator, except from the fact that it is a regularizing operator with a quantifiable regularizing effect -- it was first proved in \cite{GIP} in their Fourier setting. This sole information happens  to be insufficient to push the analysis of the (PAM) or Burgers equations far enough in a $3$-dimensional setting. As a way out of this problem, we introduce another paraproduct $\widetilde \Pi_{v}(\cdot)$, tailor-made to deal with that problem, and intertwined to $\Pi_{v}(\cdot)$ via $\mathscr{L}^{-1}$, that is
$$
\mathscr{L}^{-1}\circ\Pi_{v} = \widetilde \Pi_{v}\circ\mathscr{L}^{-1};
$$
so $\widetilde{\Pi}$ is formally the $\Pi$ operator seen in a different basis 
$$
\widetilde{\Pi} = \mathscr{L}^{-1}\circ\Pi\circ\mathscr{L}.
$$ 
We show in Section \ref{SubsectionIntertwined} that $\Pi$ and $\widetilde\Pi$ have the same analytic properties. In particular, if $f\in L^\infty_TC^\alpha_x$ with $-2<\alpha<0$, the Schauder estimate proved in proposition \ref{prop:schauder} shows that $\widetilde\Pi_v\big(\mathscr{L}^{-1} f\big)$ is an element of the parabolic H\"older space $\calC^{\alpha+2}$. In the end, we shall be working with an ansatz for the solution space of the 3-dimensional (PAM) equation given by distributions/functions of the form
$$
u = \widetilde\Pi_{u'}(Z) + (\cdots).
$$ 

\ssk

The introduction of semigroup methods for the definition and study of paraproducts is relatively new; we refer the reader to different recent works where such paraproducts have been used and studied \cite{B-T1,BS,BBR,BCF2, BB}.

\bigskip

\subsection[\hspace{-0.5cm} Intertwined paraproducts]{Intertwined paraproducts}
\label{SubsectionIntertwined}

We introduce in this section a pair of intertwined paraproducts that will be used to analyze the a priori ill-posed terms in the right hand side of the parabolic Anderson model equation and multiplicative Burgers system in the next section. We follow here for that purpose the semigroup approach developed in \cite{BB}, based on the pointwise Calder\'on reproducing formula
$$
f = \int_0^1\mcQ^{(b)}_tf \,\frac{dt}{t}+ \mcP_1^{(b)}f,
$$
where $f$ is  a bounded and continuous function. This formula says nothing else than the fact that 
$$
\underset{t\downarrow 0}{\lim}\;\mcP^{(b)}_t = \textrm{Id}.
$$
(This is a direct consequence of the fact that the operator $\varphi_t^\star$ tends to the identity operator, since $\varphi$ has unit integral.) We can thus write formally for two continuous and bounded functions $f,g$
{\small
\begin{equation}
\label{EqDecompCalderonFG}
\begin{split}
&fg = \lim_{t\to 0} \; \mcP^{(b)}_t \left( \mcP_t^{(b)} f \cdot \mcP_t^{(b)} g \right) = -\int_0^1 t \partial_t \,\left\{ \mcP_t^{(b)} \left( \mcP_t^{(b)} f \cdot \mcP_t^{(b)} g \right)\right\} \, \frac{dt}{t} + \Delta_{-1}(f,g)  \\
& = \int_0^\infty  \left\{\mcP_t^{(b)}\left( \mcQ_t^{(b)} f\cdot \mcP_t^{(b)} g\right) + \mcP_t^{(b)} \left(\mcP_t^{(b)} f\cdot \mcQ_t^{(b)} g\right) + \mcQ_t^{(b)} \left(\mcP_t^{(b)} f\cdot \mcP_t^{(b)} g\right) \right\}\, \frac{dt}{t} + \Delta_{-1}(f,g),
\end{split}
\end{equation}   }
where
$$
\Delta_{-1}(f,g):= \mcP_1^{(b)}\left(\mcP_1^{(b)} f \cdot \mcP_1^{(b)} g \right) 
$$
stands for the ``low-frequency part'' of the product of $f$ and $g$. This decomposition corresponds to an extension of Bony's well-known paraproduct decomposition \cite{Bony} to our setting given by a semigroup.

\medskip

The integral exponent $b$ has not been chosen so far. Choose it here even and no smaller than $6$. Using iteratively the Leibniz rule for the differentiation operators $V_i$ or $\partial_\tau$, generically denoted $D$,
$$
D(\phi_1)\phi_2 = D(\phi_1\cdot\phi_2) - \phi_1\cdot D(\phi_2),
$$
we see that $\mcP_t^{(b)}\left(\mcQ_t^{(b)} f\cdot \mcP_t^{(b)} g\right)$ can be decomposed as a finite sum of terms taking the form
$$ 
\mcA^{I,J}_{k,\ell}(f,g) := \mcP_t^{(b)} \Big(t^{\frac{|I|}{2}+k} V_I \partial_\tau^k\Big) \left( \mcS_t^{(b/2)} f \cdot \big(t^{\frac{|J|}{2}+\ell} V_J\partial_\tau^\ell\big) \mcP_t^{(b)} g\right) 
$$  
where $\mcS^{(b/2)} \in {\sf StGC}^{\frac{b}{2}}$ and the tuples $I,J$ and integers $k,\ell$ satisfy the constraint 
$$
\frac{|I|+|J|}{2} + k+\ell = \frac{b}{2}.
$$ 
Denote by $\mcI_b$ the set of all such $(I,J,k,\ell)$. We then have the identity
$$  
\int_0^1 \mcP_t^{(b)}\left(\mcQ_t^{(b)} f\cdot \mcP_t^{(b)} g\right) \, \frac{dt}{t} = \sum_{\mcI_b} \, a^{I,J}_{k,\ell} \int_0^1 \mcA^{I,J}_{k,\ell}(f,g) \, \frac{dt}{t},
$$
for some coefficients $a^{I,J}_{k,\ell}$. Similarly, we have
$$  
\int_0^1 \mcQ_t^{(b)}\left(\mcP_t^{(b)} f\cdot \mcP_t^{(b)} g\right) \, \frac{dt}{t} = \sum_{\mcI_b} \, b^{I,J}_{k,\ell} \int_0^1 \mcB^{I,J}_{k,\ell}(f,g) \, \frac{dt}{t},
$$
with $\mcB^{I,J}_{k,\ell}(f,g)$ of the form
$$ 
\mcB^{I,J}_{k,\ell}(f,g) := \mcS_t^{\big(\frac{b}{2}\big)} \left( \Big\{\big(t^{\frac{|I|}{2}+k} V_I\partial_\tau^k\big) \mcP_t^{(b)} f\Big\} \cdot \Big\{\big(t^{\frac{|J|}{2}+\ell} V_J\partial_\tau^\ell\big)\mcP_t^{(b)} g\Big\}\right),
$$
for some coefficients $b^{I,J}_{k,\ell}$. So we have at the end the decomposition
$$ 
fg = \sum_{\mcI_b} \, a^{I,J}_{k,\ell} \int_0^1\Big(\mcA^{I,J}_{k,\ell}(f,g) + \mcA^{I,J}_{k,\ell}(g,f)\Big) \, \frac{dt}{t} + \sum_{\mcI_b} \, b^{I,J}_{k,\ell} \int_0^1 \mcB^{I,J}_{k,\ell}(f,g) \, \frac{dt}{t},
$$
which leads us to the following definition.

\begin{defn*} \label{def:paraproduit} {\sf
Given $f\in \bigcup_{s\in (0,1)} \calC^s$ and $g\in L^\infty(\mcM)$, we define the paraproduct $\Pi^{(b)}_g(f)$ by the formula
\begin{align*}
\Pi^{(b)}_g(f) &:= \int_0^1 \Bigg\{  \sum_{\mcI_b ; \frac{|I|}{2} + k > \frac{b}{4}} \, a^{I,J}_{k,\ell} \, \mcA^{I,J}_{k,\ell}(f,g) +  \sum_{\mcI_b ; \frac{|I|}{2}+k> \frac{b}{4}} \, b^{I,J}_{k,\ell} \, \mcB^{I,J}_{k,\ell}(f,g) \Bigg\} \, \frac{dt}{t},
\end{align*}
and the resonant term $\Pi^{(b)}(f,g)$ by the formula
\begin{align*} 
\int_0^1 \Bigg\{   \sum_{\mcI_b ; \frac{|I|}{2} + k \leq  \frac{b}{4}} \, a^{I,J}_{k,\ell} \Big( \mcA^{I,J}_{k,\ell}(f,g) + \mcA^{I,J}_{k,\ell}(g,f)\Big)
+ \sum_{\mcI_b ; \frac{|I|}{2}+k = \frac{|J|}{2}+\ell = \frac{b}{4}} \, b^{I,J}_{k,\ell} \, \mcB^{I,J}_{k,\ell}(f,g) \Bigg\} \, \frac{dt}{t}.
\end{align*}   }
\end{defn*}

\medskip

With these notations, Calder\'on's formula becomes
$$ 
fg= \Pi^{(b)}_g(f)+\Pi^{(b)}_f(g)+\Pi^{(b)}(f,g)+\Delta_{-1}(f,g)
$$
with the ``low-frequency part''
$$ 
\Delta_{-1}(f,g):= \mcP_1^{(b)}\left(\mcP_1^{(b)} f \cdot \mcP_1^{(b)} g \right).
$$

\medskip

If $b$ is chosen large enough, then all of the operators involved in paraproducts and resonant term have a kernel pointwisely bounded by a kernel $\G_t$ at the right scaling. Moreover, 
\begin{enumerate}
   \item the paraproduct term $\Pi^{(b)}_g(f)$ is a finite linear combination of operators of the form
$$ 
\int_0^1 \mcQ^{1\bullet}_t\Big( \mcQ^{2}_t f \cdot \mcP^1_t g\Big) \, \frac{dt}{t}
$$ 
with $\mcQ^1,\mcQ^2\in {\sf StGC}^\frac{b}{4}$, and $\mcP^1\in {\sf StGC}^{[0,2b]}$.  \vspace{0.1cm}
   
   \item the resonant term $\Pi^{(b)}(f,g)$ is a finite linear combination of operators of the form
$$ 
\int_0^1 \mcP^1_t\Big( \mcQ^{1}_t f \cdot \mcQ^{2}_t g\Big) \, \frac{dt}{t}
$$ 
with $\mcQ^1, \mcQ^2\in {\sf StGC}^\frac{b}{4}$ and $\mcP^1\in {\sf StGC}^{[0,2b]}$.
\end{enumerate}
Note that since the operators $\mcQ^\bullet$ and $\mcP^1_t$ are of the type $\mcQ_t^{(c)}$, $\mcP_t^{(c)}$ or a $\mcP_t^{(c)} V_I $, they can easily be composed on the left with another operator $\mcQ_r^{(d)}$; this will simplify the analysis of the paraproduct and resonant terms in the parabolic H\"older spaces. Note also that $\Pi^{(b)}_f({\bf 1}) = \Pi^{(b)}(f,{\bf 1}) = 0$, and that we have the identity
$$ 
\Pi^{(b)}_{{\bf 1}}(f) = f-\mcP_{\bf 1}^{(b)}\mcP_{\bf 1}^{(b)} f,
$$
as a consequence of our choice of the normalizing constant. Therefore the paraproduct with the constant function ${\bf 1}$ is equal to the identity operator, up to the strongly regularizing operator $\mcP_{\bf 1}^{(b)}\mcP_{\bf 1}^{(b)}$.

\bigskip

One can prove the following continuity estimates in exactly the same way as in \cite{BB}. Note first that if $\omega_1,\omega_2$ are two space-time weights, then $\omega:=\omega_1 \omega_2$ is also a space-time weight.

\ssk

\begin{prop}
\label{PropRegularityParaproduct}
{\sf 
Let $\omega_1,\omega_2$ be two space-time weights, and set $\omega := \omega_1\,\omega_2$.
\begin{itemize}
   \item[{\bf (a)}] For every $\alpha,\beta\in\RR$ and every positive regularity exponent $\gamma$, we have
   $$
   \big\| \Delta_{-1}(f,g)\big\|_{\calC^\gamma_\omega} \lesssim  \|f\|_{\calC^\alpha_{\omega_1}} \|g\|_{\calC^\beta_{\omega_2}}
   $$  
   for every $f\in \calC^\alpha_{\omega_1}$ and $g\in\calC^\beta_{\omega_2}$. \vspace{0.2cm}
   \item[{\bf (b)}] For every $\alpha\in (-3,3)$ and $f\in \calC^\alpha_{\omega_1}$, we have
   \begin{equation*} 
      \Big\| \Pi^{(b)}_g(f)\Big\|_{\calC^\alpha_\omega} \lesssim \big\|\omega_2^{-1} g \big\|_\infty \|f\|_{\calC^\alpha_{\omega_1}}
      \end{equation*}   
   for every $g\in L^\infty(\omega_2^{-1})$, and 
   \begin{equation*} 
      \Big\| \Pi^{(b)}_g(f)\Big\|_{\calC^{\alpha+\beta}_\omega} \lesssim \|g\|_{\calC^\beta_{\omega_2}} \|f\|_{\calC^\alpha_{\omega_1}}
      \end{equation*}
   for every $g\in\calC^\beta_{\omega_2}$ with $\beta<0$ and $\alpha+\beta\in (-3,3)$.  \vspace{0.2cm}
   \item[{\bf (c)}] For every $\alpha,\beta\in (-\infty,3)$ with $\alpha+\beta>0$, we have the continuity estimate
$$ 
\Big\| \Pi^{(b)}(f,g) \Big\|_{\calC^{\alpha+\beta}_\omega} \lesssim \|f\|_{\calC^\alpha_{\omega_1}} \|g\|_{\calC^\beta_{\omega_2}}
$$  
for every $f\in\calC^\alpha_{\omega_1}$ and $g\in \calC^\beta_{\omega_2}$.
\end{itemize}   }
\end{prop}

\ssk

The range $(-3,3)$ for $\alpha$ (or $\alpha+\beta$) is due to the fact that all the operators involving a cancellation used in this estimate satisfy a cancellation of order at least $\nu+10 > 3$. We simply write $3$ in the above statement, which will be sufficient for our purpose. We proved similar regularity estimates for the paraproduct introduced in \cite{BB}, with a range for $\alpha$ limited to $(-2,1)$. This difference reflects the fact that the class of operators $L$ considered in \cite{BB}, characterized by the first order carr\'e du champ operators, is more general than the class of H\"ormander form operators considered in the present work, and allows only for a first order calculus.

\medskip

These regularity estimates can be refined if one uses the specific weights $\varpi$ and $p_a\varpi$ introduced in Subsection \ref{SubsectionSchauder}.

\ssk

\begin{prop} \label{prop:para-bis}  {\sf 
For every $\alpha\in (-3,3)$ and $a,\epsilon\in (0,1)$ with $\alpha-a-\epsilon\in(-3,3)$ and $f\in \calC^\alpha_{p_a}$, we have
\begin{itemize}
 \item for every $g\in L^\infty_\varpi$
\begin{equation*}  
\Big\| \Pi^{(b)}_g(f)\Big\|_{\calC^{\alpha-a-\epsilon}_\varpi} \lesssim \kappa^{-\epsilon} \big\|\varpi^{-1} g\big\|_\infty \|f\|_{\calC^\alpha_{p_a}};
\end{equation*}
 \item for every $g\in\calC^\beta_{\varpi}$ with $\beta<0$ and $\alpha+\beta-a\in (-3,3)$  
\begin{equation*} 
\Big\| \Pi^{(b)}_g(f)\Big\|_{\calC^{\alpha+\beta-2(a+\epsilon)}_\varpi} \lesssim \kappa^{-\epsilon} \|g\|_{\calC^\beta_{\varpi}} \|f\|_{\calC^\alpha_{p_a}}. 
\end{equation*}
\end{itemize}  }
\end{prop}

\ssk

The proof of this result is done along exactly the same lines as the proof of Proposition \ref{PropRegularityParaproduct}, using as an additional ingredient the elementary Lemma \ref{lemm}.

\bigskip

We shall use the above paraproduct in our study of the parabolic Anderson model equation, and multiplicative Burgers system, to give sense to the a priori undefined products $u\zeta$ and $\textrm{M}_\zeta u$ of a $\calC^\alpha$ function $u$ on $\mcM$ with a $\calC^{\alpha-2}$ distribution $\zeta$ on $\mcM$, while $2\alpha-2\leq 0$. Our higher order paracontrolled setting is developed for that purpose. As said above, and roughly speaking, we shall solve the Anderson equation 
$$ 
(\partial_\tau + L) u = u\zeta
$$
by finding a fixed point to the map $\Phi(u) = e^{-\cdot\,L} u_0 + \mathscr{L}^{-1}(u \zeta)$. We would like to set for that purpose a setting where the product $u\zeta$ can be decomposed as a sum of the form
$$
u\zeta = \sum_{i=1}^3 \Pi^{(b)}_{u_i}(Y_i) + (\cdots),
$$
for some remainder term $(\cdots)$. We would then have
$$
\mathscr{L}^{-1}(u \zeta) = \sum_{i=1}^3 \mathscr{L}^{-1}\Big(\Pi^{(b)}_{u_i}(Y_i)\Big) + (\cdots),
$$
which we would like to write in the form 
$$
\mathscr{L}^{-1}(u \zeta) = \sum_{i=1}^3 \Pi^{(b)}_{u_i}\big(\mathscr{L}^{-1}(Y_i)\big) + (\cdots),
$$
commuting the resolution operator $\mathscr{L}^{-1}$ with the paraproduct. The commutation is not perfect though and only holds up to a correction term involving the regularizing commutator operator $\big[\mathscr{L}^{-1},\Pi_g(\cdot)\big]$, whose regularizing effect happens to be too limited for our purposes. This motivates us to introduce the following operator.

\ssk

\begin{defn*}  {\sf 
We define a modified paraproduct $\widetilde\Pi^{(b)}$ setting 
$$
\widetilde \Pi^{(b)}_g( f ) := \mathscr{L}^{-1} \Big( \Pi^{(b)}_g \big( \LL f \big)\Big).
$$    }
\end{defn*}

\ssk

The next proposition shows that if one chooses the parameters $\ell_1$ that appears in the reference kernels $\G_t$, and the exponent $b$ that appears in the definition of the paraproduct, both large enough, then the modified paraproduct $\widetilde\Pi^{(b)}_g(\cdot)$ has the same algebraic/analytic properties as $\Pi^{(b)}_g(\cdot)$.

\ssk

\begin{prop} \label{prop:modifiedparaproduit}  {\sf 
If the ambient space $M$ is bounded, then for a large enough choice of constants $\ell_1$ and $b$, the modified paraproduct $\widetilde \Pi_g(f)$ is a finite linear combination of operators of the form
$$ 
\int_0^1 \mcQ^{1\bullet}_t\Big( \mcQ^{2}_t f \cdot \mcP^1_t g\Big) \, \frac{dt}{t}
$$ 
with $\mcQ^1 \in {\sf GC}^{\frac{b}{8}-2}$, $\mcQ^2\in {\sf StGC}^\frac{b}{4}$ and $\mcP^1\in {\sf StGC}$.   

\ssk 

If the space $M$ is unbounded, then the result still holds on the parabolic space $[0,T] \times M$ for every $T>0$, with implicit constants depending on $T$.   }
\end{prop}

\ssk

The operators $\mcQ^1_t$ that appears in the decomposition of $\Pi_g(f)$ are elements of ${\sf StGC}^{[0,2b]}$, while the operators $\mcQ_t^1$ that appear in the decomposition of $\widetilde\Pi_g(f)$ are mere elements of ${\sf GC}^{\frac{b}{8}-2}$. 

\medskip

\begin{Dem}
Given the structure of $\Pi_g^{(b)}(\cdot)$ as a sum of terms of the form 
$$ 
\int_0^1 \mcQ^{1\bullet}_t \Big( \mcQ^{2}_t (\cdot) . \mcP^1_t g \Big) \, \frac{dt}{t}
$$
with $\mcP^1\in {\sf StGC}$ and $\mcQ^1,\mcQ^2\in {\sf StGC}^\frac{b}{4}$, it suffices to look at 
$$ 
\int_0^1  (t^{-1} \mathscr{L}^{-1}) \mcQ^{1\bullet}_t \Big( \mcQ^{2}_t (t \LL) (\cdot) . \mcP^1_t g \Big) \, \frac{dt}{t}.
$$
We have $\mcP^1\in {\sf StGC}^{[0,2b]}$, and it is easy to check that $\big(\mcQ_t^{2} (t\LL)\big)_{0<t\leq 1}$ also belongs to ${\sf StGC}^{\frac{b}{4}+2} \subset {\sf StGC}^{\frac{b}{4}}$. Insofar as 
$$
\mathscr{L}^{-1} \mcQ^{1\bullet}_t = \big(\mcQ^1_t\mathscr{L}^{-1}\big)^\bullet,
$$
we are left with proving that the family $\widetilde{\mcQ^1} := \big(\mcQ^1_t t^{-1} \mathscr{L}^{-1}\big)_{0<t\leq 1}$ belongs to ${\sf GC}^{\frac{b}{8}-2}$, with $\mcQ^1$ essentially given here by
$$ 
\mcQ^1_t = \Big(t^{\frac{|I|}{2}+k} V_I \partial_\tau^k\Big) \mcP_t^{(b)}
$$ 
with $\frac{|I|}{2}+k > \frac{b}{4}$. Note in particular that we have either $|I|\geq \frac{b}{4}$ or $k\geq \frac{b}{8}$. We check in the first two steps of the proof that $\widetilde {\mathcal Q} \in {\sf G}$ in both cases provided $b$ is chosen big enough. The third step is dedicated to proving that $\widetilde {\mathcal Q}^1 \in {\sf GC}^{\frac{b}{4}-1}$.

\medskip

\noindent {\bf Step 1.} Assume here that $|I|\geq \frac{b}{4}$. The kernel $K$ of $\mcQ^1_t \circ (t^{-1} \mathscr{L}^{-1})$ is given by
\begin{equation}
\label{EqKernel}
K\big((x,\tau), (y,\sigma)\big) = \int_{\sigma}^\infty K_{t^{|I|/2} V_I P_{t}^{(b)} e^{-(\lambda-\sigma)L}}(x,y) (t\partial_\tau)^k \varphi_t(\tau-\lambda) \, \frac{d\lambda}{t^2}.
\end{equation}
So by the Gaussian estimates of the operator $t^\frac{|I|}{2} V_I P_{t}^{(b)} e^{-(\lambda-\sigma)L}$ at scale $\max(t,\lambda-\sigma)^{1/2}$, and since $|I|\geq \frac{b}{4}$, we deduce that
\begin{align*}
\left|K_{t^\frac{|I|}{2} V_I P_{t}^{(b)} e^{-(\lambda-\sigma)L}}(x,y)\right| & \lesssim \left(\frac{t}{t+\lambda-\sigma}\right)^\frac{|I|}{2} \G_{t+\lambda-\sigma}(x,y) \\
& \lesssim \left(\frac{t}{t+\lambda-\sigma}\right)^{\frac{b}{8} - \frac{\nu}{2}} \mu(B(x,\sqrt{t}))^{-1} \left(1+\frac{d(x,y)^2}{t+\lambda-\sigma}\right)^{-\ell_1} \\
& \lesssim \left(\frac{t}{t+\lambda-\sigma}\right)^{\frac{b}{8} - \frac{\nu}{2} - \ell_1} \mu(B(x,\sqrt{t}))^{-1} \left(1+\frac{d(x,y)^2}{t}\right)^{-\ell_1}
\end{align*}
\textit{if $b$ is chosen large enough for $\frac{b}{8} - \frac{\nu}{2} - \ell_1$ to be non-negative}. Using the smoothness of $\varphi$ we then deduce that $\Big| K \big((x,\tau), (y,\sigma) \big) \Big|$ is bounded above by 
\begin{align*}
&\mu(B(x,\sqrt{t}))^{-1} \left(1+\frac{d(x,y)^2}{t}\right)^{-\ell_1} \int_{\sigma}^\infty \left(\frac{t}{t+\lambda-\sigma}\right)^{\frac{b}{8} - \frac{\nu}{2} - \ell_1} \left(1+\frac{\tau-\lambda}{t}\right)^{-\ell_1} \, \frac{d\lambda}{t^2} \\
& \lesssim \frac{1}{t\,\mu\big(B(x,\sqrt{t})\big)} \left(1+\frac{d(x,y)^2}{t}\right)^{-\ell_1} \left(1+\frac{|\tau-\sigma|}{t}\right)^{-\ell_1}.
\end{align*}
So we get the upper bound
\begin{equation} 
\label{eq:aa1}
\Big| K \big((x,\tau), (y,\sigma) \big) \Big| \lesssim \nu\big(B_\mcM\big((x,\tau),\sqrt{t}\big)\big)^{-1} \left(1+\frac{d(x,y)^2+|\tau-\sigma|}{t}\right)^{-\ell_1}.
\end{equation}
If $d(x,y)\leq 1$, this is exactly the desired estimate. If $d(x,y)\geq 1$ and one works on a finite time interval $[0,T]$ then we keep the information that $|\lambda-\sigma|\leq T$ and so the exponentially decreasing term in the Gaussian kernel on the spatial variable allows us to keep in all the previous computations an extra coefficient of the form
$$ 
\mu(B_M(x,1))^{-1} e^{-c \frac{d(x,y)^2}{1+T}}
$$ 
which is exactly the decay required in the definition of the class ${\sf G}$. 

\medskip

\noindent {\bf Step 2.} Assume now that $k\geq \frac{b}{8}$. We work with the above formula for the kernel $K$ and use the cancellation effect in the time variable by integrating by parts in $\lambda$ for transporting the cancellation from time to space variable. So starting from formula \eqref{EqKernel}, the ``boundary term" in the integration by parts
$$ 
K_{t^\frac{|I|}{2} V_I P_{t}^{(b)} e^{-(\lambda-\sigma)L}}(x,y) (t\partial_\tau)^{k-1} \varphi_t(\tau-\lambda) 
$$
is vanishing for $\lambda\to \infty$, and equal to 
$$ 
K_{t^\frac{|I|}{2} V_I P_{t}^{(b)}}(x,y) (t\partial_\tau)^{k-1} \varphi_t(\tau-\sigma)
$$
for $\lambda=\sigma$. The latter term satisfies estimate \eqref{eq:aa1}. So up to a term denoted by $(\checkmark)$, bounded as desired, we see that $K \big((x,\tau), (y,\sigma)\big)$ is equal to 
$$ 
(\checkmark) +  \int_{\sigma}^\infty  K_{t^{\frac{|I|}{2}+1} V_I P_{t}^{(b)} Le^{-(\lambda-\sigma)L}}(x,y) (-t\partial_\lambda)^{k-1} \varphi_t(\tau-\lambda) \, \frac{d\lambda}{t^2},
$$
where we used that by analyticity of $L$ in $L^1(M)$ 
$$ 
\partial_\lambda e^{-(\lambda-\sigma)L} = -L e^{-(\lambda-\sigma)L}.
$$ 
Doing $k$ integration by parts provides an identity of the form
$$ 
K\big((x,\tau), (y,\sigma)\big) = (\checkmark) +  \int_{\sigma}^\infty  K_{t^{\frac{|I|}{2}+k} V_I P_{t}^{(b)} L^ke^{-(\lambda-\sigma)L}}(x,y)  \varphi_t(\tau-\lambda) \, \frac{d\lambda}{t^2},
$$
where $(\checkmark)$ stands for a term with \eqref{eq:aa1} as an upper bound. This procedure leaves us with a kernel which has an order of cancellation at least $\frac{b}{8}$ in space; we can then repeat the analysis of Step 1 to conclude.

\medskip

\noindent {\bf Step 3.} The proof that $\widetilde{\mcQ}^1$ actually belongs to ${\sf GC}^{\frac{b}{8}-2}$ is very similar, with details largely left to the reader. The above two steps make it clear that the study of $\widetilde\mcQ^1$ reduces to the study of operators with a form similar to that of the elements of ${\sf StGC}^{[0,2b]}$. We have provided all the details in Proposition \ref{prop:SO} of how one can estimate the composition between such operators and obtain an extra factor encoding the cancellation property. The cancellation result on $\widetilde\mcQ^1$ comes by combining the arguments of Proposition \ref{prop:SO} with the two last steps.

\ssk

Let us give some details for the particular case where the family $\mcQ$ belongs to ${\sf StGC}^{a}$ for some $a\geq \frac{b}{8}-1$ and commutes with $\mathscr{L}^{-1}$; this covers in particular the case where $\mcQ$ is built in space only with the operator $L$ with no extra $V_i$ involved. Let us then take $s,t\in (0,1)$ and consider the kernel of  the operator $\widetilde{\mcQ}^1_t Q_s^\bullet$. Note first that
\begin{align*}
\widetilde{\mcQ}^1_t Q_s^\bullet  & = \Big(\mcQ^1_t \mcQ_s^\bullet \Big)\circ \big( t^{-1} \mathscr{L}^{-1}\big) \\
& = \frac{t+s}{t} \, \Big(\mcQ^1_t \mcQ_s^\bullet \Big) \, (t+s)^{-1} \mathscr{L}^{-1}.
\end{align*}
Since $\mcQ^1 \in \mcO^{\frac{b}{4}}$, we know that $\mcQ^1_t \mcQ_s^\bullet$ is an operator with a kernel with decay at scale $(t+s)^\frac{1}{2}$ with an extra factor $\left(\frac{st}{(t+s)^2}\right)^\frac{b}{8}$. We may also consider that 
$$ 
\mcQ^1_t Q_s^\bullet = \left(\frac{st}{(t+s)^2}\right)^\frac{b}{16} \mcQ^2_{t+s} (t+s)^{-1} \mathscr{L}^{-1} 
$$
for some operator $\mcQ^2_{t+s}$ having $\frac{b}{8}$-order of cancellation and a kernel with decay at scale $\sqrt{s+t}$. So by what we did in the two first steps we also obtain that $\mcQ^2_{t+s} (t+s)^{-1} \mathscr{L}^{-1}$ has a kernel with decay at scale $(t+s)^\frac{1}{2}$, for a large enough choice of $b$. (Indeed, note that $\mcQ^2$ is very similar to the operators studied in the two first steps: easily analyzed as a function of  the space-variable, while, as far as the time-variable is concerned, the composition of convolution preserves the main properties needed on the functions -- vanishing moments.) At the end, we conclude that 
$$ 
\widetilde{\mcQ}^1_t Q_s^\bullet = \left(\frac{st}{(t+s)^2}\right)^{\frac{b}{16}-1} \mcQ^2_{t+s}
$$
with $\mcQ^2_{t+s}$ having fast decreasing kernel at scale $(s+t)^\frac{1}{2}$. That concludes the fact that $\widetilde{\mcQ}^1 \in {\sf GC}^{\frac{b}{8}-2}$.
\end{Dem}

\ssk

The following continuity estimate is then a direct consequence of Proposition \ref{prop:modifiedparaproduit}, since the latter implies that we can reproduce the same argument as for the standard paraproduct in Proposition \ref{prop:para-bis}.

\ssk

\begin{prop}
\label{PropContinuityModifParapdct}   {\sf 
For every $\alpha\in (-3,3)$ and $a,\epsilon\in (0,1)$ with $\alpha-a-\epsilon\in(-3,3)$ and $f\in \calC^\alpha_{p_a}$, we have
\begin{equation*} 
\Big\| \widetilde\Pi^{(b)}_g(f)\Big\|_{\calC^{\alpha-a-\epsilon}_\varpi} \lesssim \kappa^{-\epsilon} \big\|\varpi^{-1} g\big\|_\infty \|f\|_{\calC^\alpha_{p_a}},
\end{equation*}
for every $g\in L^\infty_\varpi$.   }
\end{prop}

\ssk

Last, note the normalization identity 
$$
\widetilde \Pi_{\bf 1}(f) = f -\mathscr{L}^{-1}\mcP^{(b)}_{\bf 1}\mcP^{(b)}_{\bf 1}(\mathscr{L} f)
$$
for every distribution $f\in\mcS'_o$; it reduces to 
$$
\widetilde \Pi_{\bf 1}(f) = f - \mcP^{(b)}_{\bf 1}\mcP^{(b)}_{\bf 1}(f)
$$
if $f_{\big| \tau =0} = 0$. (Use here the support condition on $\varphi$ in the definition of $\mcP$.) Let us also point out here the strongly regularizing effect of the two operators $\mcP_{\bf 1}^{(b)}\mcP_{\bf 1}^{(b)}$ and $\mathscr{L}^{-1} \mcP_{\bf 1}^{(b)}\mcP_{\bf 1}^{(b)}\LL$, denoted by $A$ below, that satisfy the continuity estimate
$$ 
\|A\|_{\calC^\alpha_\omega \to \calC^\beta_\omega} \lesssim 1,
$$
for any $\alpha,\beta\in (-3,3)$ and any space-time weight $\omega$. 

\ssk

We shall fix from now on the parameters $b$ and $\ell_1$, large enough for the above result to hold true.

\begin{rem} \label{rem:paraproduct} The previous Proposition is very interesting because of the following observation: the time-space paraproducts $\widetilde \Pi$ are defined in terms of parabolic cancellations and so do not differentiate the space and the time. Consequently, it is not clear if  the time-space paraproducts $\widetilde \Pi$  may be bounded on $L^\infty_TC^\alpha$ for some $\alpha<0$ (with or without weights). Such property would be very useful since the paracontrolled calculus (as shown later in the study of (PAM) for instance) needs to estimate the composition of $\mathscr{L}^{-1}$ (the resolution of heat equation) with the paraproduct. 
However, following the definition of the paraproduct we have for $f\in  L^\infty_TC^\alpha$ and $g\in \calC^\beta$
$$ \mathscr{L}^{-1} \Pi^{(b)}_g(f)= \widetilde \Pi^{(b)}_g( \mathscr{L}^{-1} f).$$
So if $f \in L^\infty_TC^\alpha$ for some $\alpha\in(-2,0)$ then Schauder estimates imply that $\mathscr{L}^{-1} f \in \calC^{\alpha+2}$ and we may then use the boundedness on H\"older spaces of the modified paraproduct $\widetilde \Pi^{(b)}$.

In conclusion, these new space-time paraproducts seem to be very natural for the paracontrolled calculus. They allow us to get around a commutation between the initial paraproduct and the resolution operator $\mathscr{L}^{-1}$ (which could be a limitation for a higher order paracontrolled calculus) and fits exactly in what paracontrolled calculus requires to solve singular PDEs, modelled on the heat equation.
\end{rem}

\bigskip

\bigskip

\subsection[\hspace{-0.5cm} Commutators and correctors]{Commutators and correctors}
\label{SubsectionIteratedCommutators}

We state and prove in this section two continuity estimates that will be useful in our study of the $3$-dimensional parabolic Anderson model equation and Burgers system in Section \ref{SectionPAM}.

\ssk

\begin{defn}
Let us introduce the following a priori unbounded trilinear operators on $\SSS$. Set
$$
R(f,g,u):= \Pi^{(b)}_u\Big(\Pi^{(b)}_g(f)\Big)  - \Pi^{(b)}_{ug}(f),
$$
and define the \textbf{corrector}
$$
\textsf{C}(f,g,u) := \Pi^{(b)}\left(\widetilde \Pi^{(b)}_{g}(f),u\right) - g\,\Pi^{(b)}(f,u).
$$
\end{defn}

\ssk

This corrector was introduced by Gubinelli, Imkeller and Perkowski in \cite{GIP} under the name of commutator. We prove in the remainder of this section that these operators have good continuity properties in some weighted parabolic H\"older spaces. 

\medskip

\begin{prop} \label{prop:R}
{\sf
Given some space-time weights $\omega_1,\omega_2,\omega_3$, set $\omega:=\omega_1 \omega_2\omega_3$. Let $\alpha,\beta,\gamma$ be H\"older regularity exponents with $\alpha\in(-3,3)$, $\beta\in (0,1)$ and $\gamma\in(-3,0]$. Then if  $\delta:=\alpha+\beta+\gamma\in(-3,3)$ with $\alpha+\beta<3$, 
we have
\begin{equation}
\label{eq:R}
\big\| R(f,g,u) \big\|_{\calC^{\delta}_\omega} \lesssim \|f\|_{\calC^\alpha_{\omega_1}} \, \|g\|_{\calC^\beta_{\omega_2}} \, \|u\|_{\calC^\gamma_{\omega_3}},
\end{equation}
for every $f\in \calC^\alpha_{\omega_1}$, $g\in\calC^\beta_{\omega_2}$ and $u\in \calC^\gamma_{\omega_3}$; so the modified commutator defines a trilinear continuous map from $\calC^\alpha_{\omega_1}\times\calC^\beta_{\omega_2}\times \calC^\gamma_{\omega_3}$ to $\calC^\delta_\omega$.  }
\end{prop}

\medskip

\begin{Dem}
Recall that $ \Pi^{(b)}_g$ is given by a finite sum of operators of the form
$$
\mcA^1_g(\cdot) := \int_0^1  \mcQ^{1\bullet}_t\Big( \mcQ^2_t(\cdot) \, \mcP^1_t(g) \Big) \,\frac{dt}{t},
$$
where $\mcQ^{1},\mcQ^2$ belong at least to ${\sf StGC}^3$. We describe similarly  $\Pi^{(b)}_u$ as a finite sum of operators of the form
$$ 
\mcA^2_u(\cdot) := \int_0^1  \mcQ^{3\bullet}_t\Big( \mcQ^4_t(\cdot)  \mcP^2_t(u) \Big) \,\frac{dt}{t}.
$$
Thus, we need to study a generic modified commutator
$$ 
\mcA^2_u\left(\mcA^1_{g}(f)\right) - \mcA^2_{ug}(f),
$$
and introduce for that purpose the intermediate quantity
$$ 
\mcE(f,g,u) := \int_0^1 \mcQ^{3\bullet}_s\Big( \mcQ^4_s(f) \cdot \mcP^1_s(g) \cdot \mcP^2_s(u) \Big) \,\frac{ds}{s}.  
$$
Note here that due to the normalization $\Pi_1 \simeq \textrm{Id}$, up to some strongly regularizing operator, there is no loss of generality in assuming that 
\begin{equation} 
\label{eq:norma}
 \int_0^1\mcQ^{1\bullet}_t \mcQ^2_t \, \frac{dt}{t} = \int_0^1\mcQ^{3\bullet}_t \mcQ^4_t \, \frac{dt}{t} = \textrm{Id}. 
\end{equation}

\medskip

\noindent \textbf{Step 1. Study of $\mcA^2_u\left( \mcA^1_{g}(f)\right) - \mcE(f,g,u)$.} We shall use a family $\mcQ$ in ${\sf StGC}^{a}$, for some $a>|\delta|$, to control the H\"older norm of that quantity. By definition, and using the normalization \eqref{eq:norma}, the quantity $\mcQ_r \Big(\mcA^2_u\left(\mcA^1_{g}(f)\right) - \mcE(f,g,u) \Big)$ is, for every $r\in(0,1)$, equal to
{\small \begin{align*}
&\int_0^1\int_0^1  \mcQ_r \mcQ^{3\bullet}_s\Big\{ \mcQ^4_s  \mcQ^{1\bullet}_t\Big(  \mcQ^2_t(f)  \mcP^1_t(g) \Big) \cdot \mcP^2_s(u)  \Big\} \,\frac{ds\,dt}{st} - \int_0^1 \mcQ_r \mcQ^{3\bullet}_s\Big( \mcQ^4_s(f) \cdot \mcP^1_s(g) \cdot \mcP^2_s(u) \Big) \,\frac{ds}{s} \\
&= \int_0^1\int_0^1  \mcQ_r \mcQ^{3\bullet}_s\Big\{ \mcQ^4_s  \mcQ^{1\bullet}_t\Big(  \mcQ^2_t(f)  \big(\mcP^1_t(g)-\mcP^1_s(g)\big) \Big) \cdot  \mcP^2_s(u)  \Big\} \,\frac{dsdt}{st},
\end{align*}   }
where in the last line the variable of $\mcP^1_s(g)$ is the one of $\mcQ^{3\bullet}_s$, and so it is frozen through the action of $ \mcQ^4_s\mcQ^{1\bullet}_t$. Then using that $g\in \calC^\beta$ with $\beta \in (0,1)$, we  know by Proposition \ref{prop:hol} that we have, for $\tau\geq \sigma$,
$$ 
\omega_2(x,\tau)^{-1}\Big|\big(\mcP^1_s g\big)(x,\tau) - \big(\mcP^1_t g\big)(y,\sigma)\Big| \lesssim \left(s+t+ \rho\big((x,\tau),(y,\sigma)\big)^2\right)^\frac{\beta}{2} e^{cd(x,y)} \|g\|_{\calC^\beta_{\omega_2}}.
$$
Note that it follows from equation \eqref{EqIteratedG} that the kernel of $\mcQ^4_s  \mcQ^{1\bullet}_t$ is pointwise bounded by $\G_{t+s}$, and allowing different constants in the definition of the Gaussian kernel $\G$, we have
\begin{equation}
\label{eq:gg}
\G_{t+s}\big((x,\tau),(y,\sigma)\big) \left(s+t+ d(x,y)^2\right)^\frac{\beta}{2}e^{cd(x,y)} \lesssim (s+t)^\frac{\beta}{2} \G_{t+s}\big((x,\tau),(y,\sigma)\big).
\end{equation}
So using Lemma \ref{lem:weight} and the cancellation property of the operators $\mcQ$ at an order no less than $a$ (resp. $3$) for $\mcQ$ (resp. the other collections $\mcQ^i$), we deduce that
\begin{align*}
& \left\| \omega^{-1} \mcQ_r \Big(\mcA^2_u\left( \mcA^1_{g}(f)\right) - \mcE(f,g,u)\Big) \right\|_\infty  \\
& \qquad \lesssim \|f\|_{\calC^\alpha_{\omega_1}} \|g\|_{\calC^\beta_{\omega_2}} \|u\|_{\calC^\gamma_{\omega_3}} \int_0^1\int_0^1 \left(\frac{sr}{(s+r)^2}\right)^\frac{a}{2} \left(\frac{st}{(s+t)^2}\right)^\frac{3}{2} t^\frac{\alpha}{2} (s+t)^\frac{\beta}{2} s^\frac{\gamma}{2}   \,\frac{ds\,dt}{st},
\end{align*}
where we used that $\gamma$ is negative to control $\mcP^2_s(u) $. The integral over $t\in(0,1)$ can be computed since $\alpha>-3$ and $\alpha+\beta<3$, and we have
\begin{align*}
& \left\|\omega^{-1} \mcQ_r \Big(\mcA^2_u\left( \mcA^1_{g}(f)\right) - \mcE(f,g,u)\Big) \right\|_\infty  \\
& \qquad \lesssim \|f\|_{\calC^\alpha_{\omega_1}} \|g\|_{\calC^\beta_{\omega_2}} \|u\|_{\calC^\gamma_{\omega_3}} \int_0^1\int_0^1 \left(\frac{sr}{(s+r)^2}\right)^\frac{a}{2}  s^\frac{\delta}{2} \,\frac{ds}{s} \\
& \qquad \lesssim \|f\|_{\calC^\alpha_{\omega_1}} \|g\|_{\calC^\beta_{\omega_2}} \|u\|_{\calC^\gamma_{\omega_3}} r^\frac{\delta}{2},
\end{align*}
uniformly in $r\in(0,1)$ because $|a|>\delta$. That concludes the estimate for the high frequency part. We repeat the same reasoning for the low-frequency part by replacing $\mcQ_r$ with $\mcQ_1$ and conclude that
$$ 
\Big\| \mcA^2_u\left( \mcA^1_{g}(f)\right) - \mcE(f,g,u) \Big\|_{\calC^{\delta}_\omega} \lesssim \|f\|_{\calC^\alpha_{\omega_1}} \|g\|_{\calC^\beta_{\omega_2}} \|u\|_{\calC^\gamma_{\omega_3}}.
$$

\medskip

\textbf{Step 2. Study of $\mcA^2_{ug}- \mcE(f,g,u)$.} This term is simpler than that of Step 1 and can be treated similarly. Note that $\mcQ_r \Big(\mcA^1_g\left(\mcA^2_{u}(f)\right) - \mcE(f,g,u) \Big)$ is equal, for every $r\in (0,1)$, to 
{\small \begin{align*}
&\int_0^1 \mcQ_r \mcQ^{3\bullet}_s\Big( \mcQ^4_s(f)  \mcP^2_s(ug)  \Big)  \,\frac{ds}{s} - \int_0^1 \mcQ_r\mcQ^{3\bullet}_s\Big(\mcQ^4_s(f) \cdot \mcP^1_s(g) \cdot \mcP^2_s(u) \Big) \,\frac{ds}{s} \\
&= \int_0^1 \mcQ_r\mcQ^{3\bullet}_s\Big( \mcQ^4_s(f) \big(\mcP^2_s(ug)-\mcP^1_s(g)\cdot \mcP^2_s(u) \big)  \Big)  \,\frac{ds}{s}.
\end{align*}   }
Now note that since $g\in \calC^\beta$ with $\beta \in (0,1)$, we  know by Proposition \ref{prop:hol}, for $\tau\geq \sigma$,
\begin{align*} 
& \omega_2(x,\tau)^{-1}\Big|g(x,\tau) - \big(\mcP^1_s g\big)(y,\sigma)\Big|  \\
& \lesssim \omega_2(x,\tau)^{-1}\Big|g(x,\tau) - g(y,\sigma)\Big| +  \omega_2(x,\tau)^{-1}\Big|g(y,\sigma) - \big(\mcP^1_t g\big)(y,\sigma)\Big| \\
& \lesssim \left(s+t+ \rho\big((x,\tau),(y,\sigma)\big)^2\right)^\frac{\beta}{2} e^{cd(x,y)} \|g\|_{\calC^\beta_{\omega_2}}.
\end{align*}
Then the same proof as in Step 1 can be repeated.
\end{Dem}

\medskip

As far as the continuity properties of the corrector 
$$
\textsf{C}(f,g,u) = \Pi^{(b)}\Big(\widetilde \Pi^{(b)}_g(f),u\Big) - g\Pi^{(b)}(f,u)
$$ 
are concerned, the next result was proved in an unweighted setting in \cite[Proposition 3.6]{BB} for a space version of the paraproduct $\Pi$; elementary changes in the proof give the following space-time weighted counterpart.

\ssk
 
\begin{prop}  \label{prop:C}   {\sf
Given space-time weights $\omega_1,\omega_2,\omega_3$, set $\omega:=\omega_1\,\omega_2\,\omega_3$. Let $\alpha,\beta,\gamma$ be H\"older regularity exponents with $\alpha\in(-3,3), \beta\in(0,1)$ and $\gamma\in (-\infty,3]$. Set $\delta := (\alpha+\beta)\wedge 3 +\gamma$. If
$$
0 < \alpha + \beta + \gamma < 1 \qquad \textrm{ and }\qquad \alpha + \gamma < 0
$$
then the corrector $\textsf{C}$ is a continuous trilinear map from $\calC^\alpha_{\omega_1}\times\calC^\beta_{\omega_2}\times\calC^\gamma_{\omega_3}$ to $\calC^\delta_\omega$.  }
\end{prop}

\bigskip

\section[\hspace{0.7cm} Anderson and Burgers equations in a 3d background]{Anderson and Burgers equations in a $3$-dimensional background}
\label{SectionPAM}

We are now ready to start our study of the parabolic Anderson model equation 
$$
(\partial_t+L) u = u\zeta
$$
and the multiplicative Burgers system
$$
(\partial_t+L) u + (u\cdot V)u = \textrm{M}_\zeta u
$$
in a $3$-dimensional manifold, using the above tools. Here for Burgers system, we consider a collection of three operators $V:=(V_1,V_2,V_3)$, so$\ell_0=3$ here. We shall study the (PAM) equation in a possibly unbounded manifold, using weighted H\"older spaces, while we shall be working in a bounded setting for the Burgers equation, as its quadratic term does not preserve any 'obvious' weighted space.

\medskip

\subsection[\hspace{-0.5cm} Getting solutions for the (PAM) equation]{Getting solutions for the (PAM) equation}
\label{Subsectionpara}

Let us take the  freedom to assume for the moment that the noise $\zeta$ in the above equations is not necessarily as irregular as white noise. We shall fix from now on a finite positive time horizon $T$. Recall the elementary result on paracontrolled distributions $u$ with derivative $u$ stated in section \ref{SectionParacontrolledNutshell}; such distributions are of the form $u=e^{-Z}v_1$, for some more regular factor $v_1$. This is indeed what happens formally for any solution to the (PAM) equation, since $u\zeta = \Pi_u(\zeta)$, up to some smoother term, and $\mathscr{L}^{-1}\Big(\Pi_u(\zeta)\Big) = \Pi_u(\mathscr{L}^{-1}\zeta)$, up to some more regular remainder. Elaborating formally on this remark leads to the introduction of the following distributions, and the choice of representation for a solution of the (PAM) equation adopted below in Proposition \ref{prop:uv}.

\ssk

For a continuous function $\zeta$ in $C^0_{p_a}$, and $1\leq i\leq 3$, define recursively the following reference distributions/functions 
$$
Z_i := \mathscr{L}^{-1}(Y_i),
$$ 
with 
\begin{align} \label{EqDefnYZ}
 Y_1 := \zeta,  \quad Y_2:= \sum_{i=1}^{\ell_0} V_i(Z_1)^2,  \quad Y_3 := 2\sum_{i=1}^{\ell_0} V_i(Z_1) V_i(Z_2) , 
\end{align}
and define
\begin{align*}
(\star) := -2 \sum_{i=1}^{\ell_0} V_i(Z_1) V_i(Z_3), \quad W_2 := - \sum_{i=1}^{\ell_0} V_i(Z_2)^2
\end{align*}
as well as for $j\in \{1,..,\ell_0\}$ 
$$ 
W_2^j:=\sum_{i=1}^{\ell_0} \Pi^{(b)}\Big(V_i(Z_1)\,, \,V_i\mathscr{L}^{-1}(V_j Z_1)\Big).
$$
Indeed in the term $(\star)$, only the resonant parts in the products have to be defined, since the parapoducts always make sense, so we focus on the resonant part of $(\star)$
$$ 
W_1 := -2\sum_{i=1}^{\ell_0} \Pi^{(b)}\big( V_i(Z_3),V_i(Z_1) \big).
$$
Defining 
\begin{itemize}
   \item the $Y_i$'s as elements of $L^\infty_TC^{\alpha-(5-i)/2} \subset L^\infty_T C^{i\alpha-2}_{p_a}$,   \vspace{0.1cm}
   
   \item the distributions $W_k$ as element of $L^\infty_T C^{k\alpha-1}_{p_a}$,   \vspace{0.1cm}
   
   \item the quantities $W^j_2$ as elements of $L^\infty_T C^{2\alpha-1}_{p_a}$,   \vspace{0.1cm}
\end{itemize}
for some $1/3<\alpha<1/2$ and $a>0$, when $\zeta$ is a space white noise, is the object of the renormalisation step, which shall be done elsewhere. These conditions ensure, by Schauder estimates, Proposition \ref{PropFirstSchauder}, that $Z_i$ is in the parabolic H\"older space $\calC^{i\alpha}_{p_a}$. Note that assuming $W_1$ is an element of $L^\infty_T C^{2\alpha-1}_{p_a}$ ensures that $(\star)$ is an element of $L^\infty_T C^{\alpha-1}_{p_a}$. There is a clear correspondence between the terms defined here and those appearing in Hairer and subsequent analyses of the KPZ equation; see e.g. \cite{FH, GP, HLR3}. In a simplified setting where the vector fields $V_i$ are constant and correspond to the derivation operator in the direction of the $i^\textrm{th}$ vector of the canonical basis, the above terms correspond to 
\begin{equation*}
\begin{split}
&W_2^j = (\partial Z_1)\cdot\partial\mathscr{L}^{-1}(\partial_j Z_1), \quad Z_2 = \mathscr{L}^{-1}\big((\partial Z_1)^2\big), \quad   W_2 = (\partial Z_2)^2   \\
&Z_3 = \mathscr{L}^{-1}\big((\partial Z_1)(\partial Z_2)\big),\quad W_1 = (\partial Z_1)(\partial Z_3).
\end{split}
\end{equation*}
Set 
$$
Z := Z_1+Z_2+Z_3 =: Z_1 + \widetilde Z.
$$
 
\ssk

\begin{prop} \label{prop:uv} {\sf 
The function $u$ is a formal solution of the (PAM) equation if and only if the function 
$$
v := e^{-Z}u
$$ 
is a solution of the equation
\begin{equation} 
\label{eq:v}
\mathscr{L} v = -U v + 2\sum_{i=1}^{\ell_0} V_i(Z)V_i(v), 
\end{equation}
with the same initial condition as $u$ at time $0$. The letter $U$ stands here for $W_1+W_2+W_3$ for an explicit distribution $W_3$ in $L^\infty_TC^{2\alpha-1}$.   }
\end{prop}

We explicitly single out $W_2$ here, and not $W_3$, even though they both belong to the same space, because $W_2$ will later have to be renormalized while $W_3$ will be well-defined as soon as the other quantities $Y_i, W_1, ...$ will be well-defined.

\ssk

\begin{Dem}
Observe that
$$ 
\partial_\tau u = e^{Z} \Big(\partial_\tau v + v\partial_\tau Z_1 + v \partial_\tau \widetilde{Z}\Big),
$$
and using the Leibniz rule on $V_i$'s 
\begin{align*} 
Lu & = e^{Z} \left(-\sum_{i=1}^{\ell_0} V_i(Z)^2v -  V_i^2(Z) v - 2V_i(Z)V_i(v)-V_i^2v \right)   \\
     & = e^{Z} \left(vLZ + Lv - \sum_{i=1}^{\ell_0} V_i(Z)^2v - 2V_i(Z)V_i(v)\right)   \\
     & = e^{Z} \left( vLZ_1 + v L\widetilde Z + Lv - v \sum_{i=1}^{\ell_0} V_i(Z)^2 - 2V_i(Z)V_i(v)\right).
\end{align*}
Due to the definition of $Y_i$'s, we have some telescoping property: 
\begin{align*}
& \mathscr{L}\widetilde Z - \sum_{i=1}^{\ell_0} V_i(Z)^2   \\
& \quad = \mathscr{L} \mathscr{L}^{-1}(Y_2+Y_3) -  \sum_{i=1}^{\ell_0} \Big(V_i(Z_1) + V_i(Z_2) +V_i(Z_3)\Big)^2   \\
& \quad = Y_2+Y_3 - \sum_{i=1}^{\ell_0} \sum_{j,k=1}^3 V_i(Z_j)V_i(Z_k)   \\  
& \quad = W_1+W_2 - \sum_{i=1}^{\ell_0} \sum_{\genfrac{}{}{0pt}{}{j,k=2}{j+k\geq 5}}^3 V_i(Z_j)V_i(Z_k).
\end{align*}
Since we assume that $Z_j\in \calC^{j\alpha}_{p_a}$, it follows that $V_i(Z_j)\in L^\infty_TC^{j\alpha-1}_{p_a}$ and $V_i(Z_k)\in L^\infty_TC^{k\alpha-1}_{p_a}$. Given that $j+k\geq 5$ and $\alpha\in(1/3,1/2)$, at least one of the two numbers $(j\alpha-1)$ and $(k\alpha-1)$ is positive and the other not smaller than $2\alpha-1$. So 
\begin{align} 
U& := \mathscr{L}\widetilde Z - \sum_{i=1}^{\ell_0} V_i(Z)^2 \in W_1+W_2+ L^\infty_TC^{2\alpha-1},
\label{eq:U}
\end{align}
and the result follows. 
\end{Dem}

\ssk 

We solve \eqref{eq:v} using paracontrolled calculus instead of solving directly (PAM).

\ssk 

\begin{defn}
Given $\frac{1}{3}<\beta<\alpha<\frac{1}{2}$ and a time-independent distribution $\zeta\in C^{\alpha-2}_{p_a}$, a \textbf{(PAM)-enhancement of $\zeta$} is a tuple $\widehat{\zeta} := \Big(\zeta, Y_2, Y_3,W_1,W_2, (W_2^j)_j\Big)$, with 
$$
Y_k\in L^\infty_T C^{\alpha-(5-k)/2}_{p_a}
$$ 
and 
$$
W_1,W_2, W_k^j\in L^\infty_T C^{2\alpha-1}_{p_a}.
$$
\end{defn}

\ssk

So the \textbf{space of (PAM)-enhanced distributions} $\widehat{\zeta}$ for the (PAM) equation is here simply the product space 
$$
C^{\alpha-2}_{p_a}\times \prod_{k=2}^3 L^\infty_T C^{\alpha-(5-k)/2}_{p_a} \times \left(L^\infty_T C^{2\alpha-1}_{p_a}\right)^{\otimes (\ell_0+2)}.
$$ 

\medskip

\subsubsection[\hspace{-0.5cm} The paracontrolled approach]{The paracontrolled approach}
\label{Subsubsectionpara}

The study of singular PDEs, such as the Anderson and Burgers equations or \eqref{eq:v}, is a \textit{four step process} from a paracontrolled point of view. Let us sketch it for equation \eqref{eq:v} as an example.   \vspace{0.2cm}

\begin{enumerate}
   \item {\bf Set yourself an ansatz for the solution space, in the form of a Banach space of paracontrolled distributions/functions.}   \vspace{0.2cm }
\end{enumerate}   
        
Given $\frac{1}{3}<\beta < \alpha<\frac{1}{2}$, we choose here to work with functions $v$ paracontrolled by the collection $\Big\{\mathscr{L}^{-1}\big(V_i(Z_1)\big)\Big\}_{i=1}^{\ell_0}$, that is with $v$ of the form
\begin{equation} 
\label{ansatz}   
 v= \sum_{i=1}^{\ell_0} \widetilde{\Pi}^{(b)}_{v_i} \Big(\mathscr{L}^{-1}(V_i Z_1)\Big) + v^\sharp
\end{equation}
for a remainder $v^\sharp\in\calC^{1+\alpha+\beta}_{\varpi p_{-a}}$ and $v_i \in \calC^{\beta}_\varpi$. We refer the reader to Subsection \ref{SubsectionSchauder} for the introduction of weights $p_a$ and $\varpi$. Note that we use the $\widetilde\Pi$ paraproduct and not the $\Pi$ paraproduct. We turn the solution space 
$$
\mcS_{\alpha,\beta}\big(\widehat\zeta\big) := \Big\{(v; v_1,\dots,v_{\ell_0}; v^\sharp)\,\textrm{satisfying the above relations}\Big\}
$$ 
into a Banach space by defining its norm as
\begin{equation}
\label{EqDefnBanachNormSolSpace}
\big\|(v;v_1,\dots,v_{\ell_0};v^\sharp)\big\|_{\alpha,\beta} := \big\|v^\sharp\big\|_{\calC^{1+\alpha+\beta}_{\varpi p_{-a}}} + \sum_{i=1}^{\ell_0} \big\|v_i\big\|_{\calC^{\beta}_{\varpi}}.  
\end{equation} 

{\bf  \begin{enumerate}
           \item[(b)] Recast the equation as a fixed point problem for a map $\Phi$ from the solution space to itself.   \vspace{0.2cm }
        \end{enumerate}  } 
        
This is where we use the continuity properties of the corrector and different paraproducts. In the specific situations of equation \eqref{eq:v}, given $(v;v_1,\dots,v_{\ell_0};v^\sharp)$ in the solution space $\mcS_{\alpha,\beta}\big(\widehat\zeta\big)$, one sets 
$$
y= \mathscr{L} \Big( -U v + 2\sum_{i=1}^{\ell_0} V_i(Z)V_i(v)\Big)
$$ 
and shows that it has a decomposition $(y;y_1,\dots,y_{\ell_0};y^\sharp)$ of the form \eqref{ansatz}. This is where we need all the extra information contained in $\widehat\zeta$. Then, given an initial data $v_0\in \calC^{1+\alpha+\beta}_{\varpi p_{-a}}$, the application $\gamma : (\tau,x) \mapsto e^{-\tau L}(v_0)(x)$, belongs to $\calC^{1+\alpha+\beta}_{\varpi p_{-a}}$ and satisfies
$$
\mathscr{L}\gamma = 0,\qquad \gamma_{\tau=0} = v_0.
$$
We define a continuous map $\Phi$ from the solution space $\mcS_{\alpha,\beta}\big(\widehat\zeta\big)$ to itself setting
$$ 
\Phi := (v;v_1,\dots,v_{\ell_0};v^\sharp) \mapsto (y+ \gamma;y_1,\dots,y_{\ell_0};y^\sharp+\gamma).  \vspace{0.1cm }
$$ 

{\bf  \begin{enumerate}
           \item[(c)] Prove that $\Phi$ is a contraction of the solution space.   \vspace{0.2cm }
\end{enumerate}}           

Recall a parameter $\kappa>1$ appears in the definition of the special weight $\varpi$. We shall see below that the function $y^\sharp$ satisfies the estimate
$$ 
\big\|y^\sharp\big\|_{\calC^{1+\alpha+\beta}_{\varpi p_{-a}}} \leq \kappa^{-\epsilon} \big\|(v;v_1,\dots,v_{\ell_0};v^\sharp)\big\|_{\alpha,\beta},
$$
for some $\epsilon>0$, and that $(y_1,\dots,y_{\ell_0})$ depends only on $v$ and not on $v_1,\dots,v_{\ell_0}$ and $v^\sharp$. These facts provide a quick proof that $\Phi\circ\Phi$ is a contraction of the solution space $\mcS_{\alpha,\beta}\big(\widehat\zeta\big)$. Indeed, given $(v;v_1,\dots,v_{\ell_0};v^\sharp)$ in $\mcS_{\alpha,\beta}\big(\widehat\zeta\big)$, set 
$$ 
\big(z+\gamma;z_1,\dots,z_{\ell_0};z^\sharp+\gamma\big) := \Phi^{\circ 2} (v;v_1,\dots,v_{\ell_0};v^\sharp)\in \mcS_{\alpha,\beta}\big(\widehat\zeta\big).
$$
We know that
\begin{align*}
\big\|z^\sharp\big\|_{\calC^{1+\alpha+\beta}_{\varpi p_{-a}}} & \leq \kappa^{-\epsilon} \big\|(y+\gamma;y_1,\dots,y_{\ell_0};y^\sharp+\gamma)\big\|_{\alpha,\beta} \\
 & \lesssim \kappa^{-\epsilon} \big\|(v;v_1,\dots,v_{\ell_0};v^\sharp)\big\|_{\alpha,\beta}.
\end{align*}
The paracontrolled structure \eqref{ansatz} of $y$ and Schauder estimates also give
\begin{align*}
\|y\|_{\calC^{1+\beta}_{\varpi p_{-a}}} & \lesssim \big\|y^\sharp\big\|_{\calC^{1+\beta}_{\varpi p_{-a}}} + \sum_{j=1}^{\ell_0} \kappa^{-\epsilon} \|y_i\|_{\calC^{1+\alpha}_{\varpi}}   \\
& \lesssim \kappa^{-\epsilon} \big\| \Phi(v,v_1,..,v_{\ell_0},v^\sharp) \big\|_{\alpha,\beta}   \\
& \lesssim \kappa^{-\epsilon} \big\| (v,v_1,..,v_{\ell_0},v^\sharp) \big\|_{\alpha,\beta}.
\end{align*} 
So we conclude that $y$ is controlled with a small bound. Since $(z_1,..,z_{\ell_0})$ will be given by $y$, we then obtain that $(z;z_1,\dots,z_{\ell_0};z^\sharp)$ will be controlled in $\mcS_{\alpha,\beta}\big(\widehat\zeta\,\big)$ with small norms, relatively to the initial $(v;v_1,\dots,v_{\ell_0};v^\sharp)$, so $\Phi \circ \Phi$ will indeed be a small perturbation of the map $(v;v_1,\dots, v_{\ell_0};v^\sharp) \mapsto (\gamma;0,\dots,0;\gamma)$. Then it is standard that if $\kappa^{-\epsilon}$ is small enough, that is $\kappa$ is large enough, then we can apply the fixed point theorem to $\Phi^{\circ 2}$ and conclude that it has a unique fixed point in the solution space  $\mcS_{\alpha,\beta}\big(\widehat\zeta\,\big)$; the same concolusion for $\Phi$ follows as a consequence.   \vspace{0.2cm}

{\bf  \begin{enumerate} 
           \item[(d)] Renormalisation step.   \vspace{0.1cm}
        \end{enumerate}  } 

The defining relations for $Z_i$ found in step \textbf{(b)} actually involve some terms that cannot be defined by purely analytical means when $\zeta$ is a white noise, but which make perfect sense for a regularized version $\zeta^\epsilon$ of $\zeta$. Their proper definition requires a renormalisation procedure that consists in defining them as \textit{limits in probability}, in some parabolic H\"older spaces, of suitably modified versions of their regularized versions (with $\zeta^\epsilon$ in place of $\zeta$), which essentially amounts \textit{in the present setting} to adding to them some deterministic functions or constants. (This may be trickier in other situations as the theory of regularity structures makes it clear.) Given the inductive construction of the $Z_i$, this renormalisation step also needs to be done inductively. At $\epsilon$ fixed, this addition of deterministic quantities in the defining relations for $Z_i$ defines another map $\Phi^\epsilon$ from the solution space to itself that can eventually be equivalent to consider a renormalised equation with noise $\zeta^\epsilon$, with $\epsilon$-dependent terms added in the equation, when compared to the initial equation. Write $u^\epsilon$ for its solution. In the end, we get, from the continuity of fixed points of parameter-dependent uniformly contracting maps, a statement of the form: \textit{Let $\Phi$ stand for the map constructed by taking as reference distributions/functions $Z_i$ the limits, in probability, of their renormalised versions. Then the functions $u^\epsilon$ converge in probability to the solution $u$ of the fixed point problem of the map $\Phi$.} 
 
\bigskip

We shall do here the first three steps of the analysis for both the Anderson and Burgers equations, leaving the probabilistic work needed to complete the renormalisation step to another work; we shall nonetheless give in Section \ref{SectionRenormalisation} some hints as to what is going on.

\bigskip

\subsection[\hspace{-0.5cm} The deterministic PAM equation]{The deterministic PAM equation}
\label{SubsectionPAM}
 
Given what was said in the preceding section, the main work for solving the (PAM) equation consists in proving the following result.

\ssk

\begin{thm} \label{thm:pointfixe} {\sf 
Let $\frac{1}{3}<\alpha<\frac{1}{2}$ be given. Choose $\beta<\alpha$, the positive parameter $a$ in the weight $p_a$, and  $\epsilon>0$, such that 
$$ 
2\alpha+ \beta>1 \qquad \textrm{and} \qquad 8(a+\epsilon)\leq \alpha-\beta.
$$
Given an enhanced distribution $\widehat{\zeta}$, one can extend the product operation 
$$
v \in C^\infty_c(M)\mapsto -U v + 2\sum_{i=1}^{\ell_0} V_i(Z)V_i(v)
$$
to the space $\mcS_{\alpha,\beta}\big(\widehat{\zeta}\,\big)$ into an operation $\widehat{v}\mapsto -\widehat{U} \widehat{v} + 2\sum_{i=1}^{\ell_0} V_i(\widehat{Z})V_i(\widehat{v})$, so that setting 
$$ 
y := \mathscr{L}^{-1}\Big[-\widehat{U} \widehat{v} + 2\sum_{i=1}^{\ell_0} V_i(\widehat{Z})V_i(\widehat{v})\Big],
$$ 
and $y_i:=2vV_i(Z_3)+2V_i(v)$, there exists $y^\sharp\in \calC^{1+\alpha+\beta}_{\varpi p_{-a}}$ such $\big(y;y_1,\dots,y_{\ell_0};y^\sharp\big)$ is an element of the solution space $\mcS_{\alpha,\beta}\big(\widehat{\zeta}\,\big)$, and
\begin{equation}  
\label{eq:fin}
\begin{split}
\Big\| \big(y;y_1,\dots,y_{\ell_0};y^\sharp\big) \Big\|_{\alpha,\beta} \lesssim  \, \Big\|\big(v;v_1,\dots,v_{\ell_0};v^\sharp\big) \Big\|_{\alpha,\beta} \\
\big\| y^\sharp \big\|_{\calC^{1+\alpha+\beta}_{\varpi p_{-a}}} \lesssim  \kappa^{-\epsilon}\, \Big\|\big(v;v_1,\dots,v_{\ell_0};v^\sharp\big) \Big\|_{\alpha,\beta}.
\end{split}
\end{equation}   } 
\end{thm}

\ssk

\begin{Dem}
First, we note that since $v$ satisfies the ansatz \eqref{ansatz} and $2a <\alpha-\beta$, we know from Schauder estimates that 
$$ 
v\in \calC^{1+\alpha}_{\varpi p_a} \cap \calC^{1+\beta}_{\varpi p_{-a}}.
$$

\ssk

\textbf{Step 1.} We first consider the part $Uv$ where we recall that $U= W_1+W_2 + W_3$ for some $W_3 \in L^\infty_T C^{2\alpha-1}_{p_a}$. 
Using the paraproduct algorithm, one gets
$$ W_3 v = \Pi^{(b)}_{W_3}(v) + \Pi^{(b)}_v(W_3) + \Pi^{(b)}(v,W_3).$$
By the boundedness of paraproducts, Proposition \ref{PropRegularityParaproduct}, and Schauder estimates, Proposition \ref{prop:schauder-bis}, we get
$$ \Pi^{(b)}_{W_3}(v)\in \calC^{2\alpha+\beta}_\varpi \qquad \textrm{so} \qquad  \mathscr{L}^{-1} \Pi^{(b)}_{W_3}(v) \in \calC^{2+3\beta}_{\varpi p_{-a}} \subset \calC^{1+\alpha+\beta}_{\varpi p_{-a}}$$
with
$$
\left\| \mathscr{L}^{-1} \Pi^{(b)}_{W_3}(v) \right\|_{\calC^{1+\alpha+\beta}_{\varpi p_{-a}}} \lesssim \kappa^{-\epsilon} \|v\|_{\calC^{1+\beta}_{\varpi p_{-a}}}$$
since $2\epsilon+2a<\alpha-\beta$ and $\alpha<1$. For the resonant part, a similar reasoning with Proposition \ref{PropRegularityParaproduct} yields
$$
\Pi^{(b)}(v,W_3)\in \calC^{2\alpha+\beta}_\varpi \qquad \textrm{so} \qquad \mathscr{L}^{-1} \Pi^{(b)}(v,W)\in \calC^{2+3\beta}_{\varpi p_{-a}}
$$
with
$$ 
\left\| \mathscr{L}^{-1} \Pi^{(b)}(v,W) \right\|_{\calC^{1+\alpha+\beta}_{\varpi p_{-a}}} \lesssim \kappa^{-\epsilon} \|v\|_{\calC^{1+\beta}_{\varpi p_{-a}}}.
$$
For the second paraproduct, we use the modified paraproduct and its boundedness, Proposition \ref{PropContinuityModifParapdct}, to have $\mathscr{L}^{-1} \Pi^{(b)}_v(W_3) = \widetilde \Pi^{(b)}_v( \mathscr{L}^{-1} W_3)$, hence since $\mathscr{L}^{-1} W_3 \in \calC^{1+2\alpha}_{p_a}$ we have $ \mathscr{L}^{-1} \Pi^{(b)}_v(W_3) \in \calC^{1+\alpha+\beta}_{\varpi p_{-a}}$ with
$$ 
\left\| \mathscr{L}^{-1} \Pi^{(b)}_v(W_3) \right\|_{\calC^{1+\alpha+\beta}_{\varpi p_{-a}}} \lesssim \kappa^{-\epsilon } \|v\|_{\calC^{1+\alpha}_{\varpi p_a}},
$$
since $4(a+\epsilon)\leq \alpha-\beta$. So we have $\mathscr{L}^{-1} \big(W_3v\big)\in \calC^{1+\alpha+\beta}_{\varpi p_{-a}}$, with an acceptable bound.

\ssk

The term $W_2$ is an element of $L^\infty_T C^{2\alpha-1}_{p_a}$, so using the same reasoning yields that $ \mathscr{L}^{-1}\big(W_2v\big) \in \calC^{1+\alpha+\beta}_{\varpi p_{-a}}$ with an acceptable bound.

\ssk

The term $W_1$ is an element of $L^\infty_T C^{\alpha-1}_{p_a}$, so it is really more singular than the two previous terms. Recall its definition
$$
W_1= - 2 \sum_{i=1}^{\ell_0} V_i(Z_1) V_i(Z_3)
$$
with $V_i(Z_3)$ in $\calC^{\alpha}_{p_a}$, since $Z_3$ is an element of $\calC^{1+\alpha}_{p_a}$. So $W_1$ is in $\calC^{\alpha-1}_{p_a}$, and since $v\in \calC^{1+\alpha}_{\varpi p_{a}}$, we have
$$ 
\Pi^{(b)}_{W_1}(v) \in \calC^{2\alpha}_{\varpi p_{2a}} \quad \textrm{and} \quad \Pi^{(b)}(W_1,v) \in \calC^{2\alpha}_{\varpi p_{2a}}.
$$
Using Schauder estimates one obtains
$$ 
\left\| \mathscr{L}^{-1} \Pi^{(b)}_{W_1}(v) \right\|_{\calC^{1+\alpha+\beta}} + \left\| \mathscr{L}^{-1} \Pi^{(b)}(W_1,v) \right\|_{\calC^{1+\alpha+\beta}} \lesssim \kappa^{-\epsilon } \|v\|_{\calC^{1+\alpha}_{\varpi p_a}}.
$$
It remains us to study the paraproduct term
$$ 
\Pi^{(b)}_{v}(W_1) = {\sf I} + {\sf II} + {\sf III},
$$
with 
\begin{align*} 
{\sf I} & := - 2 \sum_{i=1}^{\ell_0} \Pi^{(b)}_v \Big(\Pi^{(b)}_{V_i(Z_3)} \big(V_i(Z_1) \big)\Big)   \\
{\sf II} & := - 2 \sum_{i=1}^{\ell_0} \Pi^{(b)}_v \Big(\Pi^{(b)} \big( V_i(Z_3),V_i(Z_1) \big)\Big)   \\
{\sf III} & := - 2 \sum_{i=1}^{\ell_0} \Pi^{(b)}_v \Big(\Pi^{(b)}_{V_i(Z_1)} \big(V_i(Z_3) \big)\Big).
\end{align*}

By easy considerations on paraproducts, the third term ${\sf III}$ belongs to $\calC^{2\alpha-1}_{\varpi p_{a}}$ and $\mathscr{L}^{-1}({\sf III}) \in \calC^{1+\alpha+\beta}_{\varpi p_{-a}}$, with acceptable bounds, because $Z_3$ is an element of $\calC^{1+\alpha}_{p_a}$. Moreover, since we assume that $W_1=\sum_{i=1}^{\ell_0} \Pi^{(b)}\big( V_i(Z_3),V_i(Z_1) \big)$ is an element of $L^\infty_T\calC^{2\alpha-1}_{p_a}$, the second term ${\sf II}$ also satisfies $\mathscr{L}^{-1}({\sf II}) \in \calC^{1+\alpha+\beta}_{\varpi p_{-a}}$. Using the regularity of $v\in \calC^{1+\alpha}_{\varpi p_a}\subset L^\infty_{\varpi p_{a}}$ and Proposition \ref{prop:R} for the commutation property, we deduce that
$$ 
{\sf I} \in - 2 \sum_{i=1}^{\ell_0}  \Pi^{(b)}_{vV_i(Z_3)} [V_i(Z_1)] + \calC^{4\alpha-2}_{\varpi p_{3a}}
$$
and consequently
$$ 
\mathscr{L}^{-1}({\sf I}) \in - 2 \sum_{i=1}^{\ell_0}  \widetilde\Pi^{(b)}_{vV_i(Z_3)} [\mathscr{L}^{-1} V_i(Z_1)] + \calC^{1+\alpha+\beta}_{\varpi p_{-a}},
$$
with an acceptable bound for the remainder since $8(a+\epsilon)+1<3\alpha-\beta$.

\ssk

At the end, we have obtained that
$$ 
\mathscr{L}^{-1}(Uv) \in \Big\{2 \sum_{i=1}^{\ell_0}  \widetilde\Pi^{(b)}_{vV_i(Z_3)} [\mathscr{L}^{-1} V_i(Z_1)] + \calC^{1+\alpha+\beta}_{\varpi p_{-a}}\Big\},
$$
which proves that $\mathscr{L}^{-1}(Uv)$ is paracontrolled by the collection $\big(\mathscr{L}^{-1} V_i(Z_1)\big)_i$ and the remainder has a bound controlled by $\kappa^{-\epsilon}$.

\bigskip

\textbf{Step 2.} Let now focus on the term $ \sum_{i=1}^{\ell_0} V_i(Z)V_i(v)$. Fix an index $i$ and write
$$ 
V_i(Z)V_i(v) = \Pi^{(b)}_{V_i(v)}\big(V_i(Z)\big) + \Pi^{(b)}_{V_i(Z)}\big(V_i(v)\big) + \Pi^{(b)} \big(V_i(Z),V_i(v) \big).
$$
The second term is of regularity $2\alpha-1$ and using the modified paraproduct, Schauder estimate and the fact that we have $v\in \calC^{1+\alpha}_{\varpi p_{a}}$, we see that
$$ 
\mathscr{L}^{-1}\left[\Pi^{(b)}_{V_i(Z)}(V_i(v))\right] = \widetilde\Pi^{(b)}_{V_i(Z)}(\mathscr{L}^{-1} V_i(v)) \in \calC^{1+\alpha+\beta}_{\varpi p_{-a}}.
$$
We proceed as follows to study the resonant part. First, since $\alpha>1/3$, we have
\begin{align*}
 \Pi^{(b)}(V_i(Z),V_i(v)) & \in \Big\{\sum_{j=1}^{\ell_0} \Pi^{(b)}\left(V_i(Z_1),V_i \widetilde{\Pi}^{(b)}_{v_j} \big[\mathscr{L}^{-1}(V_jZ_1)\big] \right) + \calC^{3\alpha-1}_{\varpi p_{2a}}\Big\}.
\end{align*}
Consider the modified resonant part
$$ 
\overline{\Pi}_i^{(b)}(f,g) :=\Pi^{(b)}(f,V_i g)
$$
and the corresponding corrector 
$$
\overline{{\sf C}}_i(f,g,h) := \overline{\Pi}_i^{(b)} \big(\widetilde \Pi^{(b)}_g(f),h \big) - g\overline{\Pi}_i^{(b)}(f,h).
$$ 
Then since in the study of the resonant part and the commutator, we can change the {\it localization operators},  so we can integrate an extra $V_i$ operator, we get boundedness of $\overline{\Pi}_i^{(b)}$ from $\calC^\alpha \times \calC^\beta$ to $\calC^{\alpha+\beta-1}$ as soon as $\alpha+\beta-1>0$, and boundedness of the corrector $\overline{{\sf C}}_i$ from $\calC^\alpha \times \calC^\beta \times \calC^\gamma$ into $\calC^{\alpha+\beta+\gamma-1}$ as soon as $\alpha+\beta+\gamma-1>0$, proceeding exactly in the same way as above for $\Pi^{(b)}$ and ${\sf C}$. Using this commutator, we see that $ \sum_{i=1}^{\ell_0} \Pi^{(b)}(V_i(Z),V_i(v))$ is an element of the space 
\begin{align*}
\sum_{i=1}^{\ell_0} \sum_{j=1}^{\ell_0} v_j.\Pi^{(b)}\Big(V_i(Z_1), V_i\mathscr{L}^{-1}(V_jZ_1) \Big) + \sum_{i=1}^{\ell_0} \overline{C_i}\big(V_i(Z_1),v_j,V_i\mathscr{L}^{-1}(V_jZ_1) \big)+ \calC^{2\alpha+\beta-1}_{\varpi p_{2a}},
\end{align*} 
that is an element of 
$$
\sum_{j=1}^{\ell_0} v_j. W_2^j + \calC^{2\alpha+\beta-1}_{\varpi p_{2a}} \subset L^\infty_TC^{2\alpha-1}_{\varpi p_{2a}},
$$
since $W_2^j \in L^\infty_TC^{2\alpha-1}_{p_a}$ and $2\alpha+\beta>1$. In the end, we conclude that
\begin{equation*}
\begin{split}
y &:= \mathscr{L}^{-1}\left[-U v + 2\sum_{i=1}^{\ell_0} V_i(Z)V_i(v) \right]   \\
   &= 2\sum_i \widetilde\Pi^{(b)}_{vV_i(Z_3) + V_i(v)}(\mathscr{L}^{-1} V_i(Z_1)) + \calC^{1+\alpha+\beta}_{\varpi p_{-a}},
\end{split}
\end{equation*}
as expected. Observe that $V_i(Z_3)$ is of parabolic regularity $(3\alpha-1)$, so $vV_i(Z_3)$ and $V_i(v)$ belong to $\calC^{\beta}_{\varpi}$ .
\end{Dem}

\medskip

We can then apply the contraction principle, such as explaned above in Step (c) in section \ref{Subsubsectionpara}.

\medskip

Given $v_0\in C^{1+\alpha+\beta}_{\varpi_0 p_{-a}}$ , write $\mcS_{\alpha,\beta}^{v_0}\big(\widehat{\zeta}\,\big)$ for those tuples $(v;v_1,\dots,v_{\ell_0};v^\sharp)$ in $\mcS_{\alpha,\beta}\big(\widehat{\zeta}\,\big)$ with $v_{|\tau=0}=v_0$. As the function $\gamma:=(x,\tau)\mapsto \big(e^{-\tau L}\big)(v_0)(x)$ belongs to $\calC^{1+\alpha+\beta}_{\varpi p_{-a}}$ and is the solution of the equation
$$ (\partial_\tau +L)(\gamma)=0, \quad \gamma_{\tau=0}=v_0,$$
we define a map $\Phi$ from $\mcS_{\alpha,\beta}^{v_0}\big(\widehat{\zeta}\,\big)$ to itself setting
$$
\Phi(v;v_1,\dots,v_{\ell_0};v^\sharp) = \big(y+\gamma;2vV_1(Z_3)+2V_1(v),\dots,2vV_{\ell_0}(Z_3)+2V_{\ell_0}(v);y^\sharp+\gamma\big),
$$
with
$$ 
y := \mathscr{L}^{-1}\Big(-\widehat{U} \widehat{v} + 2\sum_{i=1}^{\ell_0} V_i(\widehat{Z})V_i(\widehat{v})\Big),
$$ 
and $y^\sharp$ given by the previous theorem. Note that the map $\Phi$ depends continuously on the enhanced distribution $\widehat{\zeta}$; the next global in time well-posedness result is then a direct consequence of Theorem \ref{thm:pointfixe}.

\medskip

\begin{thm} \label{thm:pointfixe-bis} {\sf 
Let us work under assumption {\bf (A)}, and let $\frac{1}{3}<\alpha<\frac{1}{2}$ be given. Choose $\beta<\alpha$, the positive parameter $a$ in the weight $p_a$, and  $\epsilon>0$, such that 
$$ 
2\alpha+ \beta>1 \qquad \textrm{and} \qquad 8(a+\epsilon)\leq \alpha-\beta.
$$ 
Then, one can choose a positive parameter $\kappa$, in the definition of the special weight $\varpi$, large enough to have the following conclusion. Given $v_0\in C^{1+\alpha+\beta}_{\varpi_0 p_{-a}}$, the map $\Phi$ has a unique fixed point $(v,v_1,..,v_{\ell_0},v^\sharp)$ in $\mcS_{\alpha,\beta}^{v_0}\big(\widehat{\zeta}\,\big)$; it depends continuously on the enhanced distribution $\widehat\zeta$, and satisfies the identity $v_i = 2vV_i(Z_3)+2V_i(v)$ for $i=1,..,\ell_0$. This distribution is the solution of the singular PDE
\begin{equation}
\label{EqColeHopfTransformed}
\mathscr{L}v = -Uv+ 2 \sum_{i=1}^{\ell_0} V_i(\widehat{Z})V_i(\widehat{v})
\end{equation}
with $v_{|\tau=0}=v_0$.  The function $u = e^{Z}v$ is then the unique solution of the (PAM) equation with initial data $v_0$.   }
\end{thm}
 
\ssk 
 
If the ambient space $M$ is bounded, then we do not have to take care of the infinity in the space variable, and one can prove a global (in time) result by considering the weight $\varpi(x,\tau) = e^{\kappa \tau}$ with a large enough parameter $\kappa$.

\bigskip

\subsection[\hspace{-0.5cm} The stochastic PAM equation]{The stochastic PAM equation}
\label{SubsectionStochasticPAM}

Recall the time-independent white noise over the measure space $(M,\mu)$ is the centered Gaussian process $\xi$ indexed by $L^2(\mu)$, with covariance
\begin{equation*} 
\E\big[\xi(f)^2\big] = \int f^2(x)\,\mu(dx). 
\end{equation*}   
It can be proved \cite{BB} to have a modification with values in the spatial H\"older space $C^{-\frac{\nu}{2}-\epsilon}_{p_a}$, for all positive constants $\epsilon$ and $a$, where $\nu$ is the Ahlfors dimension of $(M,d,\mu)$ -- its dimension in our Riemannian setting. We take $\nu=3$ here. We still denote this modification by the same letter $\xi$. The study of the stochastic singular PDE of Anderson
$$
\mathscr{L} u = u\xi
$$
can be done in the present setting. This requires a renormalisation step needed to show that the quantities $\Xi=Y_j,W_j,...$ can be defined as elements of suitable functional spaces, as limits in probability of distributions of the form $\Xi^\epsilon-\lambda^\epsilon$, where $\Xi^\epsilon$ is given by formula 
\eqref{EqDefnYZ} with $\zeta = \xi^\epsilon := e^{-\epsilon L}\xi$, the regularized version of the noise via the semigroup, and $\lambda^\epsilon$ are some deterministic functions. This renormalisation step is not done here; Section \ref{SectionRenormalisation} gives however a flavour of what is involved in this process in the present setting. Note that the two dimensional setting was studied in depth in \cite{BB}, with spatial paraproducts used there instead of space-time paraproducts. We then formulate this renormalisation step as an assumption in the present work. Recall the definition of $Z_2, W_1,W_2^j,W_2,Y_3$ and $Z_3$ given in Section \ref{Subsectionpara}. 

\medskip

\noindent \textbf{\textsf{Assumption (B) -- Renormalisation.}}   {\sf  
Let $\xi$ stand for white noise on $M$, and for $\epsilon>0$, denote by $\xi^\epsilon:=e^{-\epsilon L} \xi$ its regularized version, and by $\Xi^\epsilon$ the  distributions corresponding to $\Xi=Y_1,Z_1,Y_2,Z_2$ that one obtains by replacing $\zeta$ by $\xi^\epsilon$.   \vspace{0.15cm}

\begin{enumerate}
   \item[\textbf{(a)}] There exists a family $(\lambda_1^\epsilon)_{0<\epsilon\leq 1}$ of deterministic functions such that $Y_2^\epsilon-\lambda_1^\epsilon$ is $\epsilon$-uniformly bounded and converging in $C_TC^{\alpha-3/2}_{p_a}$, for every $a\in(0,1)$ and any $\alpha<1/2$.   \vspace{0.15cm}
   
   \item[\textbf{(b)}] Use the upper $\epsilon$-exponent in $\Xi=Z_2, W_1,W_2^j,W_2,Y_3,Z_3$  to denote the quantities that one obtains by replacing $Z_1$ by $Z_1^\epsilon$, and $Y_2$ by $Y_2^\epsilon - \lambda_1^\epsilon$. For any $\alpha<1/2$,
\begin{itemize}
   \item the distributions $Z_2^\epsilon, Y_3^\epsilon, W_2^{j,\epsilon}$, are $\epsilon$-uniformly bounded and converging in $C_TC^{\alpha-2}_{p_a}$, respectively $C_TC^{\alpha-1}_{p_a}$ and $C_TC^{2\alpha-1}_{p_a}$, for every $a\in(0,1)$;   \vspace{0.1cm}

   \item there exists deterministic functions $\lambda_{2,1}^{\epsilon}$ and $\lambda_{2,2}^\epsilon$ such that the distributions
$ W_1^\epsilon-\lambda_{2,1}^\epsilon$ and $W_2^\epsilon-\lambda_{2,2}^\epsilon$, are $\epsilon$-uniformly bounded and converging in $C_TC^{2\alpha-1}_{p_a}$ and respectively $C_TC^{2\alpha-1}_{p_a}$, for every $a\in(0,1)$.
\end{itemize}
\end{enumerate}   }

\medskip

This assumption about the renormalisation process for the above quantities practiaclly means that one can renormalise the most singular quantity $Y_2$ by substracting an $\epsilon$-dependent deterministic function, and that once this has been done, no extra renormalisation is needed for the terms $Y_3^\epsilon$ and $W_2^{j,\epsilon}$. At the same time, the quantities $W_1^\epsilon$ and $W_2^\epsilon$ have to be renormalized and this operation can be done by subtracting deterministic functions -- essentially their expectation.
 
\medskip 
 
Write $\overline{Z}^\epsilon$ and $\overline{U}^\epsilon$ for the renormalized versions of $Z^\epsilon$ and $U^\epsilon$.
By tracking in the proof of Theorem \ref{thm:pointfixe} the changes induced by such a renormalisation of $Y_2^\epsilon$, $W_1^\epsilon$ and $W_2^\epsilon$, we see that if $\big(v^\epsilon;v_1^\epsilon,\dots,v_{\ell_0}^\epsilon;v^{\epsilon,\sharp}\big)$ satisfies ansatz \eqref{ansatz} with $\zeta=\xi^\epsilon$, and setting 
$$ 
y^\epsilon := \mathscr{L}^{-1}\Big(-(\overline{U}^\epsilon+\lambda_1^\epsilon-\lambda_{2,1}^\epsilon-\lambda_{2,2}^\epsilon) v^\epsilon + 2\sum_{i=1}^{\ell_0} V_i(\overline{Z}^\epsilon) V_i(v^\epsilon)\Big),
$$
then the tuple
$$ 
\big(y^\epsilon;2v^\epsilon V_1(Z_3^\epsilon)+2V_1(v^\epsilon),\dots,2v^\epsilon V_{\ell_0}(Z_3^\epsilon)+2V_{\ell_0}(v^\epsilon);y^{\epsilon,\sharp}+\gamma\big)
$$
also satisfies the ansatz. The renormalisation quantity $\lambda_1^\epsilon-\lambda_{2,1}^\epsilon-\lambda_{2,2}^\epsilon$ comes from the definitions of $U$ and $U^\epsilon$, since we have after replacement of $Y_2^\epsilon$ by $\overline{Y_2}^\epsilon=Y_2^\epsilon-\lambda_1^\epsilon$
\begin{align*}
 U^\epsilon &= \overline{Y_2}^\epsilon+Y_3^\epsilon - \sum_{i=1}^{\ell_0} \sum_{j,k=1}^{3} V_i(Z_j^\epsilon) V_i(Z_k^\epsilon) \\
 &  = Y_2^\epsilon+Y_3^\epsilon - \sum_{i=1}^{\ell_0} \sum_{j,k=1}^{3} V_i(Z_j^\epsilon) V_i(Z_k^\epsilon) - \lambda_1^\epsilon \\
 & = W_1^\epsilon +W_2^\epsilon - \sum_{i=1}^{\ell_0} \sum_{\genfrac{}{}{0pt}{}{j,k=2}{j+k\geq 5}}^3 V_i(Z_j^\epsilon)V_i(Z_k^\epsilon) \\
 & = \overline{W_1}^\epsilon +\overline{W_2}^\epsilon - \sum_{i=1}^{\ell_0} \sum_{\genfrac{}{}{0pt}{}{j,k=2}{j+k\geq 5}}^3 V_i(Z_j^\epsilon)V_i(Z_k^\epsilon)- \lambda_1^\epsilon + \lambda_{2,1}^\epsilon + \lambda_{2,2}^\epsilon.
\end{align*}

\ssk

\begin{thm}
\label{ThmRenormalisationPAM} {\sf 
Let us work under assumptions {\bf (A)} and {\bf (B)}, and let $\frac{1}{3}<\alpha<\frac{1}{2}$ be given. Choose $\beta<\alpha$, the positive parameter $a$ in the weight $p_a$, and  $\epsilon>0$, such that 
$$ 
2\alpha+ \beta>1 \qquad \textrm{and} \qquad 8(a+\epsilon)\leq \alpha-\beta.
$$ 
One can choose a large enough parameter $\kappa$ in the definition of the special weight $\varpi$ for the following to hold. There exists a sequence of {\it deterministic functions} $\big(\lambda_j^\epsilon\big)_{0<\epsilon\leq 1}$ such that if $v^\epsilon$ stands for the solution of the renormalized equation 
   \begin{equation} 
   \label{eq:renorr}
   \mathscr{L} v^\epsilon = \Big[-(\overline{U}^\epsilon+\lambda_1^\epsilon-\lambda_{2,1}^\epsilon-\lambda_{2,2}^\epsilon) v^\epsilon + 2\sum_{i=1}^{\ell_0} V_i(\overline{Z}^\epsilon) V_i(v^\epsilon)\Big]\qquad v^\epsilon(0)=v_0  
   \end{equation}
with initial condition $v_0 \in \calC^{1+\alpha+\beta}_{\varpi_0 p_{-a}}$, then $v^\epsilon$ converges in probability to the solution $v\in \calC^{1+\alpha}_{\varpi p_a}$ of the Cole-Hopf transformed (PAM) equation \eqref{EqColeHopfTransformed} constructed from the enhancement of the noise given by assumption {\bf (B)} and Theorem \ref{thm:pointfixe}.    } 
\end{thm}  

\medskip

By reproducing the calculations of Subsection \ref{Subsectionpara}, we observe that $v^\epsilon$ is solution of equation \eqref{eq:renorr} if and only if  $u^\epsilon:=e^{\overline{Z}^\epsilon}v^\epsilon$ is solution of the equation
\begin{equation} 
\label{eq:renorr-u}
 \mathscr{L} u^\epsilon = \big(\xi^\epsilon-\lambda_1^\epsilon+\lambda_{2,1}^\epsilon+\lambda_{2,2}^\epsilon\big) u^\epsilon, \qquad u^\epsilon(0)=v_0.  
\end{equation}

\ssk

\begin{thm}
\label{ThmRenormalisationPAM-u} {\sf 
Let us work under assumptions {\bf (A)} and {\bf (B)}, and let $\frac{1}{3}<\alpha<\frac{1}{2}$ be given. Choose $\beta<\alpha$, the positive parameter $a$ in the weight $p_a$, and  $\epsilon>0$, such that 
$$ 
2\alpha+ \beta>1 \qquad \textrm{and} \qquad 8(a+\epsilon)\leq \alpha-\beta.
$$ 
One can choose a large enough parameter $\kappa$ in the definition of the special weight $\varpi$ for the following to hold. There exists a sequence of {\it deterministic functions} $\big(\lambda_j^\epsilon\big)_{0<\epsilon\leq 1}$ such that if $u^\epsilon$ stands for the solution of the renormalized equation \eqref{eq:renorr-u} with initial condition $v_0 \in \calC^{1+\alpha+\beta}_{\varpi_0 p_{-a}}$, then $u^\epsilon$ converges in probability to the solution $u\in \calC^{\alpha}_{\varpi p_{2a}}$ of the (PAM) equation constructed from the enhancement of the noise given by assumption {\bf (B)} and Theorem \ref{thm:pointfixe}.   } 
\end{thm}  

\medskip

This result is coherent with the result of Hairer and Labb\'e proved in \cite{HLR3}. Indeed, in \cite[Equation (5.3)]{HLR3}, the quantities involving an odd order of noises need no renormalisation terms, like us; nor does the term ${\widetilde W}_2^j$, which is of even order but involves an extra derivative $V_j$. We give more insights on the latter term at the end of Section \ref{SectionRenormalisation}, and explain why this extra derivative with anti-symmetry properties implies that the potential renormalisation term is actually null, as in \cite{HLR3}.

\bigskip

\subsection[\hspace{-0.5cm} The multiplicative Burgers equation]{The multiplicative Burgers equation}
\label{SubsectionDeterministicBurger}

We study in this last section the multiplicative Burgers system
$$
(\partial_t+L) u + (u\cdot V)\,u = \textrm{M}_\zeta u
$$
in the same $3$-dimensional setting as before with three operators $V:=(V_1,V_2,V_3)$ forming an elliptic system. Here the solution $u=(u^1,u^2,u^3)$ is a function with ${\mathbb R}^3$-values and $(u\cdot V)\,u$ has also $3$ coordinates with by definition
$$ \left[(u\cdot V)\,u\right]^j := \sum_{i=1}^3 u^i V_i(u^j).$$

To study this equation, we have to make the extra assumption that the ambient space $M$ is bounded. Indeed the boundeness of the ambient space is crucial here, as using weighted H\"older spaces, it would not be clear how to preserve the growth at infinity dictated by the weight when dealing with the quadratic nonlinearity. In such a bounded framework, we do not need to use spatial weights and consider instead the unweighted H\"older spaces $\calC^\gamma$ -- or rather we work for convenience with a weight in time 
\begin{equation}
\label{eq:varpi2}
\varpi(x,\tau) := e^{\kappa \tau}. 
\end{equation} 
We stick to the notations of the previous section. The study of Burgers' system requires a larger space of enhanced distributions than the study of the $3$-dimensional (PAM) equation; the additional components include those quantities that need to be renormalised to make sense of the term $(u\cdot V)\,u$, when $\zeta$ is an element of $\calC^{\alpha-2}$, such as space white noise. 

\medskip
 
We first rewrite Burgers system in a more convenient form, as we did for the (PAM) equation. For each cooordinate exponent $j=1,2,3$, we define $Z^j_\alpha, W^{j}_{\beta}$ from $\zeta^j$ as above. Then consider a function $u : M \mapsto {\mathbb R}^3$ defined by
$$ 
u^j = e^{Z^j} v^j 
$$
with $v : M \mapsto {\mathbb R}^3$. Then observe that $u$ is formally a solution of 3-dimensional Burgers system on $M$ if and only if $v$ is the solution of the system
\begin{equation} 
\label{eq:v-bis}
\mathscr{L} v^j = -U^j v + 2\sum_{i=1}^{3} V_i(Z^j)V_i(v^j) - \sum_{i=1}^3 v^i e^{Z^i} \big( V_i v^j + v^j V_i Z^j \big). 
\end{equation}
To treat the nonlinearity, we need to introduce another a priori given element in the enhancement of the noise $\zeta$. Define a $3\times 3$ matrix $\Theta$ setting formally
\begin{equation}
\label{EqDefnTheta}
\Theta^{ij} = \Pi^{(b)} \big( Z_1^i, V_i Z_1^j\big).
\end{equation}

\ssk

\begin{defn}
Given $\frac{1}{3}<\beta<\alpha<\frac{1}{2}$ and a time-independent distribution $\zeta\in C^{\alpha-2}$, a \textbf{(3d Burgers)-enhancement of $\zeta$} is a tuple $\widehat{\zeta} := \Big(\zeta, Y_2, Y_3, W_1,W_2, (W_2^j)_j,
\Theta\Big)$, with $Y_k\in L^\infty_T C^{\alpha-(5-k)/2}$, $ W_1,W_2, W_k^j\in L^\infty_T C^{2\alpha-1}$ and $\Theta\in L^\infty_T C^{2\alpha-1}$.
\end{defn} 
 
\ssk

So the \textbf{space of enhanced distributions} $\widehat{\zeta}$ for the multiplicative Burgers system is the product space 
$$
C^{\alpha-2}_{p_a}\times \prod_{k=2}^3 L^\infty_T C^{\alpha-(5-k)/2}_{p_a} \times \left(L^\infty_T C^{2\alpha-1}_{p_a}\right)^{5}   \times L^\infty_T C^{2\alpha-1}; 
$$ 
we slightly abuse notations here as the first factors in the above product refer to $\RR^3$-valued distributions/functions, while the last factr has its values in $\RR^9$. Given such an enhanced distribution $\widehat{\zeta}$, we define the Banach solution space $\mcS_{\alpha,\beta}\big(\widehat{\zeta}\,\big)$ as in Section \ref{Subsubsectionpara}, replacing the weight $p_a$ by the constant $1$. Recall the constant $\kappa>1$ appears in the time weight \eqref{eq:varpi2}.

\medskip

\begin{thm} \label{thm:pointfixe2} {\sf 
Let us work under assumption {\bf (A)}, and let $\frac{1}{3}<\alpha<\frac{1}{2}$ be given. Choose $\beta<\alpha$, the positive parameter $a$ in the weight $p_a$, and  $\epsilon>0$, such that 
$$ 
2\alpha+ \beta>1 \qquad \textrm{and} \qquad 6\epsilon\leq \alpha-\beta.
$$ 
Given an enhanced distribution $\widehat{\zeta}$ and $\widehat{v}\in\mcS_{\alpha,\beta}\big(\widehat{\zeta}\,\big)$, the nonlinear term 
$$ 
[N(v)]^j := \sum_{i=1}^3 v^i e^{Z_i} \big( \partial_i v^j + v^j V_i Z^j \big)
$$ 
is well-defined and there exists some $z^\sharp\in \calC^{1+\alpha+\beta}_\varpi$ with 
$$
\big(\mathscr{L}^{-1} [N(v)],\dots;z^\sharp\big)\in\mcS_{\alpha,\beta}\big(\widehat{\zeta}\,\big)
$$  
and
\begin{equation}  
\label{eq:fin-bis}
\Big\| \big(\mathscr{L}^{-1} [N(v)],\dots;z^\sharp\big)\Big\|_{\alpha,\beta} \lesssim \kappa^{-\epsilon} \, \Big\|\big(v,v_1,..,v_{3},v^\sharp\big) \Big\|_{\alpha,\beta}. 
\end{equation}   } 
\end{thm}

\medskip

\begin{Dem} 
We fix a coordinate $j=1,2,3$ and have to study
$$ [N(v)]^j = \sum_{i=1}^3 v^i e^{Z^i} \big( V_i v^j + v^j V_i Z^j \big).$$

The first quantity is sufficiently regular by itself, and we have $Z_i \in \calC^\alpha$, $v\in \calC^{1+\alpha}_{\varpi}$ so for every $i=1,2,3$ then
 $$ v^i e^{Z^i} V_i v^j  \in \calC^{\alpha}_{\varpi}$$
 hence
 $$ \mathscr{L}^{-1}\big[ v^i e^{Z^i} V_i v^j \big] \in \calC^{1+\alpha+\beta}_{\varpi}$$ with an acceptable norm (controlled by $\kappa^{-\epsilon}$).

Let us now focus on the second part $v^i e^{Z^i} v^j V_i Z^j$.
Since $v\in \calC^{1+\alpha}_{\varpi}$, it is very regular and the problem only relies on defining the product $e^{Z_i} V_i Z^j$.
We first decompose using paraproducts
$$ e^{Z^i} V_i Z^j = \Pi^{(b)}_{e^{Z^i}}(V_i Z^j) + \Pi^{(b)}_{V_i Z^j}(e^{Z^i}) + \Pi^{(b)}(e^{Z^i}, V_i Z^j).$$
The second term $B^{ij}$ is bounded in $\calC^{2\alpha-1}$.
The last resonant part is studied through a paralinearization formula (see \cite{BB} and references there for example)
$$ e^{Z^i} = \Pi^{(b)}_{e^{Z^i}}(Z^i) + \rest{2\alpha}$$
which implies with $\alpha>1/3$
\begin{align*}
A^{ij}:=\Pi^{(b)}(e^{Z^i}, V_i Z^j) & = \Pi^{(b)}\Big( \Pi^{(b)}_{e^{Z^i}}(Z^i) , V_i Z^j\Big) + \rest{3\alpha-1} \\
& = e^{Z^i}\Pi^{(b)}( Z^i , V_i Z^j\Big) + \rest{3\alpha-1} \\
& = e^{Z^i}\Pi^{(b)}( Z^i_1 , V_i Z^j_1\Big) + \rest{3\alpha-1} = e^{Z^i} \Theta^{ij} + \rest{3\alpha-1} ,
\end{align*} 
where we have used the commutator estimates. Since we assume that $\Theta$ is supposed to be well-defined $L^\infty_T C^{2\alpha-1}$, we conclude to $A^{ij}\in L^\infty_T C^{2\alpha-1}$.
So we observe that
$$ \Pi^{(b)}_{A^{ij}}(v^i v^j) + \Pi^{(b)}(A^{ij},v^i v^j)$$
is well-defined in $\calC^{3\alpha}_{\varpi}$ whose evaluation through $\mathscr{L}^{-1}$ is then bounded in $\calC^{1+\alpha+\beta}_{\varpi}$ with acceptable bounds. 
And since
$$ \mathscr{L}^{-1} \Pi^{(b)}_{v^i v^j}[A^{ij}] = \widetilde \Pi^{(b)}_{v^i v^j }(\mathscr{L}^{-1} A^{ij})$$
this is also controlled in $\calC^{1+\alpha+\beta}_{\varpi}$ by Schauder estimates and we conclude to
$$ \Big\| \mathscr{L}^{-1} \big(v^i v^j A^{ij}\big) \Big\|_{\calC^{1+\alpha+\beta}_\varpi} \lesssim \kappa^{-\epsilon} \|v\|_{\calC^{+\alpha}_{\varpi}}.$$

It remains the quantity with $B^{ij}$ (instead of $A^{ij}$). Here we only know that $B^{ij}$ belongs to $\calC^{2\alpha-1}$ (and not $L^\infty_T C^{2\alpha-1}$ as for $A^{ij}$) but we can take advantage of the fact that $B^{ij}$ is a paraproduct. Indeed
as before we have
$$ \Pi^{(b)}_{B^{ij}}(v^i v^j) + \Pi^{(b)}(B^{ij},v^i v^j)$$
 well-controlled in $\calC^{3\alpha}_{\varpi}$ and
 $$ \mathscr{L}^{-1} \Pi^{(b)}_{v^i v^j}[B^{ij}] = \widetilde \Pi^{(b)}_{v^i v^j }(\mathscr{L}^{-1} B^{ij}) = \widetilde \Pi^{(b)}_{v^i v^j } \big(  \widetilde  \Pi^{(b)}_{V_i Z^j}(\mathscr{L}^{-1} e^{Z^i}) \big)$$
which is well-controlled in $\calC^{1+\alpha+\beta}$ due to Schauder estimates, Proposition \ref{prop:schauder}. In conclusion, we have obtained that
\begin{align*} 
\mathscr{L}^{-1}  [N(v)]^j & = \mathscr{L}^{-1} \left[\sum_{i=1}^3 v^i v^j \Pi^{(b)}_{e^{Z^i}}(V_i Z^j) \right] + \rest{1+\alpha+\beta} \\
 & = \sum_{i=1}^3 \mathscr{L}^{-1} \left[\Pi^{(b)}_{v^i v^j}\big(\Pi^{(b)}_{e^{Z^i}}(V_i Z^j) \big) \right] + \rest{1+\alpha+\beta}  \\
 & = \sum_{i=1}^3 \widetilde \Pi^{(b)}_{v^i v^j}\widetilde \Pi^{(b)}_{e^{Z^i}} (\mathscr{L}^{-1} V_i Z^j)  + \rest{1+\alpha+\beta} \\
 & = \sum_{i=1}^3 \widetilde \Pi^{(b)}_{v^i v^j e^{Z^i}} (\mathscr{L}^{-1} V_i Z^j)  + \rest{1+\alpha+\beta},
\end{align*} 
which exactly shows that $\mathscr{L}^{-1} [N(v)]^j$ is paracontrolled by the collection 
$$
\big(\mathscr{L}^{-1} V_i Z^j\big)_{1\leq i\leq \ell_0}.
$$
\end{Dem}

\ssk

\begin{cor}
Under the assumptions of Theorem \ref{thm:pointfixe2} on the positive parameters $\alpha,\beta,a,\epsilon$, and given $\widehat{u}\in\mcS_{\alpha,\beta}\big(\widehat{\zeta}\,\big)$ with $u\in\calC_{\varpi}^{\alpha-2a-2\epsilon}$, set $v := \mathscr{L}^{-1}\big(\widehat{u}\,\widehat{\zeta} - (u\cdot V) u\big)$. Then the tuple $\big(v,u,u_1,u_2\big)$ satisfies the structure equation \eqref{ansatz}, with 
\begin{equation} 
\label{eq:finn}
\Big\| \big(v,u,u_1,u_2\big) \Big\|_{\alpha,\beta} \lesssim \kappa^{-\epsilon} \, \Big\|\big(u,u_1,u_2,u_3\big) \Big\|_{\alpha,\beta}, 
\end{equation} 
where $\kappa$ is the constant appearing in the definition \eqref{eq:varpi2} of the weight $\varpi$.
\end{cor}

We summarise in the following assumption the work about renormalisation of the ill-defined terms defining the Burgers-enhancement of $\xi$.   \vspace{0.15cm}

\noindent \textbf{\textsf{Assumption (B')}} {\sf Assumption {\bf (B)} hold and denoting by $\Theta^\epsilon$ the quantity obtained by replacing $\xi$ by $\xi^\epsilon$ in the definition \eqref{EqDefnTheta} of $\Theta$, then for any $\alpha < 1/2$, there exists some deterministic $3\times 3$ matrix-valued function $d^\epsilon$ such that $\Theta^\epsilon-d^\epsilon$ converges in probability in $L^\infty_T C^{2\alpha-1}$.   }

\medskip

This assumption is the final ingredient in the proof of Theorem \ref{thm:burgers}.
 
\medskip

\begin{Dem}[of Theorem \ref{thm:burgers}]
Well-posedness of Burgers system follows as a direct consequence of Theorem \ref{thm:pointfixe2}. Theorem \ref{thm:burgers} on the convergence of the solutions to a renormalised $\epsilon$-dependent equation to the solution of the Burgers equation is thus obtained as a direct consequence of this well-posedness result together with an additional renormalisation step that will be done in a forthcoming work. The $3\times 3$ matrix-valued functions $d^\epsilon$ is the one renormalizing the quantities $(\Theta^{i,j})^\epsilon_{1\leq i,j\leq 3}$. By tracking the changes (in the proof of Theorem \ref{thm:pointfixe2}), induced by a renormalisation of $\Theta^\epsilon$ into $\Theta^\epsilon-d^\epsilon$ in $L^\infty_T C^{2\alpha-1}$, we see that if $\big(u^\epsilon,u_1^\epsilon,u_2^\epsilon,u_3^\epsilon\big)$ satisfies Ansatz \eqref{ansatz} with $Z_i^\epsilon$, and setting $ v^\epsilon := \mathscr{L}^{-1}\big( (u^\epsilon \cdot V)u^\epsilon - d^\epsilon(u_{1}^\epsilon,u_1^\epsilon)\big)$, the tuple $\big(v^\epsilon,u^\epsilon,u^\epsilon_1,u^\epsilon_2\big)$ still satisfies the ansatz. We then complete the proof of Theorem \ref{thm:burgers}, as done for Theorem \ref{thm:pam}.
\end{Dem}

\bigskip

\section[\hspace{0.7cm} A glimpse at renormalisation matters]{A glimpse at renormalisation matters}
\label{SectionRenormalisation}

We provide in this section a flavour of the problems that are involved in proving that assumption {\bf (B)}, formulated in Section \ref{SubsectionStochasticPAM},  hold true. The analysis of the $2$ and $3$-linear terms is essentially complete, while the analysis of the $4$-linear terms is only sketched. Hairer uncovered in \cite{H} the rich algebraic setting in which renormalisation takes place within his theory of regularity structures. The full treatment of this problem was given very recently in the works \cite{BHZ} and \cite{CH} of Hairer and co-authors. They provide in particular a clear understanding of which counterterms need to/can be added in the dynamics driven by a regularized noise to get a converging limit when the regularizing parameter tends to $0$. A similar systematic treatment of renormalisation matters within the setting of high order paracontrolled calculus \cite{BB16} should be developed in a near future. We describe in this last section how things can be understood from a pedestrian point of view on the example of the (PAM) equation. \textit{We assume here as in} Theorem \ref{thm:pam} \textit{that the vector fields $V_i$ are divergence-free}; this specific assumption is used to see that the terms $W_2^j$ do not need to be renormalised.

\medskip

Basic renormalisation consists in removing from diverging random terms their expectation. While this operation is sufficient in a number of cases, such as the $2$ and $3$-dimensional (PAM) equations, or the $1$-dimensional stochastic heat equation \cite{H,HLR3}, more elaborate renormalisation procedures are needed in other examples, such as the (KPZ) or $\Phi^4_3$ equations. Hopefully, the kind of renormalisation needed here for the study of the $3$-dimensional (PAM) and Burgers equations, is essentially basic, in accordance with the work of Hairer and Labb\'e \cite{HLR3} on the (PAM) equation in $\RR^3$. 

\medskip

The a priori ill-defined terms are $2$-linear with respect to the noise
$$
\Pi^{(b)}\big(V_i Z_1, V_i Z_1\big), \quad \textrm{and }\quad \Pi^{(b)}\Big(V_i(Z_1)\,, \,V_i\mathscr{L}^{-1}(V_j Z_1)\Big),
$$
$3$-linear
$$
\Pi^{(b)}\big(V_i Z_1, V_i Z_2\big),
$$
and $4$-linear
$$
\Pi^{(b)}\big(V_i Z_1, V_i Z_3\big), \quad \textrm{and }\quad \Pi^{(b)}\big(V_i Z_2, V_i Z_2\big).
$$

\bigskip

\subsection[\hspace{-0.5cm} Renormalising the quadratic terms]{Renormalising the quadratic terms}

One takes advantage in the analysis of the renormalisation of the quadratic terms of the fact proved along the proof of Proposition \ref{prop:modifiedparaproduit} that the operator $t^{-1}\mcQ^2_t \mathscr{L}^{-1}$ is also a Gaussian operator with cancellation, an element of ${\sf GC}^{\frac{b}{8}-2}$ actually. More generally, the operators $\mcQ_t\circ V_i\circ \mathscr{L}^{-1}$ are of the form $\sqrt{t}\,\mcQ'_t$, for some Gaussian operator $\mcQ'_t$ with cancellation. Thus the term $\Pi^{(b)}\big(V_i Z_1, V_i Z_1\big)$ has the same structure as
$$
{\sf I}_2 := \int_0^1\mcP_t\Big(\mcQ_t^1 \zeta\cdot\mcQ^2_t\zeta\Big)\,dt;
$$
so does the resonent term $\Pi^{(b)}(\zeta,Z_1)$ analysed in \cite{BB} in the study of the $2$-dimensional (PAM) equation. We estimate the size of $\mcQ_r({\sf I}_2)$ in terms of $r$, to see whether or not it belongs to some parabolic H\"older space. For a space white noise $\zeta$, the expectation $\EE\left[\big|\mcQ_r({\sf I}_2)(e)\big|^2\right]$ is given by the integral on $M^2\times [0,1]^2$ of
\begin{align}
K_{\mcQ_r \mcP_{t_1}}(e,e') K_{\mcQ_r \mcP_{t_2}}(e,e'')\,\EE\Big[\mcQ^1_{t_1} \xi(e') \mcQ^2_{t_1} \xi(e') \mcQ^1_{t_2} \xi(e'') \mcQ^2_{t_2} \xi(e'')\Big]  
\label{eq:compte} \end{align}
against the measure $\nu(de')\nu(de'') dt_1dt_2$. The expectation in \eqref{eq:compte} is estimated with Wick's formula by
\begin{align*} 
\EE\big[\mcQ^1_{t_1} \xi(e') \mcQ^2_{t_1} \xi(e')\big] &\EE\big[\mcQ^1_{t_2} \xi(e'') \mcQ^2_{t_2} \xi(e'')\big] + \EE\big[\mcQ^1_{t_1} \xi(e') \mcQ^1_{t_2} \xi(e'')\big] \EE\big[ \mcQ^2_{t_2} \xi(e'') \mcQ^2_{t_1} \xi(e')\big] \\
&+ \EE\big[\mcQ^1_{t_1} \xi(e') \mcQ^2_{t_2} \xi(e'')\big] \EE\big[ \mcQ^1_{t_2} \xi(e'') \mcQ^2_{t_1} \xi(e')\big] \\
&\lesssim (t_1t_2)^{-d/2}  + \G_{t_1+t_2}(e',e'')^2,
\end{align*} 
where $d$ is the homogeneous dimension of the ambiant space $M$; the term $(t_1t_2)^{-d/2}$ comes from the first product of expectations. The quantity $\EE\left[\big|\mcQ_r\big({\sf I}_2\big)\big|^2\right]$ can thus be bounded above by the sum of two integrals, with $dm := \nu(de')\nu(de'') dt_1dt_2$ and a first integral equal to
\begin{align*}
&\int K_{\mcQ_r \mcP_{t_1}}(e,e') K_{\mcQ_r \mcP_{t_2}}(e,e'')  \G_{t_1+t_2}(e',e'')^2 dm   \\
& \quad \lesssim 
\int \left(\frac{r}{r+t_1}\right)^a \left(\frac{r}{r+t_2}\right)^a (t_1+t_2)^{-\frac{\nu}{2}} \int \G_{r+t_1}(x,y) \G_{r+t_2}(x,z)\G_{t_1+t_2}(y,z) \, dm    \\
& \quad \lesssim \int \left(\frac{r}{r+t_1}\right)^N \left(\frac{r}{r+t_2}\right)^a (t_1+t_2)^{-\frac{d}{2}} (r+t_1+t_2)^{-d/2} dt_1dt_2   \\
& \quad \lesssim r^{2-d}
\end{align*}
for $d<4$; we used the upper bound \eqref{EqIteratedG} here. We also have
\begin{align*}
&\int \int K_{\mcQ_r \mcP_{t_1}}(e,e') K_{\mcQ_r \mcP_{t_2}}(e,e'')  (t_1t_2)^{-d/2} d\nu(de')\nu(de'') dt_1dt_2   \\
& \qquad \lesssim \int \left(\frac{r}{r+t_1}\right)^a \left(\frac{r}{r+t_2}\right)^a (t_1t_2)^{-\frac{d}{2}} dt_1dt_2 
\end{align*}
for some positive exponent $a$, with a relatively sharp upper bound, which happens to be infinite in dimension $2$ or larger. This is the annoying bit. Considering ${\sf I}_2-\EE\big[{\sf I}_2\big]$ instead of ${\sf I}_2$ removes precisely this diverging part in the corresponding Wick formula for $\EE\left[\big|\mcQ_r\big({\sf I}_2-\EE[{\sf I}_2]\big)\big|^2\right]$. It follows as a consequence that one has
$$
\EE\Big[\big|\mcQ_r\big({\sf I}_2-\EE[{\sf I}_2]\big)\big|^2\Big]^{\frac{1}{2}}\lesssim r^{1-\frac{d}{2}},
$$ 
which shows that the associated distribution is almost surely in $\calC^{(2-d)^-}$, by Kolmogorov's continuity criterion. 

\medskip

While the above reasoning shows that recentering 
$$ 
W_2^{j,\epsilon}:=\sum_{i=1}^{\ell_0} \Pi^{(b)}\Big(V_i(Z_1^\epsilon), \big[V_i\mathscr{L}^{-1}(V_j Z_1^\epsilon)\big] \Big)
$$
around its expectation makes it converge in the right space, \textit{there is actually no need to renormalize this term}, as can be expected from comparing our setting with the setting of regularity structures for the 3-dimensional setting, investigated in Hairer and Labb\'e's work \cite{HLR3}. 

\ssk

One can see that point by proceeding as follows. Replace in a first step the study of the above quantity by a similar quantity where the spacetime paraproduct $\Pi^{(b)}$ and resonent term $\Pi^{(b)}(\cdot,\cdot)$ are replaced by a space paraproduct $\pi^{(b)}$ and resonent operator $\pi^{(b)}(\cdot,\cdot)$ introduced and studied in \cite{BB} -- they are defined in the exact same way as $\Pi^b$, but without the time convolution operation. Continuity properties were proved for such spatial paraproduct in \cite{BB}, and we shall use in addition an elementary comparison result between this spatial paraproduct and our space-time paraproduct proved by Gubinelli, Imkeller and Perkowski in their setting \cite[Lemma 5.1]{GIP}. A similar statement and proof holds with the two paraproducts $\Pi^{(b)}$ and $\pi^{(b)}$; we state it here for convenience.

\ssk
  
\begin{lem*} Let $\omega_1,\omega_2$ be two space-time weights. If $u\in \calC^\alpha_{\omega_1}$ for $\alpha\in(0,1)$ and $v\in L^\infty_T C^\beta_{\omega_2}$ for some $\beta\in(-3,3)$ then
$$ \pi^{(b)}_u(v) - \Pi^{(b)}_u(v) \in L^\infty_T C^{\alpha+\beta}_{\omega}$$
with $\omega=\omega_1 \omega_2$. 
\end{lem*}

\ssk
 
Setting
$$ 
w_2^{j,\epsilon}:=\sum_{i=1}^{\ell_0} \pi^{(b)}\Big(V_i(Z_1^\epsilon), \big[V_i\mathscr{L}^{-1}(V_j Z_1^\epsilon)\big] \Big)
$$
and using the comparison lemma and then the continuity estimates of each paraproduct, we see that $W_2^{j,\epsilon}-w_2^{j,\epsilon}$ is equal to
{\small \begin{align*}
\sum_{i=1}^{\ell_0} \Pi^{(b)}_{V_i(Z_1^\epsilon)}\big(V_i\mathscr{L}^{-1}(V_j Z_1^\epsilon)\big) &- \pi^{(b)}_{V_i(Z_1^\epsilon)}\big(V_i\mathscr{L}^{-1}(V_j Z_1^\epsilon)\big) + \Pi^{(b)}_{V_i\mathscr{L}^{-1}(V_j Z_1^\epsilon)}\big(V_i(Z_1^\epsilon)\big)   \\
&- \pi^{(b)}_{V_i\mathscr{L}^{-1}(V_j Z_1^\epsilon)}\big(V_i(Z_1^\epsilon)\big)    \\
 & \hspace{-1cm}\in \sum_{i=1}^{\ell_0} \left[\Pi^{(b)}_{V_i(Z_1^\epsilon)}\big(V_i\mathscr{L}^{-1}(V_j Z_1^\epsilon)\big) -\pi^{(b)}_{V_i(Z_1^\epsilon)}\big(V_i\mathscr{L}^{-1}(V_j Z_1^ \epsilon)\big)\right] + L^\infty_T C^{2\alpha-1},
 \end{align*}   } 
so it is an element of $\calC^{2\alpha-1}$. So in order to estimate $W_2^{j,\epsilon}$ is the suitable H\"older space we only need to study its "spatial" counterpart $w_2^{j,\epsilon}$. This can be done as follows.

\medskip

As $W_2^{j,\epsilon}$, the quantity $w_2^{j,\epsilon}$ is quadratic as a function of the noise, however we are going to see that its expectation is already bounded in $\calC^{2\alpha-1}$, as a consequence of some symmetry properties -- this explains why $w_2^{j,\epsilon}$ is directly converging in $\calC^{2\alpha-1}$, with no renormalisation needed along the way. The term $w_2^{j,\epsilon}$ can indeed be written as a finite sum of integrals in time of terms of the form
{\small $$
P_t \left[ Q^1_t V_i \mathscr{L}^{-1} \xi ^\epsilon \cdot  Q^2_t V_i\mathscr{L}^{-1}(V_j \mathscr{L}^{-1} \xi^\epsilon) \right](e) + P_t \left[ Q^2_t V_i \mathscr{L}^{-1} \xi ^\epsilon \cdot  Q^1_t V_i\mathscr{L}^{-1}(V_j \mathscr{L}^{-1} \xi^\epsilon) \right](e),
$$}
where the localizing operators $P_t$ and $Q_t$ are only in space. Using the above additional geometric assumptions on the operator, the previous integral can be estimated, up to a satisfying remainder term controlled in terms of $t^{2\alpha}$, by
{\small $$ 
P_t \left[ Q^1_t V_i \mathscr{L}^{-1} \xi ^\epsilon \cdot  V_j Q^2_t \mathscr{L}^{-1}(V_i \mathscr{L}^{-1} \xi^\epsilon) \right](e) + P_t \left[ Q^2_t V_i \mathscr{L}^{-1} \xi ^\epsilon \cdot  V_j Q^1_t \mathscr{L}^{-1}(V_i \mathscr{L}^{-1} \xi^\epsilon) \right](e).
$$}
Its expectation can be seen to converge in $\calC^{2\alpha-1}$ to
{\small $$ 
\int K_{P_t}(x,y) \left[ K_{ [V_j Q^2_t \mathscr{L}^{-1}(V_j \mathscr{L}^{-1})]^* Q^1_t V_i \mathscr{L}^{-1} }(y,y) + K_{ [V_j Q^1_t \mathscr{L}^{-1}(V_j \mathscr{L}^{-1})]^* Q^1_t V_i \mathscr{L}^{-1} }(y,y)\right] \mu(dy),
$$}
where $^*$ denotes the usual adjoint in  $L^2(M,d\mu)$ (in space) and where the time is fixed in the operator $\mathscr{L}^{-1}$. By symmetry, it is equal to
$$ 
\int K_{P_t}(x,y) \left[ K_{ {\mathscr{L}^{-1}}^* V_i {\mathscr{L}^{-1}}^*  Q_t V_i \mathscr{L}^{-1} }(y,y) \right] \mu(dy)
$$
where $Q_t:=Q^{2,*}_t V_j Q^1_t + Q^{1,*}_t V_j Q^2_t$ is antisymmetric. Since at time fixed, the spatial operator ${\mathscr{L}^{-1}}^*$ is self-adjoint, we deduce that ${\mathscr{L}^{-1}}^* V_i {\mathscr{L}^{-1}}^*  Q_t V_i \mathscr{L}^{-1}$ is antisymmetric in space and so its kernel is vanishing on the diagonal. This shows as a consequence that $ \EE\big[ w_2^{j,\epsilon} \big]$ is bounded in the parabolic H\"older space $\calC^{2\alpha-1}$.

 \bigskip

\subsection[\hspace{-0.5cm} Higher order terms]{Higher order terms}

The analysis of the $3$-linear term $\Pi^{(b)}\big(V_iZ_1,V_iZ_2\big)$ can be done exactly as for the $2$-linear term $\Pi^{(b)}\big(V_iZ_1,V_iZ_1\big)$, starting from the fact that the former has the same structure as 
$$
{\sf I}_3 := \int_0^1 \mcP_{t_1}\Big(\mcQ^1_{t_1}\zeta\cdot\mcQ_{t_1}^2\big(\mcP_{t_2}\big\{\mcQ_{t_2}^3\zeta\cdot\mcQ_{t_2}^4\zeta\big\}\big)\Big) \,dt_2dt_1,
$$ 
where the $\mcQ^i_t$ are Gaussian operators with cancellation. Its renormalised version ${\sf I}_3^r$ is defined by replacing $\mcQ_{t_2}^3\zeta\cdot\mcQ_{t_2}^4\zeta$ by $\mcQ_{t_2}^3\zeta\cdot\mcQ_{t_2}^4\zeta - \EE\big[\mcQ_{t_2}^3\zeta\cdot\mcQ_{t_2}^4\zeta\big]$. The quantity $\EE\Big[\big|\mcQ_r({\sf I}_3)(e)\big|^2\Big]$ is thus given by an integral with respect to some kernels with some Gaussian controls and cancellation property of the expectation of a product of six Gaussian random variables indexed by parabolic points $(e',e'',\bar{e}', \bar{e}'')$
$$
\big(\mcQ^1_{t_1}\zeta\big)(e')\,\big(\mcQ^1_{s_1}\zeta\big)(e'')\,\big(\mcQ^3_{t_2}\zeta\big)(\bar{e}')\,\big(\mcQ^4_{t_2}\zeta\big)(\bar{e}')\,\big(\mcQ^3_{s_2}\zeta\big)(\bar{e}'')\,\big(\mcQ^4_{s_2}\zeta\big)(\bar{e}'').
$$
In its renormalised version, the above product $\big(\mcQ^3_{t_2}\zeta\big)(\bar{e}')\,\big(\mcQ^4_{t_2}\zeta\big)(\bar{e}')$ is replaced by 
$$
\big(\mcQ^3_{t_2}\zeta\big)(\bar{e}')\,\big(\mcQ^4_{t_2}\zeta\big)(\bar{e}') - \EE\big[\big(\mcQ^3_{t_2}\zeta\big)(\bar{e}')\,\big(\mcQ^4_{t_2}\zeta\big)(\bar{e}')\big],
$$ 
and similarly for $\big(\mcQ^3_{s_2}\zeta\big)(\bar{e}'')\,\big(\mcQ^4_{s_2}\zeta\big)(\bar{e}'')$. Using Wick formula then shows that 
$$
\EE\Big[\big|\mcQ_r({\sf I}^r_3)(e)\big|^2\Big]
$$ 
only involves products of expectations where no two identical parameters $s,t$ appear inside each expectation, meaning that we have, after some elementary computations, an estimate of the form
\begin{align*}
\EE\left[| \mcQ_r\big({\sf I}^r_3\big)(e)|^2\right]  \lesssim \int K_{\mcQ_r \mcP_{t_1}}(e,e') K_{\mcQ_r \mcP_{t_2}}(e,e'') (t_1t_2)^{-d/2} \G_{t_1+t_2}(e',e'') t_1t_2\,dm,
\end{align*}
with $dm = \nu(de') \nu(de'')dt_1dt_2$, as above. For $d<4$, this gives the estimate
\begin{align*}
\EE\left[| \mcQ_r\big({\sf I}^r_3\big)(e)|^2\right]  & \lesssim  \int\int \left(\frac{r}{r+t_1}\right)^a \left(\frac{r}{r+t_2}\right)^a (t_1t_2)^{-d/2} (r+t_1+t_2)^{-d/2}  t_1t_2 dt_1dt_2 \\
& \lesssim r^{-3d/2+4},
\end{align*}
on which one reads that ${\sf I}^r_3$ has almost surely regularity $(-3d/2+4)^-$, so no renormalisation is required. This is coherent with \cite{HLR3} - equation (5.3), where it is shown, within the setting of regularity structures, that the  terms that are trilinear functions of the noise do not need to be renormalised. Everything happens here as if we were working with a $3$-linear term of the form 
$$
\int_0^1 \mcP_t\big(\mcQ_t^1\zeta\cdot\mcQ_t^2\zeta\cdot t\mcQ_t^3\zeta\big)\,dt;
$$
one can indeed make that comparison concrete.

\bigskip

$\bullet$ The model quantities corresponding to the $4$-linear terms are of the type
$$
{\sf I}_4 := \int_0^1 \mcP_t \big( \mcQ^1_t \zeta \cdot \mcQ^2_t \zeta \cdot t \mcQ^3_t \zeta \cdot t \mcQ^4_t \zeta \big)\,dt.
$$
One can see on such terms that a basic renormalisation procedure suffices to get objects of regularity $0^-$, in dimension $3$, such as expected.
This finishes the sketch of proof that assumption {\bf (B)} actually holds true.

\bigskip
\bigskip
\bigskip
 
\noindent {\bf Acknowledgements.} We thank the two referees for their constructive critical remarks on previous versions of this work. 

\bigskip
\bigskip
\bigskip


\end{document}